\documentclass[a4paper]{amsart} 

\newcommand{\f}{\frac}

\newcommand{\im}{\operatorname{Im}}

\newcommand{\re}{\operatorname{Re}}

\newcommand{\R}{\mathbb R}
\newcommand{\C}{\mathbb C}
\newcommand{\N}{\mathbb N}
\newcommand{\Z}{\mathbb Z}
\newcommand{\eps}{\varepsilon}
\renewcommand{\epsilon}{\varepsilon}
\newcommand{\Om}{\Omega}

\newcommand{\dist}{\operatorname{dist}}

\newcommand{\SCI}{\operatorname{SCI}}

\newcommand{\I}{\mathrm{i}}
\newcommand{\e}{\mathrm{e}}

\newcommand{\cL}{\mathcal{L}}
\newcommand{\cH}{\mathcal{H}}

\newcommand{\Q}{\mathbb{Q}}
\newcommand{\spann}{\operatorname{span}}
\newcommand{\dhaus}[2]{d_{\mathrm{H}}\!\left(#1,#2\right)}
\newcommand{\sigmaess}{\sigma_{\mathrm{ess}}}
\newcommand\dom{\operatorname{dom}}

\newenvironment{enumi}{\begin{enumerate}[(i)]}{\end{enumerate}}

\usepackage{etex}
\usepackage[english]{babel}
\usepackage{tikz}
\usetikzlibrary{arrows,decorations.pathmorphing,backgrounds,positioning,fit,petri,patterns}
\usepackage{pgfplots}
 \pgfplotsset{every tick label/.append style={font=\footnotesize},
}
\usepackage[utf8]{inputenc}
\usepackage{graphicx}
\usepackage[colorlinks=true,linkcolor=red,citecolor=blue]{hyperref}
\usepackage{amsfonts}
\usepackage{enumerate}
\usepackage{booktabs}
\usepackage{array} 
\usepackage{paralist} 
\usepackage{subfig} 
\usepackage{mathtools}
\usepackage{tabu}
\usepackage{amsthm}
\usepackage{empheq}
\usepackage{amsopn}
\usepackage{dsfont}
\usepackage{wrapfig}
\usepackage[font=small]{caption}
\usepackage{units}
\usepackage[symbol]{footmisc}
\usepackage{amssymb}
\usepackage{pdfpages}
\usepackage{setspace}
\usepackage{float}
\usepackage{soul}
\usepackage{multirow}
\usepackage{scrextend}
\usepackage[outline]{contour}
\usepackage[ruled, vlined]{algorithm2e}
\usepackage{soul}
\usepackage{tabularx}
\usepackage{mathabx}

\allowdisplaybreaks
\numberwithin{equation}{section}

\theoremstyle{definition}
\newtheorem{de}{Definition}[section]
\newtheorem{hyp}[de]{Hypothesis}
\newtheorem{notation}[de]{Notation}
\theoremstyle{plain}
\newtheorem{prop}[de]{Proposition}
\newtheorem{lemma}[de]{Lemma}
\newtheorem{theorem}[de]{Theorem}

\newtheorem{corollary}[de]{Corollary}
\newtheorem{cor}[de]{Corollary}

\theoremstyle{remark}
\newtheorem{remark}[de]{Remark}
\newtheorem{example}[de]{Example}

\newcommand{\cM}{\mathcal{M}}

\pgfplotsset{every tick label/.append style={font=\scriptsize}}
\pgfplotsset{compat = 1.3}

\newlength\figurewidth
\newlength\figureheight
\newlength\hfiguremargin
\newlength\vfiguremargin
\newlength{\tablength}
\newcounter{sectiontracker}
\newcounter{deftracker}

\title{Computing Klein-Gordon Spectra}
\author{Frank R\"{o}sler}
\address{Mathematisches Institut, Universit\"{a}t Bern, Alpeneggstrasse\ 22, 3012, Bern, Switzerland}
\email{frank.roesler@unibe.ch}
\author{Christiane Tretter}
\address{Mathematisches Institut, Universit\"{a}t Bern, Sidlerstr.\ 5, 3012, Bern, Switzerland} 
\email{tretter@math.unibe.ch}
\thanks{FR acknowledges support from the European Union's Horizon 2020 Research and Innovation Programme under the Marie Sklodowska-Curie grant agreement no.\ 885904.
CT gratefully acknowledges the support of the Swiss National Science Foundation (SNF), 
grant no.\ $200021\_204788$.}
\date\today

\begin{document}
\maketitle
\begin{abstract}
		We study the computational complexity of the eigenvalue problem for the Klein-Gordon equation in the framework of the Solvability Complexity Index Hierarchy. We prove that the eigenvalue of the Klein-Gordon equation with linearly decaying potential can be computed in a single limit with guaranteed error bounds from above. The proof is constructive, i.e.\ we obtain a numerical algorithm that can be implemented on a computer.
Moreover, we prove abstract enclosures for the point spectrum of the Klein-Gordon equation and we compare our numerical results to these enclosures.
Finally, we apply both the implemented algorithm and our abstract enclosures to several physically relevant potentials such as Sauter and cusp potentials
and we provide a convergence and error analysis.
\end{abstract}
\setcounter{tocdepth}{1}
%
%
\section{Introduction and Main Results}\label{sec:intro}
Reliably computing eigenvalues and spectra poses major challenges across physics, spanning the wide range from computing resonances in quantum mechanics over the wealth of extremely hard computational spectral problems in hydrodynamics and electromagnetics to the new era of spectral light tuning in metamaterials, see e.g.\ \cite{MR3272932}, \cite{MR3285902}, \cite{MR4155239}, \cite{MR4302658}. 
These challenges are due to inherent properties of physical systems including lack of symmetry, coupling of different phenomena or infinite dimensional effects such as non-discrete spectrum which may cause serious failures of domain truncation methods or finite dimensional approximations, see e.g.\ \cite{MR4382583}, \cite{MR2727810}. 
There are two such highly undesirable failures: spectral pollution, where finite dimensional eigenvalue approximations converge to a point that is not a true spectral point and spectral invisibility where a true spectral point is not seen by the finite dimensional approximations, see e.g.\ \cite{MR3413891}, \cite{MR4083777}, \cite{MR4209764}. Therefore there is an urgent need for tools that allow physicists and applied mathematicians to assess whether their computational problems are prone to these failures and to estimate the complexity of the required computational tasks.

A key breakthrough in this direction was the introduction of the Solvability Complexity Index (SCI) by A.\ Hansen in 2011 \cite{MR2726600} and the SCI Hierarchy by J.\ Ben-Artzi et al.\ in 2015 \cite{AHS}  which has opened up a completely new pathway to the analysis and numerics of spectral problems and sparked progress for many computational problems in physics. The SCI hierarchy offers a framework to compare the complexity of computational tasks such as approximating the spectrum in Hausdorff distance (or the Attouch-Wets metric in the unbounded case) for different classes of linear operators. This classification involves the number of successive limits required for the approximation, the availability of error bounds and, in its finest form, also the type of algorithm such as arithmetic.

The enormous impact of the SCI hierarchy in computational spectral theory has several reasons. First, it revealed why even the most common open problems such as computing the eigenvalues of Schr\"odinger operators with bounded, even real-valued potentials using point-values of the potential have remained open for more than 9 decades since quantum mechanics was created. Secondly, it helped to solve such longstanding open problems not only for compact, bounded, selfadjoint and normal operators \cite[Thm.\ 7.5]{AHS}, but also for unbounded operators with certain resolvent growth and non-empty essential spectrum. Thirdly, the corresponding results are not only theoretical, but they yield practically implementable algorithms which have been shown to perform well for large scale problems in physics~\cite{MR3980052}.

The urgent need to classify computational problems in spectral theory in the SCI Hierarchy 
is further substantiated by parallel research on challenges in spectral approximation such as spectral pollution or spectral invisibility 
and methods to avoid them, see e.g.\ \cite{MR3272932}, \cite{MR4083777}. 
In recent years, since the fundamental work \cite{AHS} which contains a detailed study on computing the spectra and pseudospectra of bounded matrix operators on $\ell^2(\N)$ 
as well as results on computing spectra and pseudospectra of Schr\"odinger operators with potentials satisfying a uniform BV bound 
(cf.\ \cite[Thm.\ 8.3, 8.5]{AHS}), the SCI theory has been further developed in multiple directions. We mention the works \cite{colbrook2019computation, Colbrook:2019aa, Webb:2021aa} on matrix operators; \cite{becker2020computing, colbrook22} on solving PDEs; \cite{colbrookhorningtownsend} on the computation of spectral measures, and \cite{ben2020computing, ben2022computing} on the computation of scattering resonances. 

While many of the above works obtained results on the (nonrelativistic) Schr\"odinger equation, 
no SCI results seem to exist so far on relativistic models, such as the Klein-Gordon or Dirac equations. The present article begins to fill this gap: we study the spectral problem for the Klein-Gordon equation in the framework of the SCI Hierarchy. The novelty of this contribution is 
computational, physically relevant as well as mathematical.
Indeed, the spectral problem for the Klein-Gordon equation is non-standard in the sense that one is lead to the study of quadratic operator pencils rather than classical eigenvalue equations which lead to linear monic pencils. 
Moreover, essential spectrum is not only present but it is not semi-bounded; even bounded symmetric potentials may create complex eigenvalues in addition to the two unbounded rays of real essential spectrum. 
Accordingly, new methods are needed to study spectral computation which require an intimate interplay of analysis and numerics with operator theory. The development and evaluation of these techniques, along with corresponding implementable algorithms, are the focus of this article.

In the next two subsections we present our main results on the computational spectral problem for the Klein-Gordon equation and we introduce the SCI hierarchy allowing us to interpret our results~therein.
\subsection{The Klein-Gordon equation}\label{sec:preliminaries_KG_intro}
In quantum mechanics the Klein-Gordon equation 
\begin{equation}
\label{KGequ} 
\left( - \Bigl( -\I \hbar\frac \partial {\partial t} - e \varphi
\Bigr)^2 + c^2 \left( -\I\hbar \nabla - \frac{e}{c}\vec A \right)^2 + m^2c^4
\right) U = 0
\end{equation}
describes the motion  of a relativistic spinless particle with mass $m$ and
charge $e$ in an electromagnetic field with scalar potential $\varphi$ and vector potential $\vec A$;
here $c$ is the speed of light and $\hbar$ is the Planck constant. 
If we separate time by setting $U(x,t) \!=:\! \e^{\I\lambda/\hbar t} u(x)$, $x\!\in\!\R^d$,~$t\!\in\!\R$, normalize $c$ to $1$, let $V$ be the multiplication operator by $e\varphi$ in $L^2(\R^d)$ and $A_0\!:=\!( -\I \hbar\nabla \!-\! e\vec A)^2$, then \eqref{KGequ} leads to a quadratic eigenvalue problem in $\lambda$, which we will cast in a more abstract~framework.

To this end, let $(\cH,\langle\cdot,\cdot\rangle)$ be a separable Hilbert space, let $A_0$ be a nonnegative operator on $\cH$ and $H_0:=A_0+m^2$ with $m>0$. If $V$ is a symmetric operator with $\dom H_0^{1/2} \subset \dom V$, the operator $VH_0^{-1/2}$ is bounded in $\cH$. Then the quadratic eigenvalue problem associated with \eqref{KGequ} is of the form $T_V(\lambda)u=0$, $\lambda\in\C$, where the 
Klein-Gordon operator polynomial (or pencil) in $\cH$ is given by
\begin{align*}
	T_V(\lambda) := H_0 - (V-\lambda)^2, \quad \dom T_V(\lambda) = \dom H_0, \qquad \lambda\in\C.
\end{align*}
If we assume that $S:=V H_0^{-1/2}=S_0+S_1$ with a strict contraction $S_0$ and compact~$S_1$, as in \cite{MR2238908}, \cite{MR2465932}, \cite{MR2268872}, then it is well-known that the essential spectrum of $T_V$ has a gap around $0$ and that the non-real spectrum of $T_V$ is discrete, see  \cite{MR2238908}, \cite{MR2465932}. Moreover, one can show that if $S_0=0$, then the essential spectrum of $T_V$ is given by the half lines $(-\infty,-m]\cup[m,\infty)$ and any other spectral points are discrete eigenvalues. We ensure that $S_0=0$ by making the following assumptions throughout the paper.
\begin{hyp}\label{hyp:assumptions_on_V}
	Unless otherwise stated, assume that $\cH = L^2(\R^d)$, and $\vec{A}=0$, i.e.\ $A_0 = -\Delta$, $\dom(A_0)=H^2(\R^d)$. Moreover, assume that
	\begin{enumi}
		\item
		$V\in W^{1,p}(\R^d)$ for some $p>d$,
		\item
		there exists a constant $M>0$ such that
		\begin{align}\label{eq:hypothesis_V_bounds}
			\|V\|_{W^{1,p}(\R^d)} &\leq M, \quad
			|V(x)| \leq M(1+|x|^2)^{-\f12} \quad \text{for all } x\in \R^d.
		\tag{$H_M$}
		\end{align}
	\end{enumi}
\end{hyp}
It is easy to see from the Fr\'echet-Kolmogorov-Riesz theorem \cite[Th. XIII.66]{RS4} that Hypothesis \ref{hyp:assumptions_on_V} implies compactness of $V H_0^{-1/2}$ and thus $S_0=0$. Therefore
$$
	\sigmaess(T_V) = \sigmaess(T_0) = \{\lambda\in\C\,|\,\lambda^2\in\sigma(H_0)\} = -\sqrt{\sigma(H_0)} \cup \sqrt{\sigma(H_0)},
$$
i.e.\ $\sigmaess(T_V)$ is independent of $V$ as long as $V$ satisfies Hypothesis \ref{hyp:assumptions_on_V}.
Our main result concerns the convergence of the numerical approximation of the spectrum $\sigma(T_V)$. Here, for subsets $A,B\subset\C$, let $\dhaus{A}{B}$ denote their Hausdorff distance.
\begin{theorem}\label{th:convergence_of_algorithm_intro}
	Let $V$ satisfy Hypothesis {\rm \ref{hyp:assumptions_on_V}}. Then there exist computational routines $\Gamma_n$ (depending only on $p$, $M$), which take their input from the set of point values of $V$ and produce sets $\Gamma_n(V)\subset\C$ such that
	\begin{align}
		d_{\textnormal{H}}\big(\Gamma_n(V) ,\, \sigma(T_V)\big) &\to 0 \quad\text{ for }n\to\infty.
		\tag{i}
		\label{eq:convergence_of_alg_intro}
	\end{align}
	If the constant $M$ in \eqref{eq:hypothesis_V_bounds} is known a-priori, then in addition to \eqref{eq:convergence_of_alg_intro} we obtain the error bound 
	\begin{align}
		\sup_{z\in\sigma(T_V)} \dist(z,\Gamma_n(V)) &\leq \f1n \quad \text{ for all }n\in\N
		\tag{ii}
		\label{eq:error_bound_intro}
	\end{align}
	$($for details, see Theorem {\rm \ref{th:convergence_of_Gamma_n})}.
\end{theorem}
\begin{remark}
\begin{enumi}
	\item
	Theorem \ref{th:convergence_of_algorithm_intro} implies a classification result in the Solvability Complexity Index Hierarchy. For details, see Corollary \ref{cor:classification} below.
	\item
	Not only does Theorem \ref{th:convergence_of_algorithm_intro} give an existence result, but its proof is constructive. In Definition \ref{de:definition_of_Gamma_n} below we give a concrete definition of an algorithm that achieves \eqref{eq:convergence_of_alg_intro}, \eqref{eq:error_bound_intro} and can be implemented on a computer. Indeed, in Section \ref{sec:numerics} we present numerical results obtained from a MATLAB implementation of our algorithm.
\end{enumi}
\end{remark}
\subsection{The SCI Hierarchy}\label{sec:preliminaries_SCI}
The SCI Hierarchy, introduced by Hansen \cite{MR2726600}, is a novel and rapidly developing field. In its initial form, the Solvability Complexity Index gave a rigorous framework to study the ``number of successive limits'' needed for the numerical approximation of spectral problems. Since then, it has evolved into a general theory of computability and error estimation for arbitrary computational problems in mathematics.
The starting point of the theory is the following rigorous definition of a computational problem.
\begin{de}[Computational problem]\label{def:computational_problem}
	A \emph{computational problem} is a quadruple $(\Omega,\Lambda,\cM,\Xi)$, where 
	\begin{enumi}
		\item $\Om$ is a set, called the \emph{primary set},
		\item $\Lambda$ is a set of complex-valued functions on $\Om$, called the \emph{evaluation set},
		\item $\mathcal M$ is a metric space,
		\item $\Xi:\Om\to \mathcal M$ is a map, called the \emph{problem function}.
	\end{enumi}
\end{de}
\begin{example}\label{ex:SCI_example}
	As an example, we show how a rigorous formulation of the task \emph{``compute the spectrum of a compact operator from its matrix elements''} fits into Definition \ref{def:computational_problem}. 

	Let $\Omega=\mathcal{S}^\infty(\ell^2(\N))$ be the space of compact operators on the Hilbert space $\ell^2(\N)$. We can formulate the computation of the spectrum of an operator from $\Omega$ as a computational problem in the following way. Let $\{e_i\}_{i\in\N}$ denote the canonical basis of $\ell^2(\N)$ and define $f_{ij}:\Omega \to \C$; $f_{ij}(A) = \langle e_i, A e_j\rangle$ and $\Lambda := \{f_{ij}(\cdot)\,|\,i,j\in\N\}$. Next, let $\cM$ be the set of closed bounded subsets of $\C$, endowed with the Hausdorff distance $d_{\textnormal{H}}$ and let $\Xi:\Omega\to\cM$; $\Xi(A)=\sigma(A)$ be the function mapping a linear operator onto its spectrum. Then the quadruple $(\Omega,\Lambda,\cM,\Xi)$ is a computational problem in the sense of Definition \ref{def:computational_problem}. 
	The above example gives an intuitive meaning to the building blocks (i)-(iv):
	\begin{enumi}
		\item $\Omega$ is the set of objects that give rise to the computational problem,
		\item $\Lambda$ is the set of information an algorithm is allowed to access during the computation,
		\item the metric on $\cM$ measures the approximation error,
		\item the image of the function $\Xi$ contains the objects to be computed.
	\end{enumi}
\end{example}
Moving on from Definition \ref{def:computational_problem}, we would like to be able to say when a computational problem is considered ``solved'' and ``how good'' it has been solved (see questions (1)-(4) below). To this end, we proceed by introducing a rigorous definition of what we mean by an algorithm.
\begin{de}[General algorithm]\label{def:Algorithm}
	Let $(\Omega,\Lambda,\cM,\Xi)$ be a computational problem. A \emph{general algorithm} is a mapping $\Gamma:\Om\to\mathcal M$ such that for each $T\in\Om$ 
	\begin{enumi}
		\item there exists a finite (non-empty) subset $\Lambda_\Gamma(T)\subset\Lambda$,
		\item the action of $\Gamma$ on $T$ depends only on $\{f(T)\}_{f\in\Lambda_\Gamma(T)}$,
		\item for every $S\in\Om$ with $f(T)=f(S)$ for all $f\in\Lambda_\Gamma(T)$ one has $\Lambda_\Gamma(S)=\Lambda_\Gamma(T)$.
	\end{enumi}
	We will sometimes write $\Gamma(T) = \Gamma(\{f(T)\}_{f\in\Lambda_\Gamma(T)})$ to emphasise point (ii) above: the output $\Gamma(T)$ depends only on the (finitely many) evaluations $\{f(T)\}_{f\in\Lambda_\Gamma(T)}$.
\end{de}
Definition \ref{def:Algorithm} above is a general abstraction of what a computer algorithm does: it takes a finite amount of input (i.e.\ $\{f(T)\}_{f\in\Lambda_\Gamma(T)}$) and returns some output in a metric space (usually a string of numbers). The term \emph{general algorithm} is used, because Definition \ref{def:Algorithm} imposes no restrictions on how the output $\Gamma(T)$ is computed. We will specify more restrictive types of algorithms in the next section by introducing a \emph{recursiveness} constraint (the reader may think of the way a Turing machine produces its output).
Given the above definitions one can ask the following questions:
\begin{enumerate}
	\item Given a computational problem $(\Omega,\Lambda,\cM,\Xi)$, does there exist a sequence of algorithms $(\Gamma_n)_{n\in\N}$ such that $\Gamma_n(T)\to \Xi(T)$ in $\cM$ for all $T\in\Omega$?
	\item If so, do we have explicit error bounds, i.e.\ $d(\Gamma_n(T),\Xi(T))<\eps_n$ for a known sequence $(\eps_n)_{n\in\N}$?
	\item If $\cM$ has certain ordering properties, do we have error bounds from above / below?
	\item If all of the above fail, does it help to ``add limits''? More precisely, does there exist a family $(\Gamma_{n_k,\dots,n_2,n_1})_{n_1,\dots,n_k\in\N}$ of algorithms with $\lim_{n_k\to\infty}\cdots\lim_{n_1\to\infty}\Gamma_{n_k,\dots,n_1}(T) \!=\! \Xi(T)$ for all $T\!\in\!\Omega$?
\end{enumerate}
These questions are nontrivial: there exist examples and counterexamples to every single one of the above questions. In particular, there exist examples for question (4), i.e.\ computational problems that inherently require more than 1 limit to solve \cite[Th. 7.5]{AHS}. This shows that it is not enough to merely consider algorithms alone and motivates the following definition.
\begin{de}[Tower of general algorithms]\label{def:Tower}
	Let $(\Om,\Lambda,\mathcal M,\Xi)$ be a computational problem. A \emph{tower of general algorithms} of height $k$ for $(\Omega,\Lambda,\cM,\Xi)$ is a family $\Gamma_{n_k,n_{k-1},\dots,n_1}:\Om\to\mathcal M$ of general algorithms (where $n_i\in\N$ for $1\leq i \leq k$) such that for all $T\in\Om$
	\begin{align*}
		\Xi(T) = \lim_{n_k\to+\infty}\cdots\lim_{n_1\to+\infty}\Gamma_{n_k,\dots,n_1}(T).
	\end{align*}
\end{de}
If a computational problem requires a tower at least of height $k$ to solve it, then we say the problem has \emph{Solvability Complexity Index} $k$.
Notice that towers of algorithms are commonplace in applications. Indeed, numerical approximation usually consists of a sequence (or a \emph{tower}) of algorithms, with increasingly large input and increasingly precise output. As a prime example, the reader may think of the Finite Element method, in which the input of the algorithm consists of the data of a PDE at the mesh points and the output consists of a sequence of numbers representing the approximate values of the solution at those mesh points.

While precise definitions will follow in the next section, we can now give an informal definition of the SCI Hierarchy. Given a computational problem $(\Omega,\Lambda,\cM,\Xi)$, we define the classes
\begin{align*}
	\Delta_2 &:= \{(\Omega,\Lambda,\cM,\Xi)\,|\,\text{the answer to (1) is yes}\},
	\\
	\Delta_1 &:= \{(\Omega,\Lambda,\cM,\Xi)\,|\,\text{the answer to (2) is yes}\},
	\\
	\Pi_1 &:= \{(\Omega,\Lambda,\cM,\Xi)\in\Delta_2\,|\,\text{there exist error bounds from above}\},
	\\
	\Sigma_1 &:= \{(\Omega,\Lambda,\cM,\Xi)\in\Delta_2\,|\,\text{there exist error bounds from below}\}.
\end{align*}
In addition, higher classes $\Delta_j$, $\Pi_j$, $\Sigma_j$ can be defined, based on towers of algorithms, instead of individual sequences (see Section \ref{sec:preliminaries}). One obtains a sequence of classes 
\begin{align*}
	\Delta_1\subset (\Sigma_1\cap\Pi_1) \subset (\Sigma_1\cup\Pi_1) \subset \Delta_2 \subset (\Sigma_2\cap\Pi_2) \subset \cdots
\end{align*}
into which computational problems fall. This is the SCI Hierarchy.

Within this so-called SCI Hierarchy our main results Theorems \ref{th:convergence_of_algorithm_intro} and \ref{th:operatororm_convergence_intro} mean that the computational spectral problem for the Klein-Gordon equation in $\R^d$ belongs to the class $\Delta_2^A$ if the potential $V$ satisfies Hypothesis \ref{hyp:assumptions_on_V} with $p>d$, and it belongs to the subclass $\Pi_1^A\subset \Delta_2^A$ if the constant $M$ in \eqref{eq:hypothesis_V_bounds} is explicitly known; here the supscript $A$ means that the algorithm to approximate the spectrum can be chosen to be arithmetic. 

\section{Preliminaries}\label{sec:preliminaries}
\subsection{Spectral theory of the Klein-Gordon equation}\label{sec:preliminaries_KG}
The spectral properties of the operator polynomial $T_V$ are intimately related to the spectral properties of the operator polynomial
\begin{align*}
 	L_V(\lambda) := I - \big(S^* - \lambda H_0^{-\f12}\big)\big(S - \lambda H_0^{-\f12}\big), \quad \lambda\in\C.
\end{align*}
which arises from $T_V$ by means of the quasi-similarity transformation $L_V = H_0^{-1/2} T_V H_0^{-1/2}$ and whose values are bounded linear operators. In particular, 
$\sigma_p(T_V)=\sigma_p(L_V)$ and $\sigmaess(T_V)\cap\R=\sigmaess(L_V)\cap\R$ by \cite[Prop.\ 2.3]{MR2268872}. Hence we can use $L_{V}$ to compute the point spectrum of $T_V$. In the sequel we will construct finite-dimensional approximations of $L_V$ in order to approximate~$\sigma_p(L_V)$.

A standard calculation shows that, for $\lambda^2\notin\sigma(H_0)$,
\begin{align}\label{eq:L_factorization}
  L_V(\lambda) = (I-\lambda^2H_0^{-1})\left[ I-K(\lambda) \right]
\end{align}
where
\begin{align}\label{eq:K(lambda)}	
  K(\lambda) := (I-\lambda^2H_0^{-1})^{-1}\Big( S^*S - \lambda\big(S^*H_0^{-\f12} + H_0^{-\f12}S \big) \Big);
\end{align}
note that compactness of $S$ implies compactness of $K(\lambda)$ for all $\lambda$ with  $\lambda^2\in\rho(H_0)$. 
The factorization \eqref{eq:L_factorization} implies that $L_V(\lambda)$ is invertible if both $I-\lambda^2H_0^{-1}$ and $I-K(\lambda)$ are invertible. The first factor is invertible if and only if $\lambda^2\notin\sigma(H_0)$, hence we obtain the  characterisation
\begin{align}
\label{eq:TVK}
  \sigma(T_V)\setminus \{\pm\sqrt{\sigma(H_0)}\}  \!=\!
	\sigma(L_V)\setminus \{\pm\sqrt{\sigma(H_0)}\} \!=\! \{\lambda\!\in\!\C | I\!-\!K(\lambda) \text{ not invertible}\}\setminus\{\pm\sqrt{\sigma(H_0)}\}.
\end{align}
The proof of Theorem \ref{th:convergence_of_algorithm_intro} relies heavily on the following theorem, which we prove in Sections \ref{sec:proofs}-\ref{sect:op_norm_bounds}.
\begin{theorem}[Computable approximation of $K(\lambda)$]\label{th:operatororm_convergence_intro}
	Assume that $V$ satisfies Hypothesis {\rm \ref{hyp:assumptions_on_V}} and let 
	$\lambda\in\C\setminus\{\pm\sqrt{\sigma(H_0)}\} =\C \setminus((-\infty,-m]\cup[m,\infty))$.
	Then there exists a sequence of matrix approximations $(K_n(\lambda))_{n\in\N}$, such that each $K_n(\lambda)$ can be computed in finitely many arithmetic operations from the values $\{V(x)\,|\,x\in\Q^d\}$ \vspace{-2mm} and 
	\begin{align*}
		\big\|K(\lambda) - K_n(\lambda)\big\|_{L^2\to L^2} \leq \f{C}{n},
	\end{align*}
	where $C>0$ is given explicitly in terms of $M$, $m$, $d$, $p$, $\lambda$ $($for details, see Theorem {\rm \ref{th:operatororm_convergence}} below$)$.
\end{theorem}
\subsection{The SCI Hierarchy in depth}
In this section we dive deeper into the theory of the SCI Hierarchy. This will allow us to fully understand the implications of Theorem \ref{th:convergence_of_algorithm_intro} for the computational complexity of the Klein-Gordon eigenvalue problem.
\begin{de}[Recursiveness]\label{def:recursive}
Suppose that for all $f\!\in\!\Lambda$ and for all $T\!\in\!\Omega$ we have $f(T)\!\in\! \R$ or~$\C$. We say that $\Gamma=\Gamma_{n_k,n_{k-1},\dots,n_1}$ is \emph{recursive} if $\Gamma_{n_k,n_{k-1},\dots,n_1}(\{f(T)\}_{f\in\Lambda_\Gamma(T)})$ can be executed by a Blum-Shub-Smale (BSS) machine \cite{BSS} that takes $(n_1,n_2,\dots,n_k)$ as input and that has an oracle that can access  $f(T)$ for any $f\in\Lambda$.
\end{de}
\setcounter{sectiontracker}{\value{section}}
\setcounter{deftracker}{\value{de}}
\setcounter{section}{1}
\setcounter{de}{5}
\begin{example}[continued]\label{ex:SCI_example_continued}
\setcounter{section}{\value{sectiontracker}}
\setcounter{de}{\value{deftracker}}
Let $(\Omega,\Lambda,\cM,\Xi)$ be as in Example \ref{ex:SCI_example}. The following sequence of algorithms has been shown to converge to $\sigma(A)$ for any $A\in\Omega$ (with respect to $d_{\textnormal{H}}$) in \cite{AHS}:
	
	For $n\in\N$ let $\cL_n:=n^{-1}(\Z+\I\Z)\cap B_n(0)$ be a finite lattice in the complex plane and, for 
	$A\in\Omega=\mathcal{S}^\infty(\ell^2(\N))$, let $\Lambda_{\Gamma_n}(A) := \{ f_{ij}\,|\, i,j\leq n \}$. Note that in this case $\Lambda_{\Gamma_n}(A)$ does not actually depend on $A$, hence condition (iii) of Definition \ref{def:Algorithm} is automatically satisfied. The algorithm $\Gamma_n$ is now defined by
	\begin{align}\label{eq:compact_alg_example}
		\Gamma_n(A) := \left\{ z\in\cL_n \,\middle|\, \|(A_n-z)^{-1}\|>n \right\},
	\end{align}
	where $A_n$ denotes the upper left $n\times n$ block obtained by truncating the matrix representation of $A$. Clearly, the action of $\Gamma_n$ only depends on those $f_{ij}$ with $i,j\leq n$, i.e.\ Definition \ref{def:Algorithm} (ii) is satisfied.
	It can also be shown (cf.\ \cite[Th. 7.5 (i)]{AHS}) that $\Gamma_n$ is recursive in the sense of Definition \ref{def:recursive} for any $n\in\N$ and that $\dhaus{\Gamma_n(A)}{\sigma(A)}\to 0$ as $n\to\infty$ for any $A\in\Omega$. 
\end{example}
We stress that convergence of the sequence $\Gamma_n$ is only guaranteed for \emph{compact} operators $A\in\Omega$ as in the example above. Indeed, for larger problem sets (e.g.\ the set of bounded operators $L(\ell^2(\N))$) it can be shown that there does \emph{not} exist \emph{any} sequence of recursive algorithms $\Gamma_n$ such that $\dhaus{\Gamma_n(A)}{\sigma(A)}\to 0$ as $n\to\infty$ for all $A\in L(\ell^2(\N))$ (cf.\ \cite[Thm.\ 7.5 (i)]{AHS}). This motivates the next definition.

\begin{de}[Tower of arithmetic algorithms]\label{def:Arithmetic-Tower}
Given a computational problem $(\Omega,\Lambda,\cM,\Xi)$ where $\Lambda$ is countable, a \emph{tower of arithmetic algorithms}  for $(\Omega,\Lambda,\cM,\Xi)$ is a general tower of algorithms where the lowest mappings $\Gamma_{n_k,\dots,n_1}:\Omega\to\mathcal{M}$ satisfy the following:
For each $T\in\Omega$ the mapping $\N^k\!\ni\!(n_1,\dots,n_k)\!\mapsto\!  \Gamma_{n_k,\dots,n_1}(T)\!=\!\Gamma_{n_k,\dots,n_1}(\{f(T)\}_{f\in\Lambda(T)})$ is recursive, and $\Gamma_{n_k,\dots,n_1}(T)$ is a finite string of complex numbers that can be identified with an element in~$\mathcal{M}$.
\end{de}
\begin{remark}[Types of towers]
One can define many types of towers, see \cite{AHS}. In this paper we write \emph{type $G$} as shorthand for a tower of \emph{general} algorithms, and \emph{type $A$} as shorthand for a tower of \emph{arithmetic} algorithms. If a tower $\{\Gamma_{n_k,n_{k-1},\dots,n_1}\}_{n_i\in\N, 1\leq i\leq k}$ is of type $\tau$ (where $\tau\!\in\!\{A,G\}$ in this paper), we~write
	\begin{equation*}
	\{\Gamma_{n_k,n_{k-1},\dots,n_1}\}\in\tau.
	\end{equation*}
\end{remark}
\begin{remark}[Computations over the reals]
The computations in this paper are assumed to take place over the real numbers, hence the appearance of a BSS machine in Definition \ref{def:recursive}. One can prove that the tower defined in \eqref{eq:compact_alg_example} is in fact a tower of arithmetic algorithms  of height 1 in the sense of Definition~\ref{def:Arithmetic-Tower}.
\end{remark}
\begin{de}[SCI]\label{de:SCI}
	A computational problem $(\Omega,\Lambda,\cM,\Xi)$ is said to have \emph{Solvability Complexity Index $(\SCI)$} $k$ with respect to a tower of algorithms of type $\tau$ if $k$ is the smallest integer for which there exists a tower of algorithms of type $\tau$ of height $k$ for $(\Omega,\Lambda,\cM,\Xi)$. Then we write 
 	$\SCI(\Omega,\Lambda,\cM,\Xi)_\tau=k$.
	If there exists a tower $\{\Gamma_n\}_{n\in\N}\in\tau$ and a finite $N_1\in\N$ with $\Xi=\Gamma_{N_1}$ we set $\SCI(\Omega,\Lambda,\cM,\Xi)_\tau :=0$.
\end{de}
\begin{de}[The SCI Hierarchy]\label{de:hierarchy}
\label{1st_SCI}
The \emph{$\SCI$ Hierarchy} is a hierarchy $\{\Delta_k^\tau\}_{k\in{\N_0}}$ of classes of computational problems $(\Omega,\Lambda,\cM,\Xi)$ where each $\Delta_k^\tau$ is defined as the collection of all computational problems satisfying
\setlength\tablength{0.2\textwidth}
\hspace{2cm} \begin{tabbing} \hspace{0.8\tablength} \= \hspace{1.2\tablength} \= \hspace{0.5\tablength} \= \hspace{\tablength} \=
\kill
	\> $(\Omega,\Lambda,\cM,\Xi)\in\Delta_0^\tau$ \> $\Longleftrightarrow$ \> $\mathrm{SCI}(\Omega,\Lambda,\cM,\Xi)_\tau= 0,$ 
	\\[1mm]
	\> $(\Omega,\Lambda,\cM,\Xi)\in\Delta_{k+1}^\tau$ \>  $\Longleftrightarrow$ \>  $\mathrm{SCI}(\Omega,\Lambda,\cM,\Xi)_\tau\leq k,\qquad k\in\N,$
\end{tabbing}
with the special class $\Delta_1^\tau$  defined as the class of all computational problems in $\Delta_2^\tau$ with known error bounds, i.e.
\hspace{2cm}  \begin{tabbing} \hspace{0.8\tablength} \= \hspace{1.2\tablength} \= \hspace{0.5\tablength} \= \hspace{\tablength} \= \kill
	\> $(\Omega,\Lambda,\cM,\Xi)\in\Delta_{1}^\tau$ \> $\Longleftrightarrow$ \>
	$
		\left\{ \!\!\! \begin{array}{c}  
		\exists\,\{\Gamma_n\}_{n\in \mathbb{N}}\!\in\!\tau\ \exists\,\epsilon_n\!\!\searrow\! 0 
		\text{ s.t. }  \forall \,T\!\in\!\Omega:\\[0.75mm]
		d(\Gamma_n(T),\Xi(T)) \leq \epsilon_n.
	\end{array} \right.$
\end{tabbing}
Hence we have  $\Delta_0^\tau\subset\Delta_1^\tau\subset\Delta_2^\tau\subset\cdots$
\end{de}
When the metric space $\mathcal{M}$ has certain ordering properties, one can define further classes that take into account convergence from below/above and associated error bounds. In order to not burden the reader with unnecessary definitions, we provide the definition that is relevant to the cases where $\mathcal{M}$ is the space of closed and bounded subsets of $\R^d$ together with the Hausdorff distance. These are the cases of relevance to us. A more comprehensive and abstract definition can be found in \cite{AHS}.

\begin{figure}
\centering
\begin{tikzpicture}[scale=0.8]
\draw (0,0) circle (1) (-0.3,.9)  node [text=black,above] {$\Sigma_k$}
      (1,0) circle (1) (1.3,.9)  node [text=black,above] {$\Pi_k$}
      (.5,0) circle (1.95) (1.5,1.9) node [text=black,above] {$\Delta_{k+1}=\{\SCI\leq k\}$}
      (.5,0) node [text=black] {$\Delta_{k}$}
      (0.5,0) circle (0.35);
\end{tikzpicture}
\caption{Schematic representation of the SCI Hierarchy. For $k\in\{1,2,3\}$ it has been proven that the innermost circle is superfluous, i.e.\ $\Delta_k = \Sigma_k\cap\Pi_k$.}
\label{fig:sci-hir}
\end{figure}
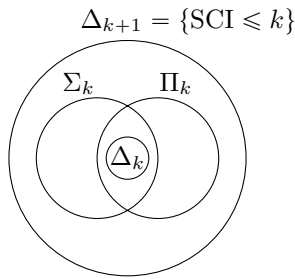

\begin{de}[The SCI Hierarchy -- Hausdorff metric]
\label{def:pi-sigma}
In the setup in Definition~\ref{1st_SCI} assume further that $\mathcal{M}\!=\!(\mathrm{cl}(\C),d)$ where $d\!=\!d_{\mathrm{H}}$. Then, for $k\!\in\!\N$, we can define the following subsets of $\Delta_{k+1}^{\tau}$:
\begin{equation*}
	\begin{split}
	\Sigma_{k}^\tau
	:=
	\Big\{(\Omega,\Lambda,\cM,\Xi) \in \Delta_{k+1}^\tau \ &\big\vert \  \exists\{ \Gamma_{n_k,\dots,n_1}\}\in\tau\text{ s.t. }    \forall T \in \Omega,\, \exists \{X_{n_k}(T)\}\subset\mathcal{M} \text{ s.t. } 
	\\
	&\; \lim_{n_k\to\infty}\cdots\lim_{n_1\to\infty}\Gamma_{n_k,\dots,n_1}(T)=\Xi(T),
	\\
	& \lim_{n_{k-1}\to\infty}\!\cdots\lim_{n_1\to\infty}\Gamma_{n_k,\dots,n_1}(T)\subset X_{n_k}(T),
	\\
	&\; d\left(X_{n_k}(T),\Xi(T)\right)\leq n^{-1} \Big\},
	\\[1mm]
	\Pi_{k}^\tau
	:=
	\Big\{(\Omega,\Lambda,\cM,\Xi) \in \Delta_{k+1}^\tau \ &\big\vert \ \exists\{ \Gamma_{n_k,\dots,n_1}\}\in\tau\text{ s.t. }    \forall T \in \Omega,\, \exists \{X_{n_k}(T)\}\subset\mathcal{M}, \text{ s.t. } 
	\\
	& \lim_{n_k\to\infty}\cdots\lim_{n_1\to\infty}\Gamma_{n_k,\dots,n_1}(T)=\Xi(T),
	\\
	&\; \Xi(T)\subset X_{n_k}(T),
	\\
	&\; d\Bigl(X_{n_k}(T),\lim_{n_{k-1}\to\infty}\cdots\lim_{n_1\to\infty}\Gamma_{n_k,\dots,n_1}(T)\Bigr)\leq n^{-1} \Big\}.
	\end{split}
\end{equation*}
\end{de}
It follows immediately from the definitions that $\Delta_k\subset\Sigma_k\cap\Pi_k$, see Figure \ref{fig:sci-hir}.
In fact, it can be shown that $\Delta_k=\Sigma_k\cap\Pi_k$ for $k\in\{1,2,3\}$, but for higher $k$ equality remains an open problem. 
We  refer to \cite{AHS} for a detailed treatise.
\subsection{Classification of the Klein-Gordon eigenvalue problem}
Let $d\in\N$, $p>d$ and define
\begin{align*}
	\Omega_{p,M} &:= \{ V\in W^{1,p}(\R^d) \,|\, V \text{ satisfies Hypothesis \ref{hyp:assumptions_on_V} with \eqref{eq:hypothesis_V_bounds}} \},
	\\
	\Omega_{p} &:= \bigcup\nolimits_{M>0} \Omega_{p,M},
	\\
	\Lambda &:= \{V \mapsto V(x) \,|\, x\in\Q^d\},
	\\[1mm]
	\cM &:= (\mathrm{cl}(\C),\, d_{\textnormal{H}}),
	\\[1mm]
	\Xi &: \Omega_{p} \to \cM; \quad \Xi(V) = \sigma(T_V). 
\end{align*}
We will commit slight abuse of notation by not distinguishing between $\Xi$ and its restriction to $\Omega_{p,M}$.
The quadruples $(\Omega_{p,M},\Lambda,\cM,\Xi)$, $(\Omega_{p},\Lambda,\cM,\Xi)$ define computational problems in the sense of Definition~\ref{def:computational_problem}. Loosely speaking, $\Omega_{p,M}$ refers to the situation where $M$ is known a-priori and $\Omega_{p}$ refers to the situation where it is not. It is intuitively clear that the former problem should be simpler than the latter. The following corollary, which follows immediately from Theorem \ref{th:convergence_of_algorithm_intro} and Definition \ref{de:hierarchy}, shows that this intuition is in fact correct in a rigorous~sense.
\begin{corollary}\label{cor:classification}
	For $d\in\N$, $p>d$, and $M>0$ as in Hypothesis {\rm \ref{hyp:assumptions_on_V}}, 
	\begin{align*}
		(\Omega_{p},\Lambda,\cM,\Xi) &\in \Delta_2^{\mathrm{A}},
		\\
		(\Omega_{p,M},\Lambda,\cM,\Xi) &\in \Pi_1^{\mathrm{A}}.
	\end{align*}
\end{corollary}
\begin{proof}
	The first inclusion follows immediately from Theorem \ref{th:convergence_of_algorithm_intro} \eqref{eq:convergence_of_alg_intro}. 
	To prove the second claim, let $V\in\Omega_{p,M}$ and let $\Gamma_n$ denote the algorithm given by Theorem \ref{th:convergence_of_algorithm_intro} \eqref{eq:error_bound_intro}. According to Definition~\ref{def:pi-sigma}, we need to find subsets $X_n(V)\subset\C$ such that for all $n\in\N$
	\begin{subequations}
	\begin{align}
		&\sigma(T_V) \subset X_n(V),
		\label{eq:Pi_1_proof1}
		\\
		&\dhaus{X_n(V)}{\Gamma_{n}(V)} \leq n^{-1}.
		\label{eq:Pi_1_proof2}
	\end{align}
	\end{subequations}
	An obvious candidate to satisfy \eqref{eq:Pi_1_proof1}, \eqref{eq:Pi_1_proof2} is the neighbourhood
	$$X_n(V) := \Bigl\{z\in\C\,\Big|\,\dist(z,\Gamma_n(V))\leq\f1n\Bigr\}.$$
	Now \eqref{eq:Pi_1_proof2} is obvious from the definition.\ Moreover, Theorem \ref{th:convergence_of_algorithm_intro} \eqref{eq:error_bound_intro} implies that $\sigma(T_V) \!\subset\! X_n(V)$ for all $n\in\N$, hence \eqref{eq:Pi_1_proof1} is satisfied. Thus the existence of an adequate sequence $X_n(V)\!\subset\!\C$ is shown, which concludes the proof.
\end{proof}
\begin{remark}
	Whenever we talk about an \emph{explicitly known} constant in the following, we mean a constant for which an explicit formula can be given, which involves only quantities that are known a-priori, such as $m$, $d$, $p$, $\lambda$. Constants that additionally depend on $M$ are regarded as explicitly known only if the computational problem at hand is $\Omega_{p,M}$.
\end{remark}

\section{Proofs of the main theorems}\label{sec:proofs}
The next sections are devoted to the proof of Theorem \ref{th:operatororm_convergence_intro}. The strategy is to first choose a basis and truncate the matrix representation of $K(\lambda)$ to finite size (say $l\times l$). Then, in a second step, the matrix elements of $K(\lambda)$ are approximated by computable quantities (depending on a discretisation parameter~$n$) and convergence is proved as $l,n\to\infty$ jointly.

To carry out the outlined strategy, we begin by constructing a finite-dimensional subspace of $\cH$ with an explicit orthonormal basis. To this end, let $n,r\in\N$ and define
\begin{align*}
	\widehat G_{n,r} := (n^{-1}\Z)\cap S_{r}
\end{align*} 
where $S_{r} \!:=\! \{x\in\R^d\,|\,|x_i|\leq r\text{ for all }i=1,\dots,d\}$ (the hat indicates that $\widehat G_{n,r}$ is a lattice in Fourier space).
Clearly, $l_n:=|G_{n,r}| =(2nr+1)^d$. 
For any enumeration $(i_k:k\!=\!1,\dots (2nr\!+\!1)^d)$ of the set~$\widehat G_{n,r}$, we define the~functions
\begin{align}\label{eq:ekn}
	e_k^{(n)}:=n^{\f d2}\cdot\widehat\chi_{i_k+[0,\f1n)^d}, \quad k=1,\dots,l_n,
\end{align}
where $\widehat\cdot$ denotes Fourier transform and the normalisation constant $n^{\f d2}$ is chosen such that $\|e_k^{(n)}\|_{L^2(\R^d)}=1$ for all $n\in\N$. The $e_k^{(n)}$ are linearly independent, smooth functions in $L^2(\R^d)$ and it is easily checked that their first and second derivatives are again in $L^2(\R^d)$. 
We denote their span by $\cH_n$ and the orthogonal projection onto $\cH_n$ by $P_{l_n}$, \vspace{-2mm} i.e.\ 
\begin{align}\label{eq:P_l_def}
	\cH_n \!:=\! \spann \{ e_k^{(n)} \!: k=1,\dots,l_n\} \subset H^1(\R^d), \quad
	P_{l_n}:=\sum_{k=1}^{l_n} \langle e_k^{(n)},\cdot\rangle e_k^{(n)}. 
\\[-8mm] \nonumber
\end{align}
Note that the functions $e_k^{(n)}$ and their $L^{\!\infty}(\R^d)$-norms can be calculated explicitly; indeed, 
\begin{align}
\label{eq:ekn_2}
	e_k^{(n)}(\xi) 
	\!=\! \left(\f{n}{2\pi}\right)^{\f d2}\prod_{j=1}^d\f{\e^{\text{i}\xi_j((i_k)_j+\f1n)} \!-\! \e^{\text{i}\xi_j (i_k)_j}}{\xi_j},
	\quad 
	\|e_k^{(n)}\|_{L^{\!\infty}(\R^d) 
	= (2\pi)^{-\f d2} n^{-\f d2},}
\end{align}
where $(i_k)_j$ denotes the $j$'th component of the vector $i_j$ and $\xi=(\xi_1,\dots,\xi_d)\in\R^d$.
\begin{lemma}\label{lemma:(I-P_n)_bound}
	The set $\{e_k^{(n)}\}_{k=1}^{l_n}$ forms an orthonormal basis of $\cH_n$ and, for any $f\in H^1(\R^d)$ with $xf\in L^2(\R^d)$, $P_{l_n}$ one has the error bound
	\begin{align}\label{eq:(I-P_n)_bound}
		\|(I-P_{l_n})f\|_{L^2(\R^d)}^2 \leq (n\pi)^{-2} \|x f\|_{L^2(\R^d)}^2 + r^{-2} \|\nabla f\|_{L^2(\R^d)}^2.
	\end{align}
	Moreover, if $r = r_n \to \infty$ as $n\to \infty$, then $P_{l_n}\to I$ strongly in $L^2(\R^d)$.
\end{lemma}
\begin{proof}
	Orthonormality follows immediately from the definition of the $e_k^{(n)}$ and the unitarity of the Fourier transform. Moreover, strong convergence in $L^2(\R^d)$ follows from \eqref{eq:(I-P_n)_bound} and smooth approximation. Hence it remains to prove \eqref{eq:(I-P_n)_bound}. The latter follows by Parseval's identity and Poincar\'e's inequality, as the following calculation shows, where $\widecheck\cdot$ denotes the inverse Fourier transform and $\langle\cdot\rangle_{U}$ denotes the mean value over $U$.
	\begin{align*}
		\bigg\|f - \sum_{k=1}^{l_n} \langle e_k^{(n)},f\rangle_{L^2(\R^d)} e_k^{(n)} \bigg\|_{L^2(\R^d)}^2 
		&= \bigg\|\widecheck f - \sum_{k=1}^{l_n} \langle \widecheck e_k^{(n)},\widecheck f\rangle_{L^2(\R^d)} \widecheck e_k^{(n)} \bigg\|_{L^2(\R^d)}^2 
		\\[-2mm]
		&= \bigg\|\widecheck f - n^d \sum_{k=1}^{l_n} \langle\chi_{i_k+[0,\f1n)^d}, \widecheck{f} \rangle_{L^2(\R^d)}\chi_{i_k+[0,\f1n)^d} \bigg\|_{L^2(\R^d)}^2
		\\[-2mm]
		&= \bigg\|\widecheck f - \sum_{k=1}^{l_n} \langle\widecheck{f}\rangle_{i_k+[0,\f1n)^d}\chi_{i_k+[0,\f1n)^d} \bigg\|_{L^2(\R^d)}^2 
		\\[-1mm]
		&= \sum_{k=1}^{l_n} \int_{i_k+[0,\f1n)^d} \big|\widecheck f(\xi) - \langle\widecheck{f}\rangle_{i_k+[0,\f1n)^d}\big|^2\,d\xi + \int_{\R^d\setminus S_{r}} |\widecheck f(\xi)|^2\,d\xi
		\\
		&\leq (n\pi)^{-2}\sum_{k=1}^{l_n} \|\nabla\widecheck f\|_{L^2(i_k+[0,\f1n)^d)}^2 + \int_{\R^d\setminus S_{r}} \f{|\xi|^2}{r^2}|\widecheck f(\xi)|^2\,d\xi
		\\[1mm]
		&\leq (n\pi)^{-2} \|\nabla\widecheck f\|_{L^2(\R^d)}^2 + r^{-2} \|\xi\widecheck f\|_{L^2(\R^d)}^2
		\\[2mm]
		&= (n\pi)^{-2} \|x f\|_{L^2(\R^d)}^2 + r^{-2} \|\nabla f\|_{L^2(\R^d)}^2.
\qedhere
	\end{align*}
\end{proof}
The following lemma shows that $(-\Delta)^{-\f12}$ preserves linear decay of functions.
\begin{lemma}\label{lemma:xDf_bound}
	Let $f\in L^2(\R^d)$ with $\|xf\|_{L^2(\R^d)}<\infty$. Let $a\in W^{1,\infty}(\R^d)$ and define a pseudodifferential operator $D_a$ \vspace{-2mm} by
	\begin{align*}
		D_a f := \big(a\widehat f\big)\widecheck{\;}_.
	\end{align*}
	\vspace{-1mm}Then
	\begin{align}\label{eq:xDf_bound}
		\|xD_a f\|_{L^2(\R^d)} \leq \|\nabla a\|_{L^{\!\infty}(\R^d)}\|f\|_{L^2(\R^d)} + \|a\|_{L^{\!\infty}(\R^d)} \|xf\|_{L^2(\R^d)}.
	\end{align}
\end{lemma}
\begin{proof}
	Repeated use of Parseval's identity gives
	\begin{align*}
		\|x D_a f\|_{L^2(\R^d)} &= \|x (a\widehat f)\widecheck{\;}\|_{L^2(\R^d)} = \|(x (a\widehat f)\widecheck{\;})\widehat{\;}\|_{L^2(\R^d)} = \|\nabla(a\widehat f)\|_{L^2(\R^d)} 
		\\
		&=  \|\widehat f\nabla a + a\nabla\widehat f\|_{L^2(\R^d)} 
		\leq \|\nabla a\|_{L^{\!\infty}(\R^d)} \|\widehat f\|_{L^2(\R^d)} + \|a\|_{L^{\!\infty}(\R^d)} \|\nabla\widehat f\|_{L^2(\R^d)}
		\\
		&\leq \|\nabla a\|_{L^{\!\infty}(\R^d)} \|f\|_{L^2(\R^d)} + \|a\|_{L^{\!\infty}(\R^d)} \|x f\|_{L^2(\R^d)}.
\qedhere
	\end{align*}
\end{proof}
\vspace{1mm}
\subsection{Matrix truncation}
Now we are ready to define a finite-dimensional approximation converging in operator norm to the compact operator $K(\lambda)$ in $L^2(\R^d)$ defined in \eqref{eq:K(lambda)}, 
\begin{align*}
	K(\lambda) &= (I-\lambda^2 H_0^{-1})^{-1}H_0^{-\f12}\big( V^2 - 2\lambda V \big)H_0^{-\f12}, \quad
	\lambda \in \C\setminus \setminus\{\pm\sqrt{\sigma(H_0)}\}.
\end{align*}
An obvious choice for this 
is the compression $P_{l_n}K(\lambda)P_{l_n}$
to the subspace $\cH_n$ spanned by the functions $e_k^{(n)}$ in \eqref{eq:ekn}. By the triangle inequality one has
\begin{align}\label{eq:K-PKP_12}
	\|K(\lambda)-P_{l_n}K(\lambda)P_{l_n}\|_{L^2\to L^2} 
	\leq \|(I-P_{l_n})K(\lambda)P_{l_n}\|_{L^2\to L^2} + \|K(\lambda)(I-P_{l_n})\|_{L^2\to L^2},
\end{align}
where $\|\cdot\|_{L^2\to L^2}$ denotes the operator norm of a bounded linear operator in $L^2(\R^d)$. 
For later reference, we introduce the following notation.
\begin{notation}\label{notation25}
	We define
	\begin{alignat*}{2}
	W_{\!\lambda} &:= V^2 - 2\lambda V, \quad &&\lambda\in\C, \\
	a_\lambda(\xi) &:= \f{(\xi^2+m^2)^{\f12}}{\xi^2+m^2-\lambda^2}, \quad &&\lambda\in\C\setminus\{\pm\sqrt{\sigma(H_0)}\}, \ \xi\in\R^d.
	\end{alignat*}%
\end{notation}
Note that, by Hypothesis \ref{hyp:assumptions_on_V}, we have $V\!\in\! W^{1,p}(\R^d)\cap L^{\!\infty}(\R^d)$ which implies $W_{\!\lambda}\!\in\! W^{1,p}(\R^d)\cap L^{\!\infty}(\R^d)$ and condition $(H_M)$ therein yields $xW_{\!\lambda} \in L^{\!\infty}(\R^d)$ for any $\lambda\in\C$. 
In fact, the following bounds follow immediately from the definition of $W_{\!\lambda}$.
\begin{lemma}\label{lemma:W_bounds}
	If $\,V\!$ satisfies Hypothesis {\rm \ref{hyp:assumptions_on_V}}, then for $W_{\!\lambda} \!=\! V^2 \!-\! 2\lambda V$, $\lambda\!\in\!\C$, the following bounds hold.
\begin{equation}
	\begin{aligned}
		\|W_{\!\lambda}\|_{L^{\!\infty}(\R^d)} &\leq M(M+2|\lambda|),
		&
		\|W_{\!\lambda}\|_{L^p(\R^d)} &\leq M(M+2|\lambda|),
		\\
		\|x W_{\!\lambda}\|_{L^{\!\infty}(\R^d)} &\leq M\Big(\f{M}{2}+2|\lambda|\Big),
		&
		\|\nabla W_{\!\lambda}\|_{L^p(\R^d)} &\leq 2M(M+|\lambda|).
	\end{aligned} 
\end{equation}
\end{lemma}
\begin{proof}
	The first three inequalities follow readily from \eqref{eq:hypothesis_V_bounds} and the triangle inequality. To~prove the last one, we note that $\|\nabla W_{\!\lambda}\|_{L^p(\R^d)} \!=\! \|2V\nabla V \!-\! 2\lambda\nabla V\|_{L^p(\R^d)} \!\leq\! 2(\|V\|_{L^{\!\infty}(\R^d)} \!+\! |\lambda|) \|\nabla V\|_{L^p(\R^d)}$.
\end{proof}
In addition to the above bounds on $W_{\!\lambda}$, we have the following explicit estimates for $a_\lambda$.
\begin{lemma}\label{lemma:a-W-bounds}
	For $\lambda\in \C\setminus\{\pm\sqrt{\sigma(H_0)}\}$ the following bounds  hold.
\begin{equation}\label{eq:a_lambda_bounds}
	\begin{aligned}
		\|a_\lambda\|_{L^{\!\infty}(\R^d)} &\leq \f{|\lambda|+m}{\dist(\lambda^2,[m^2,\infty))},
		&
		\|\nabla a_\lambda\|_{L^{\!\infty}(\R^d)} &\leq \f{(m+|\lambda|)^2+|\lambda|^2}{ \dist(\lambda^2,[m^2,\infty))^2 },
		\\
		\|\xi a_\lambda\|_{L^{\!\infty}(\R^d)} &\leq 1 + \f{|\lambda|^2}{\dist(\lambda^2,[m^2,\infty))},
		&
	\end{aligned}
\end{equation}
	where the choice $\lambda=0$ yields the corresponding estimates for $a_0$.
\end{lemma}
\begin{proof}
Let $\xi \in \R^d$. 
	First inequality: 
	Using the definition in Notation \ref{notation25} we have
	\begin{align*}
		a_\lambda(\xi) &= \f{(\xi^2+m^2)^{\f12}}{\xi^2+m^2-\lambda^2}
		= \f{(\xi^2+m^2)^{\f12}}{((\xi^2+m^2)^{\f12}-\lambda)((\xi^2+m^2)^{\f12}+\lambda)}.
	\end{align*}
	Necessarily, either $\re(\lambda)\geq 0$ or $\re(-\lambda)\geq 0$. Assume without loss of generality that $\re(\lambda)\geq 0$. Then $(\xi^2+m^2)^{\f12}/|(\xi^2+m^2)^{\f12}+\lambda| \leq 1$ and hence
	\begin{align*}
		|a_\lambda(\xi)| &\leq \f{1}{|(\xi^2+m^2)^{\f12}-\lambda|}
		\leq \f{1}{\dist(\lambda,\{\pm\sqrt{\sigma(H_0)}\})}
		\leq \f{|\lambda|+m}{\dist(\lambda^2,[m^2,\infty))}.
	\end{align*}
	Second inequality:
	By explicit calculation and using the above bound for $|a_\lambda(\xi)|$, we find 
	\begin{align*}
		|\nabla a_\lambda(\xi)| &= \left|\f{\xi}{(\xi^2+m^2)^\f12}\f{\xi^2+m^2+\lambda^2}{ (\xi^2+m^2-\lambda^2)^2}\right|
		\leq \left|\f{\xi^2+m^2+\lambda^2}{ (\xi^2+m^2-\lambda^2)^2}\right|
		\\
		&\leq \f{\xi^2+m^2}{ |\xi^2+m^2-\lambda^2|^2} + \f{|\lambda|^2}{ \dist(\lambda^2,[m^2,\infty))^2 }
		\leq \|a_\lambda\|_{L^{\!\infty}(\R^d)}^2 + \f{|\lambda|^2}{ \dist(\lambda^2,[m^2,\infty))^2 }
		\\
		&\leq \f{(|\lambda|+m)^2+|\lambda|^2}{ \dist(\lambda^2,[m^2,\infty))^2 }.
	\end{align*}
	Third inequality: 
	By the definition of $a_\lambda(\xi)$ we have
	\begin{align*}
		|\xi a_\lambda(\xi)| &\!=\! \f{|\xi| (\xi^2\!+\!m^2)^{\f12}}{|\xi^2\!+\!m^2\!-\!\lambda^2|} 
		\!\leq\! \f{\xi^2\!+\!m^2}{|\xi^2\!+\!m^2\!-\!\lambda^2|}
		\!\leq\! 1 \!+\! \f{|\lambda|^2}{|\xi^2\!+\!m^2\!-\!\lambda^2|}
		\!\leq\! 1 \!+\! \f{|\lambda|^2}{\!\dist(\lambda^2,[m^2,\infty))\!}.
	\qedhere	
	\end{align*}
\end{proof}
\vspace{1mm}
The next lemma gives explicit error bounds for the two terms on the right-hand side of \eqref{eq:K-PKP_12}. 
\begin{lemma}\label{lemma:matrix_trunctation_2_bounds}
	For any $\lambda\in\C\setminus\{\pm\sqrt{\sigma(H_0)}\}$ and $n\in\N$,  
	\begin{align}
		\|(I-P_{l_n})K(\lambda)P_{l_n}\|_{L^2\to L^2}^2 &\leq n^{-2}C_1(a_\lambda,\lambda)^2 + r^{-2}C_2(a_\lambda,\lambda)^2,
		\tag{i}
		\label{eq:12_1}
		\\
		\|K(\lambda)(I-P_{l_n})\|_{L^2\to L^2}^2 &\leq \Big(n^{-2}C_1(a_0,\lambda)^2 + r^{-2}C_2(a_0,\lambda)^2\Big) \dist(\lambda^2,[m^2,\infty))^{-2}, 
		\tag{ii}
		\label{eq:12_2}
	\end{align}
	\vspace{-1mm}where 
	\begin{equation}
	\label{eq:C1C2}
	\begin{aligned}
	  C_1(a_\lambda,\lambda) 
	  &:= (m\pi)^{-1}\big( \|\nabla a_\lambda\|_{L^{\!\infty}(\R^d)} \|W_{\!\lambda}\|_{L^{\!\infty}(\R^d)} + \|a_\lambda\|_{L^{\!\infty}(\R^d)} \|xW_{\!\lambda}\|_{L^{\!\infty}(\R^d)} \big), \\ 
		C_2(a_\lambda,\lambda) 
		&:=m^{-1}\|\xi a_\lambda\|_{L^{\!\infty}(\R^d)} \|W_{\!\lambda}\|_{L^{\!\infty}(\R^d)}.
		\end{aligned}
	\end{equation}	
\end{lemma}

\begin{proof}
	\eqref{eq:12_1} \ 
	Let $f\in L^2(\R^d)$ and $\lambda\in\C\setminus\{\pm\sqrt{\sigma(H_0)}\}$. If we define 
	\begin{align*}
		D_{a_\lambda} &:= (I-\lambda^2 H_0^{-1})^{-1}H_0^{-\f12},
	\end{align*}
	then $D_{a_\lambda}$ maps $H^k(\R^d)$ bijectively onto $H^{k+1}(\R^d)$ for any $k\!\in\!\N_0$ and
	$K(\lambda) \!=\! D_{a_\lambda} W_{\!\lambda} H_0^{-\f12}\!$. 
	Clearly, $a_\lambda\!\in\! W^{1,\infty}(\R^d)$ and $\|\xi a_\lambda\|_{L^{\!\infty}(\R^d)}\!<\!\infty$.
	It is not difficult to check, e.g.\ by symbolic differentiation, that
	\begin{align}
	\label{eq:g-alambda}
		\nabla a_\lambda(\xi) = -\xi \f{\xi^2+m^2+\lambda^2}{(\xi^2+m^2)^\f12 (\xi^2+m^2-\lambda^2)^2}, 
		\quad \xi \in \R^d.
	\end{align}
	To simplify notation, we write $g:=W_{\!\lambda} H_0^{-\f12}P_{l_n}f$ and compute
	\begin{align*}
		\|(I-P_{l_n})K(\lambda)P_{l_n}f\|_{L^2\to L^2} 
		&= \|(I-P_{l_n})D_{a_\lambda}g\|_{L^2\to L^2}. 
	\end{align*}
		Since $P_{l_n}f \in \cH_n \subset H^1(\R^d)$, $H_0^{-\f12}P_{l_n}f \in H^2(\R^d)$ and $W_{\!\lambda} \in W^{1,p}(\R^d)\cap L^{\!\infty}(\R^d)$, $xW_{\!\lambda} \in L^{\!\infty}(\R^d)$ by Hypothesis \ref{hyp:assumptions_on_V}, we have
		$g \in L^2(\R^d)$ and $xg \in L^2(\R^d)$. It follows that $D_{a_\lambda}g \in H^1(\R^d)$ and, by Lemma~\ref{lemma:xDf_bound}, that $x D_{a_\lambda}g \in L_2(\R^d)$. 
	Hence we can apply Lemmas \ref{lemma:(I-P_n)_bound} and \ref{lemma:xDf_bound} to obtain
	\begin{align*}
		\|(I\!-\!P_{l_n})D_{a_\lambda}g\|_{L^{\!2}\to L^{\!2}}^2 
		\!&\leq\! (n\pi)^{-2} \|x D_{a_\lambda}g\|_{L^2(\R^d)}^2 \!+\! r^{-2} \|\nabla(D_{a_\lambda}g)\|_{L^2(\R^d)}^2
		\\
		&\leq\! (n\pi)^{-2} \Big( \|\nabla a_\lambda\|_{L^{\!\infty}(\R^d)}\|g\|_{L^{\!2}(\R^d)} \!+\! \|a_\lambda\|_{L^{\!\infty}(\R^d)} \|xg\|_{L^{\!2}(\R^d)} \Big)^{\!2} \!\!\!\!+\! r^{-2} \|\xi(D_{a_\lambda}g)\widehat{\;}\|_{L^{\!2}(\R^d)}^2
		\\
		&\leq\! (n\pi)^{-2} \Big( \|\nabla a_\lambda\|_{L^{\!\infty}(\R^d)}\|g\|_{L^{\!2}(\R^d)} \!+\! \|a_\lambda\|_{L^{\!\infty}(\R^d)} \|xg\|_{L^{\!2}(\R^d)} \Big)^{\!2} \!\!\!+\! r^{-2} \|\xi a_\lambda\widehat g\|_{L^{\!2}(\R^d)}^2
		\\
		&\leq\! (n\pi)^{-2} \Big( \|\nabla a_\lambda\|_{L^{\!\infty}(\R^d)}\|g\|_{L^{\!2}(\R^d)} \!+\! \|a_\lambda\|_{L^{\!\infty}(\R^d)} \|xg\|_{L^{\!2}(\R^d)} \Big)^{\!2} \!\!\!+\! r^{-2} \|\xi a_\lambda\|_{L^{\!\infty}(\R^d)}^2 \|g\|_{L^{\!2}(\R^d)}^2.
	\end{align*}
	The bounds $\|g\|_{L^2(\R^d)} \!\leq\! \|W_{\!\lambda}\|_{L^{\!\infty}(\R^d)} \|H_0^{-\f12}P_{l_n}f\|_{L^2(\R^d)}$, $\|xg\|_{L^2(\R^d)} \!\leq\! \|xW_{\!\lambda}\|_{L^{\!\infty}(\R^d)} \|H_0^{-\f12}P_{l_n}f\|_{L^2(\R^d)}$ now yield that
	\begin{align*}
		\|(I-P_{l_n})Dg\|_{L^2\to L^2}^2 
		&\leq \f{1}{(n\pi)^{2}} \Big( \|\nabla a_\lambda\|_{L^{\!\infty}(\R^d)} \|W_{\!\lambda}\|_{L^{\!\infty}(\R^d)} \|H_0^{-\f12}P_{l_n}f\|_{L^2(\R^d)}
		\\[-1mm]
		&\qquad\qquad + \|a_\lambda\|_{L^{\!\infty}(\R^d)} \|xW_{\!\lambda}\|_{L^{\!\infty}(\R^d)} \|H_0^{-\f12}P_{l_n}f\|_{L^2(\R^d)} \Big)^2 
		\\[-1mm]
		&\qquad\qquad + r^{-2} \, \|\xi a_\lambda\|_{L^{\!\infty}(\R^d)}^2 \|W_{\!\lambda}\|_{L^{\!\infty}(\R^d)}^2 \|H_0^{-\f12}P_{l_n}f\|_{L^2(\R^d)}^2
		\\
		&\leq m^2\Big(n^{-2}C_1(a_\lambda,\lambda)^2 + r^{-2}C_2(a_\lambda,\lambda)^2\Big) \|H_0^{-\f12}\|_{L^2\to L^2}^2 \|f\|_{L^2(\R^d)}^2
		\\[1mm]
		&\leq \Big(n^{-2}C_1(a_\lambda,\lambda)^2 + r^{-2}C_2(a_\lambda,\lambda)^2\Big)\|f\|_{L^2(\R^d)}^2
	\end{align*}
	with $C_i(a_\lambda,\lambda)$, $i=1,2$, as in \eqref{eq:C1C2}.
	
	\eqref{eq:12_2}
	Let $f\in L^2(\R^d)$. To simplify notation, we denote $h := K(\lambda)(I-P_{l_n})f$.
	Using the selfadjointness of $P_{l_n}$ and $H_0$, we compute
	\begin{align}
		\|h\|_{L^2(\R^d)}^2 
		&= \langle K(\lambda)(I-P_{l_n})f,h\rangle_{L^2(\R^d)}
		\nonumber
		\\
		&= \langle f,(I-P_{l_n})K(\lambda)^*h\rangle_{L^2(\R^d)}
		\nonumber
		\\
		&= \big\langle f, \, (I-P_{l_n}) H_0^{-\f12} W_{\bar\lambda} H_0^{-\f12}(I-\overline\lambda^2 H_0^{-1})^{-1}h \big\rangle_{L^2(\R^d)}
		\nonumber
		\\
		&\leq \|f\|_{L^2(\R^d)} \big\|(I\!-\!P_{l_n}) H_0^{-\f12} W_{\bar\lambda} H_0^{-\f12}(I\!-\!\overline\lambda^2 H_0^{-1})^{-1}h\big\|_{L^2(\R^d)}.
		\label{eq:K(I-P)_calculation}
	\end{align}
	The last term above can be estimated by a similar calculation as in the proof of \eqref{eq:12_1}. 
	In fact, since $h\!\in\! L^2(\R^d)$ and $W_{\!\lambda} \!\in\! W^{1,p}(\R^d) \cap L^{\!\infty}(\R^d)$, we have 
	$H_0^{-\f12} W_{\bar\lambda} H_0^{-\f12} (I\!-\!\overline\lambda^2 H_0^{-1})^{-1}h \!\in\! H^1(\R^d)$.
	Further, since $x W_{\bar\lambda} \!\in\! L^{\!\infty}(\R^d)$, we have 
	$x W_{\bar\lambda} H_0^{-\f12} (I-\overline\lambda^2 H_0^{-1})^{-1}h \!\in\! L^2(\R^d)$
	and so, by Lemma~\ref{lemma:xDf_bound} for $D_{a_0}\!=\!H_0^{-1/2}\!$, also
	$x H_0^{-\f12} W_{\bar\lambda} H_0^{-\f12} (I\!-\!\overline\lambda^2 H_0^{-1})^{-1}h \!\in\! L^2(\R^d)$ .
	Then, as in \eqref{eq:12_1}, one can show~that
	\begin{align}
		\big\|(I-P_{l_n}) H_0^{-\f12} W_{\bar\lambda} H_0^{-\f12} & (I-\overline\lambda^2 H_0^{-1})^{-1}h\big\|_{L^2(\R^d)}^2 \leq
		\nonumber
		\\
		&\leq m^2\Big(n^{-2}C_1(a_0,\overline\lambda)^2 + r^{-2}C_2(a_0,\overline\lambda)^2\Big) \big\|H_0^{-\f12}(I-\overline\lambda^2 H_0^{-1})^{-1}h\big\|_{L^2(\R^d)}^2
		\nonumber
		\\
		&\leq \Big(n^{-2}C_1(a_0,\overline\lambda)^2 + r^{-2}C_2(a_0,\overline\lambda)^2\Big) \big\|(I-\overline\lambda^2 H_0^{-1})^{-1}\big\|_{L^2(\R^d)}^2 \|h\|_{L^2(\R^d)}^2.
		\label{eq:h_norm_eq}
	\end{align}
	Using \eqref{eq:h_norm_eq} in \eqref{eq:K(I-P)_calculation} and dividing by $\|h\|_{L^2(\R^d)}$, we obtain
	\begin{align*}
		\|h\|_{L^2(\R^d)}
		&\leq \|f\|_{L^2(\R^d)} \Big(n^{-2}C_1(a_0,\overline\lambda)^2 + r^{-2}C_2(a_0,\overline\lambda)^2\Big)^{\f12} \big\|(I-\overline\lambda^2 H_0^{-1})^{-1}\big\|_{L^2(\R^d)}. 
	\end{align*}
	Therefore,
	\begin{align*}
		\|K(\lambda)(I - P_{l_n})f\|_{L^2(\R^d)}^2 &= \|h\|_{L^2(\R^d)}^2 \\
		&\leq \|f\|_{L^2(\R^d)}^2 \Big(n^{-2}C_1(a_0,\overline\lambda)^2 + r^{-2}C_2(a_0,\overline\lambda)^2\Big) \big\|(I-\overline\lambda^2 H_0^{-1})^{-1}\big\|_{L^2(\R^d)}^2
		\\
		&\leq \|f\|_{L^2(\R^d)}^2 \Big(n^{-2}C_1(a_0,\overline\lambda)^2 + r^{-2}C_2(a_0,\overline\lambda)^2\Big)\dist(\lambda^2,[m^2,\infty))^{-2}.
	\end{align*}
	To conclude the proof, we note that $C_i(a_0,\overline\lambda) = C_i(a_0,\lambda)$ for $i=1,2$.
\end{proof}
Using Lemma \ref{lemma:matrix_trunctation_2_bounds} in \eqref{eq:K-PKP_12}, we immediately conclude our first main error bound.
\begin{prop}[Matrix truncation error]
	\label{prop:matrix_trunctation_error}
	If $\,V\!$ satisfies Hypothesis {\rm \ref{hyp:assumptions_on_V}}, the following error bound holds for all $\lambda\in\C\setminus\{\pm\sqrt{\sigma(H_0)}\}$;
	\begin{align}
		\|K(\lambda)-P_{l_n}K(\lambda)P_{l_n}\|_{L^2\to L^2}^2
		&\leq
		n^{-2}\Big( C_1(a_\lambda,\lambda)^2 + C_1(a_0,\lambda)^2\dist(\lambda^2,[m^2,\infty))^{-2} \Big)
		\nonumber
		\\
		&\quad + r^{-2}\Big( C_2(a_\lambda,\lambda)^2 + C_2(a_0,\lambda)^2\dist(\lambda^2,[m^2,\infty))^{-2}\Big).
		\label{eq:matrix_trunctation_error}
	\end{align}
	In particular, if $r\xrightarrow{n\to\infty} \infty$, then $\|K(\lambda)-P_{l_n}K(\lambda)P_{l_n}\|_{L^2\to L^2}\to 0$ as $n\to\infty$. 
\end{prop}
\begin{proof}
	This follows immediately from Lemma \ref{lemma:matrix_trunctation_2_bounds}. Note that if $V$ satisfies Hypothesis \ref{hyp:assumptions_on_V}, then both $\|W_{\!\lambda}\|_{L^{\!\infty}(\R^d)}$ and $\|xW_{\!\lambda}\|_{L^{\!\infty}(\R^d)}$ are finite.
\end{proof}
\begin{de}\label{de:C_trunc}
	For later reference we define the \emph{truncation constant}
	\begin{align}\label{eq:truncation_constant}
		C_{\textnormal{trunc}}^\lambda &:= \bigg(C_1(a_\lambda,\lambda)^2 + C_2(a_\lambda,\lambda)^2 + \f{C_1(a_0,\lambda)^2}{\dist(\lambda^2,[m^2,\infty))^{2}} + \f{C_2(a_0,\lambda)^2}{\dist(\lambda^2,[m^2,\infty))^{2}} \bigg)^{\f12},
	\end{align}
	where $\lambda\in\C\setminus\{\pm\sqrt{\sigma(H_0)}\}$. Recall that $C_1(a_\lambda,\lambda)$, $C_2(a_\lambda,\lambda)$ were defined in Lemma \ref{lemma:matrix_trunctation_2_bounds}.
\end{de}
\section{Error bounds for matrix elements}\label{sec:error_bounds_matrix_elements}
In this section we prove several technical lemmas containing error bounds for the computation of the matrix elements $\langle e_k^{(n)},K(\lambda)e_m^{(n)}\rangle_{L^2(\R^d)}$ of the compression $P_{l_n} K(\lambda) P_{l_n}$ to the finite dimensional subspaces $\cH_n$ defined in \eqref{eq:P_l_def}, \eqref{eq:ekn}.
To simplify the presentation, we fix  some notation.
\begin{notation}
	For $k,n\in\N$ and $i_k\in\f1n\Z^d$, let
	\begin{enumi}
		\item $Q_k := i_k+[0,\f1n)^d$,
		\item $\chi_k := \chi_{Q_k}$.
	\end{enumi}
\end{notation}
An explicit formula for the  matrix elements
$\langle e_k^{(n)},K(\lambda)e_m^{(n)}\rangle_{L^2(\R^d)}$, $k,m=1,\dots,l_n$
is provided by the following calculation. 
\begin{align}
	\big\langle e_k^{(n)},K(\lambda)e_m^{(n)} \big\rangle_{L^2(\R^d)} &= \big\langle e_k^{(n)},(I-\lambda^2 H_0^{-1})^{-1}H_0^{-\f12}W_{\!\lambda} H_0^{-\f12} e_m^{(n)} \big\rangle_{L^2(\R^d)} 
	\nonumber
	\\
	&= \big\langle (I-\lambda^2 H_0^{-1})^{-1}H_0^{-\f12} e_k^{(n)},W_{\!\lambda} H_0^{-\f12} e_m^{(n)} \big\rangle_{L^2(\R^d)} 
	\nonumber
	\\
	&= \big\langle (a_\lambda n^{\f d2}\chi_k)\widecheck{\;},W_{\!\lambda} (a_0 n^{\f d2}\chi_m)\widecheck{\;} \big\rangle_{L^2(\R^d)} 
	\nonumber
	\\
	&= \int_{\R^d} \overline{E^\lambda_k(x)}E_m^0(x)W_{\!\lambda}(x)\,dx 
	\label{eq:matrix_elements_of_K}
\end{align}
where we have defined
\begin{align}\label{eq:E_k_def}
	E^\lambda_k(x) := (2\pi)^{-\f d2} n^{\f d2} \int_{Q_k} \!\!a_\lambda(\xi)\,\e^{\I\xi x}\,d\xi, \quad 
	x\in\R^d,
\end{align}
for $\lambda\in\C\setminus\{\pm\sqrt{\sigma(H_0)}\}$.
The integral above is not computable in finitely many arithmetic operations, therefore, we need the following.
\begin{lemma}\label{lemma:E_bounds}
	There exist computable approximations $E^\lambda_{k,N}$ of $E^\lambda_{k}$ such that
	\begin{align}
	\begin{split}\label{eq:E_bounds}
		\|E^\lambda_{k} - E^\lambda_{k,N}\|_{L^2(\R^d)} &\leq \f{\sqrt 2}{nN}\|\nabla a_\lambda\|_{L^{\!\infty}(\R^d)},
		\\
		\|E^\lambda_{k,N}\|_{L^2(\R^d)} &\leq \|a_\lambda\|_{L^{\!\infty}(\R^d)},
	\end{split}
	\end{align}
	for any $k,m,N\in\N$ with $a_\lambda$ defined as in Notation {\rm \ref{notation25}}.
\end{lemma}
\begin{proof}
	Let $k\!\in\!\N$ and $\lambda\!\in\!\C\setminus\{\pm\sqrt{\sigma(H_0)}\}$. For $N\!\in\!\N$ choose a new lattice inside $Q_k$ as $L_N^k \!:=\! Q_k\cap (\f{1}{Nn}\Z^d)$ (see Figure \ref{fig:L_N}). Then, clearly, \vspace{-2mm}$|L_N^k|=N^d$.
	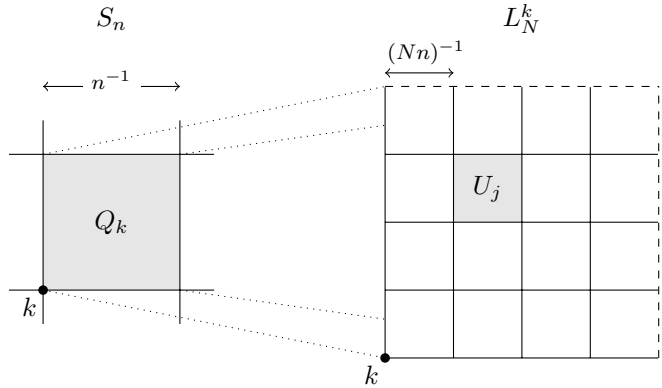
\begin{figure}[htbp]
		\centering
%
%
%

\begin{tikzpicture}[scale=0.9]

\draw (1,4) node{$S_n$};
\draw (7,4) node{$L_N^k$};

\fill[fill=gray!20!white] (0,0) rectangle (2,2);
\draw (1,1) node{$Q_k$};
\fill[fill=black] (0,0) circle(0.07);
\draw (-0.2,-0.25) node{$k$};
\draw[step=2] (-0.5,-0.5) grid (2.5,2.5);

\draw[dotted] (0,2) -- (5,3)
			  (0,0) -- (5,-1)
			  (2,2) -- (9,3)
			  (2,0) -- (9,-1);
\fill[white] (5,-1) rectangle (9,3);
\fill[fill=gray!20!white] (6,1) rectangle (7,2);
\draw (5,-1) grid (8.99,2.99);
\draw[dashed] (5,3) -- (9,3)
			  (9,-1) -- (9,3);
\fill[fill=black] (5,-1) circle(0.07);
\draw (4.8,-1.25) node{$k$};
\draw[<->] (0,3) -- (2,3);
\draw (1,3.1) node[fill=white]{\footnotesize $n^{-1}$};
\draw[<->] (5,3.2) -- (6,3.2);
\draw (5.6,3.5) node{\footnotesize $(Nn)^{-1}$};

\draw (6.5,1.5) node{$U_j$};

\end{tikzpicture}%
 
		\vspace{-2mm}
		\caption{Sketch of the lattice $L_N^k$.}
		\label{fig:L_N}
	\end{figure}
	
	Next, we define the approximation~by
	\begin{align}
		E^\lambda_{k,N}(x) &:= (2\pi)^{-\f d2} n^{\f d2} \sum_{j\in L_N^k} a_\lambda(j)\int_{U_j} \e^{-\I\xi x}\,d\xi
		\label{eq:E_kN1}
		\\[-1.5mm]
		&= n^{\f d2} \bigg( \sum_{j\in L_N^k} a_\lambda(j) \chi_{U_j} \bigg)^{\!\wedge} (x)
		\\[-0.5mm]
		&= N^{-\f d2} \sum_{j\in L_N^k} a_\lambda(j) e^{(Nn)}_j.
		\label{eq:E_kN3}
	\end{align}
	Now the error bound is obtained by the estimates
	\begin{align*}
		\|E^\lambda_{k} - E^\lambda_{k,N}\|_{L^2(\R^d)}^2 &= \|\widecheck E^\lambda_{k} - \widecheck E^\lambda_{k,N}\|_{L^2(\R^d)}^2
		= n^d \Big\| a_\lambda - \sum_{j\in L_N^k} a_\lambda(j)\chi_{U_j} \Big\|_{L^2(\R^d)}^2
		\\
		&\leq n^d \sum_{j\in L_N^k} |U_j|\, \|a_\lambda - a_\lambda(j)\|_{L^{\!\infty}(U_j)}^2
		\leq n^d \f{2}{(nN)^2} \|\nabla a_\lambda\|_{L^{\!\infty}(\R^d)}^2 \sum_{j\in L_N^k} (Nn)^{-d}
		\\
		&\leq \f{2}{n^2N^2} \|\nabla a_\lambda\|_{L^{\!\infty}(\R^d)}^2. 
	\end{align*}
	The corresponding stability bound is obtained similarly,
	\begin{align*}
		\|E_{k,N}^\lambda\|_{L^2(\R^d)} &= n^{\f d2}\bigg\|\bigg ( \sum_{j\in L_N^k} a_\lambda(j)\chi_{U_j} \bigg)^{\!\wedge}\bigg\|_{L^2(\R^d)}
		= n^{\f d2}\bigg\| \sum_{j\in L_N^k} a_\lambda(j)\chi_{U_j} \bigg\|_{L^2(\R^d)}
		\\
		&\leq n^{\f d2}\|a_\lambda\|_{L^{\!\infty}(\R^d)} \bigg(\sum_{j\in L_N^k} \int_{U_j} dx \bigg)^\f12
		= \|a_\lambda\|_{L^{\!\infty}(\R^d)}.
		\qedhere
	\end{align*}
\end{proof}
%
The next step in computing the integrals in \eqref{eq:matrix_elements_of_K} is to pass to a bounded domain. To this end, we let $R>0$ and decompose $\R^d = [-R,R)^d\cup\big( \R^d\setminus [-R,R)^d \big)$. Due to \eqref{eq:matrix_elements_of_K}, this gives us
\begin{align}
	\big\langle e_k^{(n)},K(\lambda)e_m^{(n)} \big\rangle_{L^2(\R^d)} 
	&= \int_{\R^d} \overline{E^\lambda_k(x)}E_m^0(x)W_{\!\lambda}(x)\,dx 
	\nonumber
	\\
	&= \int_{[-R,R)^d} \overline{E^\lambda_k(x)}E_m^0(x)W_{\!\lambda}(x)\,dx + \int_{\R^d\setminus [-R,R)^d}\! \overline{E^\lambda_{k}(x)}E_m^0(x)W_{\!\lambda}(x)\,dx
	\nonumber
	\\
	&= \int_{[-R,R)^d} \overline{E^\lambda_{k,N}(x)}E_{m,N}^0(x)W_{\!\lambda}(x)\,dx
	 + F_1 + F_2 + F_3 + F_4
	 \label{eq:F4_def}
\end{align}
where, to keep the notation simple, we introduce the error terms 
\begin{equation}
\label{eq:errterms}
   \eps_{k,N}^\lambda := E_k^\lambda - E_{k,N}^\lambda
\end{equation}	
and define $F_i$, $i=1,\dots,4$, as
\begin{equation}
\label{eq:F1234}
\begin{aligned}
	 F_1 &:= \int_{[-R,R)^d} \overline{\eps^\lambda_{k,N}(x)}E_{m,N}^0(x)W_{\!\lambda}(x)\,dx,
	 \quad
	 F_2 := \int_{[-R,R)^d} \overline{E^\lambda_{k,N}(x)}\eps_{m,N}^0(x)W_{\!\lambda}(x)\,dx,
	 \\
	 F_3 &:= \int_{[-R,R)^d} \overline{\eps^\lambda_{k,N}(x)}\eps_{m,N}^0(x)W_{\!\lambda}(x)\,dx,
	\quad
	 F_4 := \int_{\R^d\setminus [-R,R)^d}\! \overline{E^\lambda_{k}(x)}E_m^0(x)W_{\!\lambda}(x)\,dx.
\end{aligned}
\end{equation}
The new error terms $F_i$, $i=1,\dots,4$, are indeed small as the next lemma shows.
\begin{lemma}\label{lemma:F_bounds}
	With the bound $M>0$ as in condition $(H_M)$ on $V$, the following error bounds hold for all \vspace{-2mm}$k,m\in\{1,\dots,l_n\}$.
	\begin{align*}
		|F_1| &\leq \f{\sqrt 2}{nN} \|\nabla a_\lambda\|_{L^{\!\infty}(\R^d)} \| a_0\|_{L^{\!\infty}(\R^d)} \|W_{\!\lambda}\|_{L^{\!\infty}(\R^d)},
		\\
		|F_2| &\leq \f{\sqrt 2}{nN} \|\nabla a_0\|_{L^{\!\infty}(\R^d)} \| a_\lambda\|_{L^{\!\infty}(\R^d)} \|W_{\!\lambda}\|_{L^{\!\infty}(\R^d)},
		\\
		|F_3| &\leq \f{2}{n^2N^2} \|\nabla a_0\|_{L^{\!\infty}(\R^d)} \| \nabla a_\lambda\|_{L^{\!\infty}(\R^d)} \|W_{\!\lambda}\|_{L^{\!\infty}(\R^d)},
		\\
		|F_4| &\leq 2M\f{M+|\lambda|}{R}\|a_\lambda\|_{L^{\!\infty}(\R^d)} \|a_0\|_{L^{\!\infty}(\R^d)}.
	\end{align*}
\end{lemma}
\begin{proof}
	Applying H\"older's inequality to $F_1$ and  Lemma \ref{lemma:E_bounds}, we conclude that
	\begin{align*}
		|F_1| \leq \|E_k^\lambda - E_{k,N}^\lambda\|_{L^2(\R^d)} \|E_{m,N}^0\|_{L^2(\R^d)} \|W_{\!\lambda}\|_{L^{\!\infty}(\R^d)}
		\leq \f{\sqrt 2}{nN} \|\nabla a_\lambda\|_{L^{\!\infty}(\R^d)} \|a_0\|_{L^{\!\infty}(\R^d)} \|W_{\!\lambda}\|_{L^{\!\infty}(\R^d)}.
	\end{align*}
	The calculations for $F_2$ and $F_3$ are analogous. The bound for $F_4$ follows readily from the bound $\|E_k^\lambda\|_{L^2(\R^d)}\leq \|a_\lambda\|_{L^{\!\infty}(\R^d)}$ and from Hypothesis \ref{hyp:assumptions_on_V}.
\end{proof}
The first term on the right-hand side of \eqref{eq:F4_def} has a computable integrand and can thus be computed to arbitrary precision using any standard quadrature formula. 
Finally, we are  able to define our approximation to $P_{l_n} K(\lambda) P_{l_n}$ and prove operator-norm convergence. 
\begin{de}[Approximate matrix elements]\label{de:approximate_matrix_elements}
	For $\lambda\!\in\!\C$, $n,N,s\!\in\!\N$, $k,m\!\in\!\{1,\dots,l_n\}$ and $R\!>\!0$ define 
	\begin{align*}
		K_{km}^{(N,R,s)}(\lambda) := \operatorname{Quad}_{R,s}\big(\overline{E^\lambda_{k,N}} E_{m,N}^0 W_{\!\lambda}\big)
	\end{align*}
	where $\operatorname{Quad}_{R,s}$ stands for any quadrature formula converging to $\int_{[-R,R)^d} \overline{E^\lambda_{k,N}}E_{m,N}^0W_{\!\lambda}\,dx$ as $s\to \infty$. By $K^{(N,R,s)}(\lambda)$ we denote the corresponding linear operator on $\cH_n$ with matrix elements $K_{km}^{(N,R,s)}(\lambda)$, $k,m=1,\dots,l_n$.
\end{de}
In practice, to improve numerical performance, one would choose a quadrature rule with a high order of convergence. In our proofs below, however, we will take $\operatorname{Quad}_{R,s}$ to be a standard Riemann sum. Naturally, this approximation achieves only linear order of convergence, but it applies to a wide range of functions because of its low regularity requirements. More concretely, we will use the formula
\begin{align}
	\operatorname{Quad}_{R,s}(f) &:= \sum_{i\in (\f1s\Z^d)\cap[-R,R)^d} f(i)\big|i+[0,\tfrac1s)^d\big|
	= s^{-d}\sum_{i\in (\f1s\Z^d)\cap[-R,R)^d} f(i).
	\label{eq:quadrature_def}
\end{align}
Lemma \ref{lemma:integral_approx} below shows that it is indeed enough to assume $f\in W^{1,p}([-R,R)^d)$ with $p>d$ for $\operatorname{Quad}_{R,s}(f)$ to converge.
\begin{lemma}\label{lemma:integral_approx}
 	If $R\in\N$, then, for every $f\in W^{1,p}([-R,R)^d)$ with $p>d$, 
 	\begin{align*}
 		\bigg| s^{-d}\!\!\!\!\sum_{i\in (\f1s\Z^d)\cap[-R,R)^d}\!\! f(i) \,-\, \int_{(-R,R)^d} f(x)\,dx \bigg| \leq \f{2s^{-1}}{1-\f dp} (2R)^{d-\f dp}  \|\nabla f\|_{L^p([-R,R)^d)}.
 	\end{align*}
 \end{lemma}
 \begin{proof}
 	To simplify  notation, we denote $G_s := (\f1s\Z^d)\cap[-R,R)^d$.
 	Since $R\in\N$, we can write
 	\begin{align*}
 		\int_{[-R,R)^d} f(x)\,dx = \sum_{i\in G_s} \int_{[0,\f1s)^d+i} f(x)\,dx.
 	\end{align*}
 	Then, comparing the two sums term by term,  for $i\in (\f1s\Z^d)\cap[-R,R)^d$ we have
 	\begin{align*}
 		\bigg|\int_{[0,\f1s)^d+i} f(x)\,dx - s^{-d} f(i)\bigg|
 		&= \bigg|\int_{[0,\f1s)^d+i} f(x) - f(i) \,dx\bigg|
 		\\
 		&\leq \int_{[0,\f1s)^d+i} |f(x) - f(i)| \,dx
 		\leq \int_{[0,\f1s)^d+i} \f{2s^{-1+\f dp}}{1-\f dp}\|\nabla f\|_{L^p([0,\f1s)^d+i)} \,dx
 		\\
 		&= \f{2s^{-1+\f dp}}{1-\f dp}\|\nabla f\|_{L^p([0,\f1s)^d+i)} \int\limits_{[0,\f1s)^d+i} \,dx
 		= \f{2s^{-1+\f dp}}{1-\f dp} s^{-d} \|\nabla f\|_{L^p([0,\f1s)^d+i)} 
 	\end{align*}
 	where Morrey's inequality was used in the third line (cf.\ (28) in the proof of \cite[Thm.\ 9.12]{Brezis}). Summing these inequalities over $(\f1s\Z^d)\cap[-R,R)^d$, we finally obtain
 	\begin{align*}
 		\bigg| s^{-d}\sum_{i\in G_s} f(\xi) - \!\int_{[-R,R)^d}\! f(x)\,dx \bigg| 
 		&\leq \sum_{i\in G_s} \bigg|\int_{[0,\f1s)^d+i}\! f(x)\,dx - s^{-d} f(i)\bigg|
 		\\
 		&\leq \sum_{i\in G_s} \f{2s^{-1+\f dp}}{1-\f dp} s^{-d} \|\nabla f\|_{L^p([0,\f1s)^d+i)} 
 		\\
 		&= \f{2s^{-1+\f dp}}{1-\f dp} s^{-d} \sum_{i\in G_s} 1\cdot \|\nabla f\|_{L^p([0,\f1s)^d+i)} 
 		\\
		&\leq \f{2s^{-1+\f dp}}{1-\f dp} s^{-d} (2Rs)^{\f d{p'}} \bigg( \sum_{i\in G_s} \|\nabla f\|_{L^p([0,\f1s)^d+i)}^p\bigg)^{\f1p}
		\\
		&= \f{2s^{-1+\f dp}}{1-\f dp} s^{-d} (2Rs)^{d-\f dp}  \|\nabla f\|_{L^p([-R,R)^d)}
		\\
		&= \f{2s^{-1+\f dp}}{1-\f dp} s^{-\f dp}(2R)^{d-\f dp}  \|\nabla f\|_{L^p([-R,R)^d)}
		\\
		&= \f{2s^{-1}}{1-\f dp} (2R)^{d-\f dp}  \|\nabla f\|_{L^p([-R,R)^d)}
 	\end{align*}
 	where we have used H\"older's inequality in the fourth line.
 \end{proof}
The next lemma contains our fifth and final error bound. Note that this is the first time that Hypothesis \ref{hyp:assumptions_on_V} (i) is used, i.e.\ $V\in W^{1,p}(\R^d)$.
\begin{lemma}\label{lemma:F5_bound}
	Denote by $F_5 := \int_{[-R,R)^d} \overline{E^\lambda_{k,N}}E_{m,N}^0W_{\!\lambda}\,dx - \operatorname{Quad}_{R,s}(\overline{E^\lambda_{k,N}} E_{m,N}^0 W_{\!\lambda})$ the quadrature error. Then
	\begin{align*}
		|F_5| \leq \f{2s^{-1}}{1-\f dp} n^{-d} r(2R)^{d-\f dp} 
		\|a_\lambda\|_{L^{\!\infty}(\R^d)} & \|a_0\|_{L^{\!\infty}(\R^d)}  (2\pi)^{-d}
		\Big( 2\sqrt{2} \|W_{\!\lambda}\|_{L^p([-R,R)^d)} + r^{-1}\|\nabla W_{\!\lambda}\|_{L^p([-R,R)^d)} 
		\Big).
	\end{align*}
\end{lemma}
\begin{proof}
	Recall formula \eqref{eq:E_kN3} for $E_{k,N}^{\lambda}$ and note that
	\begin{align*}
		\|E_{k,N}^\lambda\|_{L^{\!\infty}(\R^d)} 
		&= N^{-\f d2} \bigg\| \sum_{j\in L_N^k} a_\lambda(j) e^{(Nn)}_j \bigg\|_{L^{\!\infty}(\R^d)} 
		\leq N^{-\f d2}\|a_\lambda\|_{L^{\!\infty}(\R^d)} \sum_{j\in L_N^k} \|e^{(Nn)}_j\|_{L^{\!\infty}(\R^d)}
		\\
		&\leq (2\pi)^{-\f d2}N^{\f d2}(Nn)^{-\f d2} 
		= \|a_\lambda\|_{L^{\!\infty}(\R^d)}(2\pi)^{-\f d2} n^{-\f d2},
		\\
		\|\nabla E_{k,N}^\lambda\|_{L^{\!\infty}(\R^d)} 
		&\leq N^{-\f d2} \bigg\| \sum_{j\in L_N^k} a_\lambda(j) \nabla e^{(Nn)}_j \bigg\|_{L^{\!\infty}(\R^d)}
		\leq N^{-\f d2}\|a_\lambda\|_{L^{\!\infty}(\R^d)} \sum_{j\in L_N^k} \|\nabla e^{(Nn)}_j\|_{L^{\!\infty}(\R^d)}
		\\
		&\leq (2\pi)^{-\f d2}N^{\f d2}(Nn)^{-\f d2} r
		= \|a_\lambda\|_{L^{\!\infty}(\R^d)}(2\pi)^{-\f d2} n^{-\f d2}\sqrt 2 r
	\end{align*}
	where the extra factor $\sqrt{2} r$ comes from the fact that $|\xi|\!\leq\! \sqrt{2}r$ in \eqref{eq:E_kN1}. By Hypothesis \ref{hyp:assumptions_on_V} we have $W_{\!\lambda}\!\in\! W^{1,p}(\R^d)$ with $p\!>\!d$. Hence applying Lemma \ref{lemma:integral_approx} with $f\!=\!\overline{E^\lambda_{k,N}} E_{m,N}^0 W_{\!\lambda}$, we~obtain
	\begin{align*}
		|F_5| 
		&\leq\! \f{2s^{-1}}{1\!-\!\f dp} (2R)^{d-\f dp}  \big\|\nabla (\overline{E^\lambda_{k,N}}E_{m,N}^0W_{\!\lambda})\big\|_{L^p([-R,R)^d)}
		\\
		&\leq\! \f{2s^{-1}}{1\!-\!\f dp} (2R)^{d-\f dp}  \Big( 
			\|\nabla E^\lambda_{k,N}\|_{L^{\!\infty}(\R^d)} \|E_{m,N}^0\|_{L^{\!\infty}(\R^d)} \|W_{\!\lambda}\|_{L^p([-R,R)^d)} 
			\\[-2mm]
			&\hspace{3cm} + \|E^\lambda_{k,N}\|_{L^{\!\infty}(\R^d)} \|\nabla E_{m,N}^0\|_{L^{\!\infty}(\R^d)} \|W_{\!\lambda}\|_{L^p([-R,R)^d)} 
			\\
			&\hspace{3cm} + \|E^\lambda_{k,N}\|_{L^{\!\infty}(\R^d)} \|E_{m,N}^0\|_{L^{\!\infty}(\R^d)} \|\nabla W_{\!\lambda}\|_{L^p([-R,R)^d)} 
		\Big)
		\\
		&\leq\! \f{2s^{-1}}{1\!-\!\f dp} (2R)^{d-\f dp}  \Big( 
			\|a_\lambda\|_{L^{\!\infty}(\R^d)}(2\pi)^{-\f d2} n^{-\f d2}\sqrt 2 r \|a_0\|_{L^{\!\infty}(\R^d)} (2\pi)^{-\f d2} n^{-\f d2} \|W_{\!\lambda}\|_{L^p([-R,R)^d)} 
			\\[-2mm]
			&\hspace{3cm} + \|a_\lambda\|_{L^{\!\infty}(\R^d)}(2\pi)^{-\f d2} n^{-\f d2} \|a_0\|_{L^{\!\infty}(\R^d)}(2\pi)^{-\f d2} n^{-\f d2}\sqrt 2 r \|W_{\!\lambda}\|_{L^p([-R,R)^d)} 
			\\
			&\hspace{3cm} + \|a_\lambda\|_{L^{\!\infty}(\R^d)}(2\pi)^{-\f d2} n^{-\f d2} \|a_0\|_{L^{\!\infty}(\R^d)}(2\pi)^{-\f d2} n^{-\f d2} \|\nabla W_{\!\lambda}\|_{L^p([-R,R)^d)} 
		\Big)
		\\
		&=\! \f{2s^{-1}}{1\!-\!\f dp} (2R)^{d-\f dp} 
		\|a_\lambda\|_{L^{\!\infty}(\R^d)} \|a_0\|_{L^{\!\infty}(\R^d)} (2\pi)^{-d} n^{-d}
		 \Big( 
			2\sqrt{2} r \|W_{\!\lambda}\|_{L^p([-R,R)^d)} \!+\! \|\nabla W_{\!\lambda}\|_{L^p([-R,R)^d)} 
		\Big)
		\\
		&=\! \f{2s^{-1}}{1\!-\!\f dp} n^{\!-d} r(2R)^{d\!-\!\f dp} 
		\|a_\lambda\|_{L^{\!\infty}(\R^d)} \|a_0\|_{L^{\!\infty}(\R^d)} (2\pi)^{\!-d}
		 \Big( 
			2\sqrt{2} \|W_{\!\lambda}\|_{L^{\!p}([-\!R,R)^d)} \!\!+\! r^{-1}\|\nabla W_{\!\lambda}\|_{L^{\!p}([-\!R,R)^d)}\! 
		\Big). 
		\\[-14mm]
	\end{align*}
\end{proof}
%
%
\section{Operator norm bounds for the approximation of $K(\lambda)$}
\label{sect:op_norm_bounds}
In this section we use the technical bounds from the previous section to obtain operator norm error bounds for the difference\footnote{
 We commit slight abuse of notation by not distinguishing between $K^{(N,R,s)}(\lambda)$ (acting in $\cH_n$) and $K^{(N,R,s)}(\lambda)\oplus 0$ (acting in $\cH$).
 }
\begin{align}\label{eq:PKP-Kn}
	P_{l_n}K(\lambda)P_{l_n} - K^{(N,R,s)}(\lambda).
\end{align}
Combined with Proposition \ref{prop:matrix_trunctation_error}, this will provide a fully computable approximation of $K(\lambda)$ with explicitly known error bounds.
To this end, we use the following version of Young's inequality, which follows from the Riesz-Thorin interpolation theorem (see e.g.\ \cite[Thm. 0.3.1]{Sogge}).
\begin{lemma}
\label{lemma:young}
	If $A$ is a bounded operator on $\ell^2(\N)$ with matrix elements $A_{ij}$, then
	\begin{align*}
		\|A\| &\leq F^{\f12} \tilde F^{\f12}
	\vspace{-2mm}	
	\end{align*}%
	\vspace{-2mm}where
	\begin{align*}
		F := \sup_{i\in\N} \sum_{j=1}^\infty |A_{ij}|,\;
		\qquad
		\widetilde F := \sup_{j\in\N} \sum_{i=1}^\infty |A_{ij}|.
	\end{align*}
\end{lemma}
Combining Lemma \ref{lemma:young} with the results from Section \ref{sec:error_bounds_matrix_elements}, we obtain the following error bound,
\begin{prop}[Matrix error bound]\label{prop:matrix_error_bound}
	Let $n,N,s\!\in\!\N$, $R,r>0$ and $\lambda\!\in\!\C\setminus\{\pm\sqrt{\sigma(H_0)}\}$.~Then 
	\begin{align*}
		\big\|P_{l_n}K(\lambda)P_{l_n} \!-\! K^{(N,R,s)}(\lambda)\big\|_{L^2\to L^2}
		\!\leq\! 
		D_1\f{(2nr\!+\!1)^d}{nN} \!+\! D_2\f{(2nr\!+\!1)^d\!}{R} \!+\! D_3\f{(2nr\!+\!1)^d r R^{d-\f dp}\!}{n^ds}
	\end{align*}
	where
	\begin{align}\label{eq:C1C2C3}
	\begin{split}
		D_1 \!&:=\! \sqrt{2}M(M+2|\lambda|) \Bigl(\|\nabla a_\lambda\|_{L^{\!\infty}(\R^d)} \| a_0\|_{L^{\!\infty}(\R^d)} \!+\! \| a_\lambda\|_{L^{\!\infty}(\R^d)} \|\nabla a_0\|_{L^{\!\infty}(\R^d)} 
		 \!+\! \sqrt 2\|\nabla a_\lambda\|_{L^{\!\infty}(\R^d)} \|\nabla a_0\|_{L^{\!\infty}(\R^d)}  \Bigr) ,
		\\[1mm]
		D_2 \!&:=\! 2M(M+|\lambda|)\|a_\lambda\|_{L^{\!\infty}(\R^d)} \|a_0\|_{L^{\!\infty}(\R^d)}, 
		\\
		D_3 \!&:=\! \f{2^{1+d-\f dp}}{(1-\f dp)(2\pi)^{d}} 4\sqrt{2} M(M+|\lambda|) 
		\|a_\lambda\|_{L^{\!\infty}(\R^d)} \|a_0\|_{L^{\!\infty}(\R^d)}.
	\end{split}
	\end{align}
\end{prop}
\begin{proof}
	Denote $K_{km}(\lambda) := \langle e_k^{(n)},K(\lambda)e_m^{(n)}\rangle_{L^2(\R^d)}$, where $k,l=1,\dots,l_n$.
	Applying Lemma \ref{lemma:young} to $A=P_{l_n}K(\lambda)P_{l_n} - K^{(N,R,s)}(\lambda)$, we find that we need to estimate the sums 
	\begin{align*}
		F = \sup_{k\in\{1,\dots,l_n\}} \sum_{m=1}^{l_n} \big|K_{km}(\lambda) - K_{km}^{(N,R,s)}(\lambda)\big|,
		\qquad
		\widetilde F = \sup_{m\in\{1,\dots,l_n\}} \sum_{k=1}^{l_n} \big|K_{km}(\lambda) - K_{km}^{(N,R,s)}(\lambda)\big|.
	\end{align*}
	By \eqref{eq:F4_def}, we have
	\begin{align*}
		\big|K_{km}(\lambda) - K_{km}^{(N,R,s)}(\lambda)\big| 
		&\leq |F_1|+|F_2|+|F_3|+|F_4|+|F_5|
	\end{align*}
	and we can estimate $|F_i|$, $i=1,\dots,5$, by Lemmas \ref{lemma:F_bounds} and \ref{lemma:F5_bound}, 
	respectively. Since all bounds therein are 
	independent of $k,m$, summation over $m$ (or $k$) will multiply the right-hand side with a factor of $l_n = (2nr+1)^d$. Thus we obtain the same bound for both $F$ and $\tilde F$, namely
	\begin{align*}
		F, \tilde F
		&\leq (2nr+1)^d(|F_1|+|F_2|+|F_3|+|F_4|+|F_5|) 
		\\
		&\leq \sqrt 2\f{(2nr+1)^d}{nN} \bigg(\|\nabla a_\lambda\|_{L^{\!\infty}(\R^d)} \| a_0\|_{L^{\!\infty}(\R^d)}  + \| a_\lambda\|_{L^{\!\infty}(\R^d)} \|\nabla a_0\|_{L^{\!\infty}(\R^d)} 
		\\[-2mm]
		&\hspace{7.02cm} + \f{\sqrt 2}{nN}\|\nabla a_\lambda\|_{L^{\!\infty}(\R^d)} \|\nabla a_0\|_{L^{\!\infty}(\R^d)}  \bigg) \|W_{\!\lambda}\|_{L^{\!\infty}(\R^d)}
		\\[-1mm]
		& \quad+ (2nr+1)^d2 M\f{M+|\lambda|}{R}\|a_\lambda\|_{L^{\!\infty}(\R^d)} \|a_0\|_{L^{\!\infty}(\R^d)}
		\\
		& \quad + \f{(2nr+1)^d}{(2\pi n)^{d}} \f{2 r(2R)^{d-\f dp}}{s(1-\f dp)}  
		\|a_\lambda\|_{L^{\!\infty}(\R^d)} \|a_0\|_{L^{\!\infty}(\R^d)} 
		 \Big( 2\sqrt{2} \|W_{\!\lambda}\|_{L^{\!p}([-R,R)^d)} + r^{-1}\|\nabla W_{\!\lambda}\|_{L^{\!p}([-R,R)^d)} 
		\!\Big). 
	\end{align*}
	Equations \eqref{eq:C1C2C3} now follow immediately from the $W_{\!\lambda}$ bounds in Lemma \ref{lemma:W_bounds}.
\end{proof}
\vspace{1mm}
\begin{cor}[Convergence of matrix approximation]\label{cor:matrix_convergence}
	Let $(r_n)_{n\in\N}$ be any sequence in $\N$ and let $\lambda\in\C\setminus\{\pm\sqrt{\sigma(H_0)}\}$. If we define sequences $R_n,N_n,s_n$ by
	\begin{align}\label{eq:parameter_choices}
	\begin{split}
		N_n &:=  (2nr_n\!+\!1)^d,  
		\\[1mm]
		R_n &:=  (2nr_n\!+\!1)^d n, 
		\\
		s_n &:= (2nr_n\!+\!1)^d n^{1-d} r_n R_n^{d-\f dp}.  
	\end{split}
	\end{align}
	then there exists an explicit constant $C_{\textnormal{mat}}^\lambda\!=\!D_1 \!+\! D_2 \!+\! D_3\!>\!0$ 
	with $D_i$, $i\!=\!1,2,3$, as in \eqref{eq:C1C2C3} such~that
	\begin{align*}
		\big\|P_{l_n}K(\lambda)P_{l_n} - K^{(N_n,R_n,s_n)}(\lambda)\big\|_{L^2\to L^2} \leq \f{C_{\textnormal{mat}}^\lambda}{n}.
	\end{align*}
	In particular $\|P_{l_n}K(\lambda)P_{l_n} - K^{(N_n,R_n,s_n)}(\lambda)\|_{L^2\to L^2}\to 0$ as $n\to\infty$.
\end{cor}
\begin{proof}
	The claim follows immediately from Proposition \ref{prop:matrix_error_bound} applied with the chosen sequences $N_n$, $R_n$ and $s_n$, $n\in\N$.
\end{proof}
Finally, we can make our first main result, Theorem \ref{th:operatororm_convergence_intro} stated earlier, more precise and, at the same time, give a constructive proof of it.
\begin{theorem}[Approximation of $K(\lambda)$]\label{th:operatororm_convergence}
Suppose $V$ satisfies Hypothesis {\rm \ref{hyp:assumptions_on_V}}.
	For $n\in\N$ let $r_n = n$ and let $N_n,R_n,s_n$ be chosen as in \eqref{eq:parameter_choices}. Then, for every $\lambda\in\C\setminus\{\pm\sqrt{\sigma(H_0)}\}$, there exists a constant $C_K^\lambda=C_{\textnormal{trunc}}^\lambda+C_{\textnormal{mat}}^\lambda>0$ such \vspace{-1mm} that
	\begin{align}\label{eq:matrix_appr_bound}
		\big\|K(\lambda) - K^{(N_n,R_n,s_n)}(\lambda)\big\|_{L^2\to L^2} \leq \f{C_K^\lambda}{n}
	\end{align}
	where $K^{(N_n,R_n,s_n)}(\lambda) = \big(K_{km}^{(N_n,R_n,s_n)}(\lambda)\big)_{k,m=1}^{l_n}$ is as in Definition {\rm\ref{de:approximate_matrix_elements}} and the constants 	$C_{\textnormal{trunc}}^\lambda$, $C_{\textnormal{mat}}^\lambda$ are as in Definition {\rm\ref{de:C_trunc}} and Corollary {\rm\ref{cor:matrix_convergence}}.
	If the bound $M$ on $V$ in $(H_M)$ is known a-priori $($i.e.\ if $\,V\in\Omega_{p,M})$, then $C_K^\lambda$ is explicitly known.
\end{theorem}
\begin{proof}[Proof of Theorem {\rm \ref{th:operatororm_convergence}} $($and Theorem {\rm \ref{th:operatororm_convergence_intro}}$)$]
	The estimate \eqref{eq:matrix_appr_bound} follows if we combine Proposition \ref{prop:matrix_trunctation_error} and Corollary \ref{cor:matrix_convergence}. Because $r_n = n$ by assumption, we have $N_n ={\rm O}(n^{2d}) \to \infty$, $R_n ={\rm O}(n^{2d+1})\to \infty$
	and, since $p>d$, $s_n =  {\rm O}(n^{d+2} n^{(2d+1) (d-\f dp})) = {\rm O}(n^{2d^2+1}) \to \infty$ as $n\to\infty$. Now \eqref{eq:matrix_appr_bound} follows if we
choose $C_K^\lambda = C_{\textnormal{trunc}}^\lambda + C_{\textnormal{ma}}^\lambda$.
Theorem \ref{th:operatororm_convergence_intro} is obtained by taking $K_n = K^{(N_n,R_n,s_n)}$ for $n\in \N$.
The last claim follows if we note that the bounds for $C_{\textnormal{trunc}}^\lambda$ and $C_{\textnormal{mat}}^\lambda$ depend only on $M$ and quantities that can be estimated by $M$ by means of condition $(H_M)$ such as $\|W_{\!\lambda}\|_{L^{\!\infty}(\R^d)}$ or $\|x W_{\!\lambda}\|_{L^{\!\infty}(\R^d)}$ since $W_{\!\lambda}=V^2-2\lambda V$.
\end{proof} 
\subsection{Improvements under stronger assumptions on $V$}
\label{sec:improvements}
The two-layered approximation scheme in Lemma \ref{lemma:E_bounds} and Figure \ref{fig:L_N} gives precise error control and is convenient in the proofs of Theorems~\ref{th:operatororm_convergence_intro} and~\ref{th:convergence_of_algorithm_intro}. In practice, however, the double discretisation in both $n$ and $N$ may be cumbersome to implement and computationally expensive. Lemma \ref{lemma:V_in_L2_bounds} below shows that, under stronger assumptions on $V$, it is sufficient to choose $N=1$, i.e.\ the second layer of approximation is unnecessary.
\begin{lemma}\label{lemma:V_in_L2_bounds}
	Assume, in addition to Hypothesis {\rm \ref{hyp:assumptions_on_V}}, that $V\!\in\! L^1(\R^d)$ with $\|V\|_{L^1(\R^d)}\!\leq\! M$. Then 
	\begin{align}\label{eq:V_in_L2_bounds}
	\begin{split}
		|F_1| &\leq \f{\sqrt 2}{n} (2\pi n)^{-d} \|\nabla a_\lambda\|_{L^{\!\infty}(\R^d)} \|a_0\|_{L^{\!\infty}(\R^d)} M(M+2|\lambda|),
		\\
		|F_2| &\leq \f{\sqrt 2}{n} (2\pi n)^{-d} \|a_\lambda\|_{L^{\!\infty}(\R^d)} \|\nabla a_0\|_{L^{\!\infty}(\R^d)} M(M+2|\lambda|),
		\\
		|F_3| &\leq \f{2}{n^2} (2\pi n)^{-d} \|\nabla a_\lambda\|_{L^{\!\infty}(\R^d)} \|\nabla a_0\|_{L^{\!\infty}(\R^d)} M(M+2|\lambda|).
	\end{split}
	\end{align}
\end{lemma}
\begin{proof}
	For $N=1$, we have $L_N^K=Q_K$ and hence from the definitions \eqref{eq:E_k_def} and \eqref{eq:E_kN1}
	it immediately follows \vspace{-2mm}that
	\begin{align*}
    		\|E_k^\lambda\|_{L^{\!\infty}(\R^d)} &\leq (2\pi n)^{-\f d2} \|a_\lambda\|_{L^{\!\infty}(\R^d)}, \\
		\|E^\lambda_{k} - E^\lambda_{k,1}\|_{L^{\!\infty}} &\leq (2\pi n)^{-\f d2} \f{\sqrt 2}{n} \|\nabla a_\lambda\|_{L^{\!\infty}(\R^d)}.	
	\end{align*}
	Using these two inequalities in the definition of $F_1$, see \eqref{eq:F1234}, \eqref{eq:errterms}, and applying H\"older's inequality,  we obtain
	\begin{align*}
		|F_1| &\leq \|E^\lambda_{k} - E^\lambda_{k,1}\|_{L^{\!\infty}} \|E_k^0\|_{L^{\!\infty}(\R^d)} \|W_{\!\lambda}\|_{L^1}
		\\
		&\leq \f{\sqrt 2}{n} (2\pi n)^{-d} \|\nabla a_\lambda\|_{L^{\!\infty}(\R^d)} \|a_0\|_{L^{\!\infty}(\R^d)} \big(\|V^2\|_{L^1} + 2|\lambda|\|V\|_{L^1}\big)
		\\
		&\leq \f{\sqrt 2}{n} (2\pi n)^{-d} \|\nabla a_\lambda\|_{L^{\!\infty}(\R^d)} \|a_0\|_{L^{\!\infty}(\R^d)} M(M+2|\lambda|)
	\end{align*}
	where we have used the inequality $\|V^2\|_{L^1} \leq \|V\|_{L^{\!\infty}(\R^d)}\|V\|_{L^1} \leq M^2$ in the last line.
	The bounds for $|F_2|$ and $|F_3|$ follow from analogous estimates.
\end{proof}
Repeating the proof of Proposition \ref{prop:matrix_error_bound} with the bounds for $F_1$, $F_2$ and $F_3$ given by \eqref{eq:V_in_L2_bounds} and $N=1$, we find that, with appropriate parameter choices, a convergence rate of $n^{-1/2}$ can be achieved in \eqref{eq:matrix_appr_bound}. 

\begin{remark}
\label{rem:improved} 
If, in addition to Hypothesis {\rm \ref{hyp:assumptions_on_V}}, $V\!\in\! L^1(\R^d)$ with $\|V\|_{L^1(\R^d)}\!\leq\! M$ and we choose
\begin{align}\label{eq:improved_parameter_choices}
\begin{split}
	r_n &= n^{\f1{2d}} - (2n)^{-1}
\end{split}
\end{align}%
	and $R_n, s_n$  as in \eqref{eq:parameter_choices}, 
	then there exists an explicit constant $C_{\textnormal{mat}}^\lambda=D_1 + D_2 + D_3>0$ with 
	\begin{align*}
		D_1 = \f{\sqrt{2}}{\pi^d} M(M+2|\lambda|) \Big(  
			\|\nabla a_\lambda\|_{L^{\!\infty}(\R^d)} \|a_0\|_{L^{\!\infty}(\R^d)}  
			&+  \|a_\lambda\|_{L^{\!\infty}(\R^d)} \|\nabla a_0\|_{L^{\!\infty}(\R^d)}
			+  \sqrt{2}\|\nabla a_\lambda\|_{L^{\!\infty}(\R^d)} \|\nabla a_0\|_{L^{\!\infty}(\R^d)}
		\Big)
	\end{align*}
	and $D_2$, $D_3$, as in \eqref{eq:C1C2C3}
	such \vspace{-2mm} that
	\begin{align*}
		\big\|P_{l_n}K(\lambda)P_{l_n} - K^{(N_n,R_n,s_n)}(\lambda)\big\|_{L^2\to L^2} \leq \f{C_{\textnormal{mat}}^\lambda}{n^{\f 12}}.
	\end{align*}
	In particular $\|P_{l_n}K(\lambda)P_{l_n} - K^{(N_n,R_n,s_n)}(\lambda)\|_{L^2\to L^2}\to 0$ as $n\to\infty$.
\end{remark}
\section{Definition of the algorithm}\label{sec:def_of_algorithm}
In this section we define an algorithm $\Gamma_n$ to compute $\sigma(T_V)$ for any $V$ satisfying Hypothesis~\ref{hyp:assumptions_on_V}. In order to do this, we have to control the dependence of our above estimates on the spectral parameter. 
\subsection{Preparations}\label{sec:preps}
We begin with a series of technical lemmas.
\begin{lemma}\label{lemma:Lipschitz_continuity}
	The operator-valued function $K$ is locally Lipschitz continuous on $\C\setminus\{\pm\sqrt{\sigma(H_0)}\}$, i.e.\ for any compact subset $B\subset\C\setminus\{\pm\sqrt{\sigma(H_0)}\}$ there exists $C_L>0$ such that 
	\begin{align*}
		\|K(\lambda)-K(\mu)\|_{L^2\to L^2} \leq C_L|\lambda-\mu|
	\end{align*}
	for all $\lambda,\mu\!\in\! B$. Moreover, if $\varrho,\rho\!>\!0$ are such that $|\lambda|,|\mu| \!\leq\! \rho$ and $\dist(\lambda^2,[m^2,\infty))$, $\dist(\mu^2,[m^2,\infty)) \!\geq\! \varrho$, then one can choose
	\begin{align}\label{eq:Lipschitz_constant}
		C_L = \f{M}{m^2}\Big(1+\f{\rho^2}{\varrho}\Big)\left[ \f{2R}{m^2}\Big( 1+\f{\rho^2}{\varrho} \Big)(M+2\rho) + 2 \right].
	\end{align}
\end{lemma}
\begin{proof}
	Denote $R(\lambda^2) := (I-\lambda^2 H_0^{-1})^{-1}$ and $W_{\!\lambda} := W_1 - \lambda W_2$, where $W_1 := V^2$ and $W_2 := 2V$. Then $K(\lambda) = R(\lambda^2)H_0^{-\f12}(W_1-\lambda W_2)H_0^{-\f12}$ and thus, for $\lambda,\mu\in\C$,
	\begin{align*}
		K(\lambda) - K(\mu) &= R(\lambda^2)H_0^{-\f12}(W_1-\lambda W_2)H_0^{-\f12} - R(\mu^2)H_0^{-\f12}(W_1-\mu W_2)H_0^{-\f12}
		\\
		&= (R(\lambda^2) - R(\mu^2))H_0^{-\f12}(W_1-\lambda W_2)H_0^{-\f12} + R(\mu^2)(\mu-\lambda)H_0^{-\f12} W_2 H_0^{-\f12}
		\\
		&= (\lambda^2-\mu^2)H_0^{-1}R(\lambda^2)R(\mu^2)H_0^{-\f12}(W_1-\lambda W_2)H_0^{-\f12} - R(\mu^2)(\lambda-\mu)H_0^{-\f12} W_2 H_0^{-\f12}
		\\
		&= (\lambda-\mu)R(\mu^2)\big[ (\lambda+\mu)H_0^{-1}R(\lambda^2)H_0^{-\f12}(W_1-\lambda W_2) - H_0^{-\f12}W_2 \big]H_0^{-\f12};
	\end{align*}
	here we have used a resolvent identity for operators of the form $(I-\lambda^2 H_0^{-1})^{-1}$ in the third line. Now the assertion follows because $R(\mu^2)[ (\lambda+\mu)H_0^{-1}R(\lambda^2)H_0^{-1/2}(W_1-\lambda W_2) - H_0^{-1/2}W_2 ]H_0^{-1/2}$ is a bounded operator for $\lambda,\mu$ in any compact subset of $\C\setminus\{\pm\sqrt{\sigma(H_0)}\}$.
	
	The explicit bound for $C_L$ follows from the estimates
	\begin{align*}
		C_L &\leq \big\|R(\mu^2)\big[ (\lambda+\mu)H_0^{-1}R(\lambda^2)H_0^{-\f12} W_{\!\lambda} H_0^{-\f12} - H_0^{-\f12} W_2 H_0^{-\f12}\big]\big\|_{L^2\to L^2}
		\\
		&\leq \|R(\mu^2)\|_{L^2\to L^2} \left[ \f{|\lambda+\mu|}{m^4}\|R(\lambda^2)\|_{L^2\to L^2} \|W_{\!\lambda}\|_{L^2\to L^2} + \f1{m^2}\|W_2\|_{L^2\to L^2} \right]
		\\
		&\leq\f 1{m^2}\left( 1+\f{|\mu|^2}{\dist(\mu^2,[m^2,\infty))} \right) \left[ \f{|\lambda+\mu|}{m^2}\left( 1+\f{|\lambda|^2}{\dist(\lambda^2,[m^2,\infty))} \right)M(M+2|\lambda|) + 2M \right]
		\\
		&\leq \f{M}{m^2}\Big(1+\f{\rho^2}{\varrho}\Big)\left[ \f{2R}{m^2}\Big( 1+\f{\rho^2}{\varrho} \Big)(M+2\rho) + 2 \right];
	\end{align*}
	here we have used $\|R(\mu^2)\|\!\leq\! 1\!+\!|\mu|^2\dist(\mu^2,[m^2,\infty))^{-1}$ and the bounds from Lemma \ref{lemma:W_bounds}.
\end{proof}
The $\lambda$-dependence of $C_K^\lambda$ can now be controlled by means of Lemma \ref{lemma:a-W-bounds}. The following proposition yields a precise result.
\begin{prop}\label{prop:C_K_bound}
	For any $\rho>0$ there exists a constant $A>0$ $($independent of $\lambda)$ such that,
	for all $\lambda\in \overline{B_\rho(0)}\setminus\{\pm\sqrt{\sigma(H_0)}\}$,
	the constants $C_{\textnormal{trunc}}^\lambda$ in Definition {\rm \ref{de:C_trunc}} and 
	$C_{\textnormal{mat}}^\lambda$ in Corollary {\rm \ref{cor:matrix_convergence}} satisfy
	\begin{align*}
		C_{\textnormal{trunc}}^\lambda &\leq \f{A}{2\dist(\lambda^2,[m^2,\infty))^2},
		\quad
		C_{\textnormal{mat}}^\lambda \leq \f{A}{2\dist(\lambda^2,[m^2,\infty))^2},
	\end{align*}
	and \vspace{-1mm} therefore
	\begin{align}
		C_K^\lambda &\leq \f{A}{\dist(\lambda^2,[m^2,\infty))^2}.
		\label{eq:C_K_bound}
	\end{align}
		If $M$ is known a-priori $($i.e.\ if $V\in\Omega_{p,M})$, then $A$ is known explicitly.
\end{prop}
\begin{proof}
	Let $\rho>0$ be fixed. For the sake of readability we use the notation
	\begin{align*}
		d(\lambda) := \dist(\lambda^2,[m^2,\infty)), \quad \lambda \in \C \setminus\{\pm\sqrt{\sigma(H_0)}\}.
	\end{align*}
	We begin with $C_{\textnormal{trunc}}^\lambda$. By Definition \ref{de:C_trunc}, we have
	\begin{align*}
		C_{\textnormal{trunc}}^\lambda &:= \bigg(C_1(a_\lambda,\lambda)^2 + C_2(a_\lambda,\lambda)^2 + \f{C_1(a_0,\lambda)^2}{d(\lambda)^{2}} + \f{C_2(a_0,\lambda)^2}{d(\lambda)^{2}} \bigg)^{\f12}
	\end{align*}
	where, as in \eqref{eq:C1C2}, we have $C_1(a_\lambda,\lambda) = (m\pi)^{-1}\big( \|\nabla a_\lambda\|_{L^{\infty}(\R^d)} \|W_{\lambda}\|_{L^{\infty}(\R^d)} + \|a_\lambda\|_{L^{\infty}(\R^d)} \|xW_{\lambda}\|_{L^{\infty}(\R^d)} \big)$ and $C_2(a_\lambda,\lambda) = m^{-1}\|\xi a_\lambda\|_{L^{\infty}(\R^d)} \|W_{\lambda}\|_{L^{\infty}(\R^d)}$. 
	Applying Lemmas \ref{lemma:W_bounds} and \ref{lemma:a-W-bounds}, we can estimate
	\begin{align}\label{eq:C1C2_bounds}
	\begin{split}
		C_1(a_\lambda,\lambda) &\leq \f{M}{m\pi}\left( 
			\f{(m+|\lambda|)^2+|\lambda|^2}{ d(\lambda)^2 } 
			(M+2|\lambda|)
			+ \f{|\lambda|+m}{d(\lambda)}
			\Big(\f{M}{2}+2|\lambda|\Big)  
		\right),
		\\
		C_2(a_\lambda,\lambda) &\leq \f{M}{m}\left(1 + \f{|\lambda|^2}{d(\lambda)}\right)(M+2|\lambda|),
		\\
		C_1(a_0,\lambda) &\leq \f{M}{\pi m^2}\left( 
			\f{M+2|\lambda|}{m} + 
			\f{M}{2}+2|\lambda|
		\right),
		\\
		C_2(a_0,\lambda) &\leq \f{M}{m}(M+2|\lambda|).
	\end{split}
	\end{align}
	Inserting equations \eqref{eq:C1C2_bounds} into the formula for $C_{\textnormal{trunc}}^\lambda$ 
	and using $d(\lambda) \leq |\lambda|^2+m^2 \leq \rho^2+m^2$ for $\lambda\in \overline{B_\rho(0)}$, 
	we readily obtain that there exists a constant $\widetilde C$ such that, for all  $\lambda\in \overline{B_\rho(0)}\setminus\{\pm\sqrt{\sigma(H_0)}\}$ one has $C_{\textnormal{trunc}}^\lambda \leq \widetilde C/d(\lambda)^2$
%
%
Similarly, to bound $C_{\textnormal{mat}}^\lambda$,  we insert the bounds \eqref{eq:a_lambda_bounds} into the formulas \eqref{eq:C1C2C3} for $D_1$, $D_2$, $D_3$. 
	we arrive at
	\begin{align*}
		D_1 &\leq M(M+2|\lambda|)\left( \f{(m+|\lambda|)^2 + |\lambda|^2}{md(\lambda)^2} + \f{|\lambda|+m}{m^{2} d(\lambda)} + \sqrt{2}\f{(m+|\lambda|)^2 + |\lambda|^2}{m^2 d(\lambda)^2} \right),
		\\
		D_2 &\leq 2M(M+2|\lambda|)\f{|\lambda|+m}{md(\lambda)},
		\\[-1mm]
		D_3 &\leq \f{2^{\f72-\f dp}}{(1-\f dp)\pi^{d}} M(M+|\lambda|) \f{|\lambda|+m}{md(\lambda)}.
	\end{align*}
	Thus, for $C_{\textnormal{mat}} = D_1+D_2+D_3$, there exists a constant $C'$ such that $C_{\textnormal{mat}} \leq C'/d(\lambda)^2$
	for all $\lambda\in \overline{B_\rho(0)}\setminus\{\pm\sqrt{\sigma(H_0)}\}$.
	The proof is complete if we set $A:=2\max\{\widetilde C,C'\}$.
\end{proof}
\subsection{Abstract bounds on the eigenvalues}\label{sec:spectral_bound}
\subsubsection{Bounds on the non-real eigenvalues}
In this subsection we return to the general setting of Sections~\ref{sec:preliminaries_KG_intro} and \ref{sec:preliminaries_KG}. We will derive a bound for $\sigma_p(T_V)\setminus \R$, which holds whenever $VH_0^{-\f12}$ decomposes into a sum of a compact operator and a strict contraction.

In a Hilbert space $\cH$, consider a self-adjoint, uniformly positive operator $H_0 \!\ge\! m \!>\!0$ and a symmetric operator $V$ with $\dom H_0^{1/2} \subset \dom V$. Then the operator $VH_0^{-1/2}$ is bounded in $\cH$ and we assume that $S\!:=\!V H_0^{-1/2}\!=\!S_0\!+\!S_1$ with a strict contraction $S_0$ and compact~$S_1$, as in \cite{MR2238908}, \cite{MR2465932}, \cite{MR2268872}.  
The assumption on $V$ implies that $V\!=\!V_0\!+\!V_1$ where $V_0\!:=\!S_0H_0^{1/2}$ is $H_0^{1/2}$-bounded with $H_0^{1/2}$-bound $\le\! \|S_0\|\!<\!1$ and $V_1\!:=\!S_1H_0^{1/2}$ is $H_0^{1/2}$-compact.
Hence, by \cite[Thm.~III.7.6]{MR89b:47001}, it follows that $V$ is $H_0^{1/2}$-bounded with $H_0^{1/2}$-bound $< 1$; more precisely, for every $\epsilon\!>\!0$ with $\epsilon^2  <  1  -  \|S_0\|^2$ there exists $a_\epsilon > 0$ such that
\begin{alignat}{2}
\label{eq:relbd1}
\|V_0x\|^2 \le \|S_0\|^2 \|H_0^{1/2} x\|^2, \quad \|V_1 x\|^2 \le a_\epsilon^2 \|x\|^2 + \epsilon^2 \|H_0^{1/2} x\|^2, &&\qquad x \in \dom H_0^{1/2},
\intertext{and hence}
\label{eq:relbd2}
\|Vx\|^2 \le a_\epsilon^2 \|x\|^2 + b_\epsilon^2 \|H_0^{1/2} x\|^2, \quad b_\epsilon^2:= \|S_0\|^2+\epsilon^2 < 1, 
&&\qquad x \in \dom H_0^{1/2}.
\end{alignat}
As an example in the physical case in $\R^d$ with $d\ge 3$, any $V_1 \in L_p(\R^d)$ with $d \le p < \infty$ is $(-\Delta+m^2)^{1/2}$-compact, see e.g.\ 
\cite[Thm.\ 6.1]{MR2238908}, and also its proof for estimates \eqref{eq:relbd1}. 
In the sequel $W(V) := \{ (Vx,x):x \in  \dom V, \|x\| = 1\}$ is the numerical range of $V$.
\begin{prop}
\label{prop:nr-spec}
	The non-real $($point$)$ spectrum of $T$ is bounded and satisfies
	\[
	  \sigma_p(T_V) \setminus \R \subset 
		\left\{ z\in\C: \re z \in W(V), \, |z|^2 \le \frac{a_\epsilon^2}{1-(\|S_0\|^2 + \epsilon^2)} - m^{2} \right\}
	\]
for every $0 < \epsilon  <  (1 - \|S_0\|^2)^{1/2}$, $a_\epsilon \ge 0$ as in \eqref{eq:relbd1} bounding the relatively compact part~$V_1$ 
and with $S_0=V_0 H_0^{-1/2}$ the strictly contractive part of $V$.
\end{prop}
\begin{proof} Let $\lambda \in \sigma_p(T_V) \setminus \R$. Then there exists an $x\in \dom H_0 \subset \dom V^2 \subset \dom V$, $\|x\|=1$, such that $T_V(\lambda)x=0$ and  \vspace{-1mm} hence
$$
  0 = (T_V(\lambda) x,x) = (H_0x,x) - (V^2x,x) + 2 \lambda (Vx,x) - \lambda^2.
$$
Since $H_0$ is self-adjoint and $V$ is symmetric, all coefficients in the quadratic equation above are real. Taking real and imaginary part, we thus \vspace{-1mm} obtain
\begin{align}
\label{eq:ev-re1}
0 & = (H_0^{1/2}x,H_0^{1/2}x) - (Vx,Vx) + 2 \re \lambda (Vx,x) - ( (\re \lambda)^2 - (\im \lambda)^2), \\
\label{eq:ev-im1}
0 & = 2 \im \lambda (Vx,x) - 2 \re \lambda \im \lambda.
\end{align}
Since $\lambda\notin \R$ by assumption, we can solve \eqref{eq:ev-im1} for $\lambda$ and insert the resulting identity for $(Vx,x)$ into \eqref{eq:ev-re1} to conclude  \vspace{-1mm}that
\begin{align}
\label{eq:ev-im2}
\re \lambda & = (Vx,x)  \in W(V), \\ 
\label{eq:ev-re2}
|\lambda|^2 & = \|Vx\|^2 - \|H_0^{1/2}x\|^2. 
\end{align}
By \eqref{eq:ev-re2} it follows that $0 \le - \|H_0^{1/2}x\|^2 + \|Vx\|^2$ and thus, together with \eqref{eq:relbd2},  
\[
  \|H_0^{1/2}x\|^2 \le \|Vx\|^2 \le a_\epsilon^2 \|x\|^2 + b_\epsilon^2 \|H_0^{1/2} x\|^2 \le  a_\epsilon^2 \|x\|^2 + b_\epsilon^2 \|V x\|^2.
\]
Since $b_\epsilon^2=\|S_0\|^2+ \epsilon^2<1$, this \vspace{-1mm}implies
\[
  \|Vx\|^2 \le \frac {a_\epsilon^2}{1-b_\epsilon^2} \|x\|^2 = \frac {a_\epsilon^2}{1-b_\epsilon^2} .
\]
Using the estimate $\|H_0^{1/2}x\|^2=(H_0x,x) \ge m ^2\|x\|^2 = m^2$ in \eqref{eq:ev-re2}, we finally \vspace{-1mm}obtain
\begin{align*}
	 \|\lambda\|^2 &\le \frac {a_\epsilon^2}{1-b_\epsilon^2} - m^2. 
\qedhere	
\end{align*}
\end{proof}
\begin{example}
	If $V$ is bounded and relatively compact, one can choose $S_0=0$, $\epsilon>0$ and $a_\epsilon=\|V\|$ in \eqref{eq:relbd2}, and hence in this case Proposition \ref{prop:nr-spec} yields
	\begin{equation}\label{eq:a-priori_inclusion}
		\sigma_p(T) \setminus \R \subset \left\{ z\in\C: |z|^2 \le \|V\|^2 - m^2, \, \re z \in W(V) \right\};
	\end{equation}
	if e.g.\ the potential is non-positive $V\le 0$ as in \cite{schiff}, \cite[Ex.\ V.2]{MR2268872}, then
	\begin{equation}\label{eq:a-priori_inclusion2}
	\hspace{7mm} \sigma_p(T) \setminus \R \subset \left\{ z\in\C: |z|^2 \le \|V\|^2 - m^2, \,  -\|V\| \le \re z \le 0 \right\}.
	\end{equation}
\end{example}
\begin{remark}
\begin{enumi}
	\item
	The above estimate for $\sigma_p(T_V) \setminus \R$ holds for any choice of $S_0$ and $S_1$ such that $VH_0^{-1/2} = S_0 + S_1$ and for any pair $0 < \epsilon < (1-\|S_0\|^2)^{1/2}$, $a_\epsilon \ge 0$ as in \eqref{eq:relbd2}. Hence $\sigma_p(T_V)$ is contained in the intersection over all corresponding enclosures obtained by Proposition \ref{prop:nr-spec}. However, for bounded $V$, one can show that the choice $\epsilon=0$, $a_\epsilon=\|V\|$ yields the tightest enclosure.
	\item
	For the real eigenvalues in the essential spectral gap, two-sided eigenvalue estimates were derived in \cite[Sect.~IV]{MR2268872}.
\end{enumi}
\end{remark}
\subsubsection{Bounds in one space dimension}
In view of our numerical examples, we also derive eigenvalue bounds by means of a Birman-Schwinger argument. 
To this end, we adapt an argument from \cite{AD} to obtain explicit bounds on the eigenvalues of $T_V$ in terms of $V$,
see also \cite{MR3177918}, \cite{MR3216811} for another case of non-semi-bounded spectra.

\begin{lemma}\label{lemma:AD_bound}
	Assume that $V\in L^1(\R)\cap L^{\!\infty}(\R)$. Then every $\lambda\in\sigma_p(T_V)$ satisfies
	\begin{align}\label{eq:AD_bound1}
		4|m^2-\lambda^2| \leq \|V^2-2\lambda V\|_{L^1(\R)}^2
	\end{align}
	and, as a consequence,
	\begin{align}\label{eq:AD_bound2}
		2|m^2-\lambda^2|^{\f12} \leq \|V\|_{L^2\R}^2 + 2|\lambda| \|V\|_{L^1(\R)}.
	\end{align}
\end{lemma}
\begin{proof}
	Let $\lambda\in\sigma_p(T_V)$ be an eigenvalue of $T_V$ with corresponding eigenfunction  $u\in W^{2,2}(\R)$. 
	Then
	\begin{align*}
		(-\Delta+m^2-\lambda^2)u = W_{\!\lambda} u
		\ \ \iff \ \ 
		u = (-\Delta+m^2-\lambda^2)^{-1}W_{\!\lambda} u.
	\end{align*}
	Define $g:= \f{W_{\!\lambda}}{|W_{\!\lambda}|^{1/2}}u$. Then $g\in L^2(\R)$ because $V\in L^{\!\infty}(\R)$ and
	\begin{align}\label{eq:g=Ag}
		g &= \f{W_{\!\lambda}}{|W_{\!\lambda}|^{\f12}} (-\Delta+m^2-\lambda^2)^{-1} |W_{\!\lambda}|^{\f12} g.
	\end{align}
	To simplify notation, we introduce $A(\lambda):=\f{W_{\!\lambda}}{|W_{\!\lambda}|^{1/2}} (-\Delta+m^2-\lambda^2)^{-1} |W_{\!\lambda}|^{\f12}$. 
	Equation \eqref{eq:g=Ag} implies that if $\lambda\in\C\setminus\{\pm\sqrt{\sigma(H_0)}\}$ is an eigenvalue, then $-1\in\sigma(A(\lambda))$.
	It follows that $\|A(\lambda)\|_{\mathrm{HS}} \geq 1$ (otherwise $I-A(\lambda)$ would be invertible). Using the fundamental solution for $-\Delta+z^2$ in one dimension, $A(\lambda)$ can be written as an explicit integral operator and its Hilbert-Schmidt norm satisfies 
	\begin{align}
		1\leq \|A\|_{\mathrm{HS}}^2 &\leq \int_{\R^2} \left| \f{W_{\!\lambda}(x)}{|W_{\!\lambda}(x)|^{\f12}} \f{\e^{-\sqrt{m^2-\lambda^2}|x-y|}}{2\sqrt{m^2-\lambda^2}} |W_{\!\lambda}(y)|^{\f12} \right|^2\,dx\,dy
		\label{eq:AD_fund_sol}
		\\
		\nonumber
		\iff \quad 4|m^2-\lambda^2| &\leq \int_{\R^2}  |W_{\!\lambda}(x)| \e^{-2\re\sqrt{m^2-\lambda^2}|x-y|} |W_{\!\lambda}(y)| \,dx\,dy
		\,\leq \,\|V^2-2\lambda V\|_{L^1\R}^2
	\end{align}
	This proves \eqref{eq:AD_bound1}. Equation \eqref{eq:AD_bound2} follows readily from  \eqref{eq:AD_bound1} and the triangle inequality.
\end{proof}%
\subsection{Definition of the algorithm}\label{sec:definition_of_algo}
Having established full control over the $\lambda$-dependence of all constants (see Section \ref{sec:preps}) and the location of the point spectrum (see Section \ref{sec:spectral_bound}), we are finally ready to define the algorithm and prove convergence. 

Let $\rho>0$ be large enough to ensure that $B_\rho(0)$ contains all non-embedded eigenvalues of $T_V$ and denote $B_n:= \{z\in \overline{B_\rho(0)}\,|\, \dist(z^2,[m^2,\infty)) > n^{-\f14}\}$; note that 
by \eqref{eq:a-priori_inclusion} and \eqref{eq:hypothesis_V_bounds} one can choose $\rho = \max\{\sqrt{|M^2-m^2|}, m\}$. 
For $n\!\in\!\N$ let $\cL_n \!:=\! n^{-1}(\Z\!+\!\I\Z)\cap B_n$; note that $\cL_n$ can be constructed using finitely many arithmetic operations. Moreover, choose $A\!>\!0$ as in Proposition \ref{prop:C_K_bound}. Then, for \vspace{-1mm}  any~$z\!\in\! \cL_n$,
\begin{align}\label{eq:Cn^12}
		C_K^\lambda &\leq \f{A}{\dist(z^2,[m^2,\infty))^2}
		\leq An^{\f12}.
\end{align}
Applying Theorem \ref{th:operatororm_convergence}, we conclude that 
\begin{align*}
	\big\|K(z) - K^{(N_n,R_n,s_n)}(z)\big\|_{L^2\to L^2} \leq \f{A}{n^\f12}
\end{align*}
for all $z\!\in\! \cL_n$ and  $n\!\in\!\N$. We emphasise that the constant $A$ is explicit and can be computed in finitely many arithmetic operations from $\rho,m,M$ (this follows from the proof of Proposition~\ref{prop:C_K_bound}).

Similarly, from \eqref{eq:Lipschitz_constant} we conclude that there exists an explicit constant $B>0$ such that
\begin{align}\label{eq:C'n^12}
	C_L &\leq \f{B}{\dist(z^2,[m^2,\infty))^2} \leq Bn^{\f12}.
\end{align}
for all $z\in \cL_n$ and all $n\in\N$.
\begin{de}[Spectral algorithm]\label{de:definition_of_Gamma_n}
	For $n\in\N$ let
	\begin{equation}\label{eq:def_algorithm}
		\Gamma_n(V) := \left\{ z\in \cL_n \,\middle|\, \|(I-K_n(z))^{-1}\|_{L^2\to L^2} \geq \f{n^{\f12}}{2(A+B)} \right\}
	\end{equation}
	where $K_n = K^{(N_n,R_n,s_n)}$ is as in Theorem \ref{th:operatororm_convergence} and $A$, $B$ are as in \eqref{eq:Cn^12}, \eqref{eq:C'n^12}, respectively.
\end{de}
\begin{lemma}
	The sequence $(\Gamma_n)_{n\in\N}$ is a tower of arithmetic algorithms of height $1$ in the sense of Definition {\rm \ref{def:Arithmetic-Tower}}. 
\end{lemma}
\begin{proof}
	For $n\in\N$ and $V\in\Omega_{p,M}$ define the subset $\Lambda_{\Gamma_n}(V)\subset\Lambda$ by $\Lambda_{\Gamma_n}(V) := \{V\mapsto V(x) \,|\, x\in (s_n^{-1}\Z^d)\cap[-R_n,R_n)^d\}$ where $s_n,R_n$ were defined in \eqref{eq:parameter_choices}. Clearly, this is a finite set and $(s_n^{-1}\Z^d)\cap[-R_n,R_n)^d\subset\Q^d$.
	We need to prove that $\Gamma_n(V)$ can be computed in finitely many arithmetic operations from the $f\in\Lambda_{\Gamma_n}(V)$. First we note that the lattice $\cL_n$ and the constants $A,B$ can be computed in finitely many arithmetic operations from $m,M$. This is clear from \eqref{eq:Lipschitz_constant} and from the proof of Proposition \ref{prop:C_K_bound}. Moreover, the matrix elements of $K^{(N_n,R_n,s_n)}$ can be computed in finitely many arithmetic operations from the point values $V(i)$ for $i\in s^{-1}\Z$ using Definition \ref{de:approximate_matrix_elements}. Finally, for any matrix $X$ and any $\eta>0$, testing whether $\|X^{-1}\|>\eta$ can be done in finitely many arithmetic operations by \cite[Prop.~10.1]{AHS}.
\end{proof}
Now the main result of this paper is the following. 
\begin{theorem}\label{th:convergence_of_Gamma_n}
	For any $V:\R^d\to\R$ satisfying Hypothesis {\rm \ref{hyp:assumptions_on_V}} we have
	\begin{align}
	\label{eq:conv}
		d_{\textnormal{H}}\Bigl(\Gamma_n(V)\big\backslash\{\pm\sqrt{\sigma(H_0)}\} ,\, \sigma_p(T_V)\big\backslash\{\pm\sqrt{\sigma(H_0)}\}\Bigr) \to 0
	\end{align}
	as $n\to\infty$. Moreover, the explicit error \vspace{-1.5mm} bound
	\begin{align}\label{eq:pi_1_bound}
		\sup_{z\in\sigma(T_V)\setminus\{\pm\sqrt{\sigma(H_0)}\}} \dist(z,\Gamma_n(V)) \leq \f1n
\\[-8mm] \nonumber
	\end{align}
	holds for all $n\in\N$.
\end{theorem}
\begin{proof}
	For notational convenience, we will not distinguish between a bounded operator $O$ on $\cH_n$ and its extension $O\oplus0$ on $\cH$. This is justified because, clearly, $\|O\|_{\cH_n\to\cH_n} = \|O\oplus0\|_{\cH\to\cH}$. Moreover, the notation $\lim_{n\to\infty}E_n$ for $E_n\subset\C$ refers to limits in the Hausdorff distance $d_{\textnormal{H}}$.
	
	The proof of \eqref{eq:conv} is in two steps. In the first step we exclude spectral pollution (i.e.\ we show $\lim_{n\to\infty}\Gamma_n(V)\!\subset\!\sigma(T_V)$); in the second step we prove spectral inclusion (i.e.\ we show $\lim_{n\to\infty}\Gamma_n(V)\!\supset\!\sigma(T_V)$).
	
	\underline{Step 1: Excluding spectral pollution.}
	Let $z_0\!\in\! \lim_{n\to\infty}\Gamma_n(V)$, i.e.\ suppose that there exist $z_n\!\in\!\Gamma_n(V)$ such that $z_n\!\to\! z_0$ as $n\!\to\!\infty$. By definition of $\Gamma_n$ we have $\|(I-K_n(z_n))^{-1}\|_{L^2\to L^2} \geq n^\f12/(2(A+B))$ and, applying \cite[Lemma 3.1]{ben2020computing}, we conclude that
	\begin{align*}
		\left(1-An^{-\f12}\|(I-K(z_n))^{-1}\|_{L^2\to L^2}\right)\|(I-K_n(z_n))^{-1}\|_{L^2\to L^2} &\leq \|(I-K(z_n))^{-1}\|_{L^2\to L^2}
		\\[-1mm]
		\implies \qquad\qquad\qquad 
		\left(1-An^{-\f12}\|(I-K(z_n))^{-1}\|_{L^2\to L^2}\right) \f{n^\f12}{2(A+B)} &\leq \|(I-K(z_n))^{-1}\|_{L^2\to L^2}.
	\end{align*}
	Rearranging terms, we arrive at
	\begin{align*}
		 \f{n^\f12}{2(A\!+\!B)} &\leq \|(I-K(z_n))^{-1}\|_{L^2\to L^2} + \f{A}{2(A\!+\!B)}\|(I-K(z_n))^{-1}\|_{L^2\to L^2}
		 \\[-1mm]
		&= \left( 1 + \f{A}{2(A\!+\!B)} \right)\|(I-K(z_n))^{-1}\|_{L^2\to L^2}
		 \\[-1mm]
		 \implies \qquad \|(I-K(z_n))^{-1}\|_{L^2\to L^2} &\geq \f{n^\f12}{3A+2B}.
	\end{align*}%
	Thus $\|(I\!-\!K(z_n))^{-1}\|_{L^2\to L^2} \!\to\! \infty$ as $n\!\to\!\infty$ and thus $1\!\in\!\sigma(K(z_0))$, or equivalently, $z_0\!\in\!\sigma(T_V)$ by \eqref{eq:TVK}.
	
	\underline{Step 2: Proving spectral inclusion.}
	Let $z\in\C\setminus\{\pm\sqrt{\sigma(H_0)}\}$ be such that $1\in\sigma(K(z))$ and, for $n\in\N$, let $z_n\in\cL_n$ be an arbitrary point with $|z-z_n|<n^{-1}$. We show that necessarily $z_n\in\Gamma_n(V)$. By the Lipschitz continuity of $K$ (see Lemma \ref{lemma:Lipschitz_continuity}), \cite[Lemma 3.4]{ben2020computing} and \eqref{eq:C'n^12} we \vspace{-1.5mm} have
	\begin{align*}
		\|(I-K(z_n))^{-1}\|_{L^2\to L^2} \geq \f{1}{C_L|z-z_n|}
		\geq \f{1}{C_L n^{-1}}
		\geq \f{n^{\f12}}{B}
	\end{align*}
	for $n\in\N$.	Similarly to Step 1, this implies a lower bound for $\|(I-K_n(z_n))^{-1}\|_{L^2\to L^2}$. Applying \cite[Lemma 3.1]{ben2020computing} (with $K$ and $K_n$ swapped) we obtain
	\begin{align*}
		\left(1 - An^{-\f12}\|(I-K_n(z_n))^{-1}\|_{L^2\to L^2} \right)\|(I-K(z_n))^{-1}\|_{L^2\to L^2} &\leq \|(I-K_n(z_n))^{-1}\|_{L^2\to L^2}
		\\[-1mm]
		\implies\qquad\qquad\qquad\qquad
		\left(1 - An^{-\f12}\|(I-K_n(z_n))^{-1}\|_{L^2\to L^2} \right)\f{n^{\f12}}{B} &\leq \|(I-K_n(z_n))^{-1}\|_{L^2\to L^2}.
	\end{align*}
	Rearranging terms, we conclude \vspace{-1.5mm}that 
	\begin{align}
		\f{n^{\f12}}{B} &\leq \left(1 + \f{A}{B}\right)\|(I-K_n(z_n))^{-1}\|_{L^2\to L^2}
		\nonumber
		\\[-1mm]
		\iff\qquad
		\f{n^{\f12}}{A+B} &\leq \|(I-K_n(z_n))^{-1}\|_{L^2\to L^2}.
		\label{eq:spectral_inclusion}
	\end{align}
	If we compare \eqref{eq:spectral_inclusion} to \eqref{eq:def_algorithm}, we find that $z_n\in \Gamma_n(V)$ for $n\in\N$. 
	By construction we have 
	$|z-z_n|<n^{-1}$ for $n\in\N$ and hence \eqref{eq:pi_1_bound} follows. In particular, \eqref{eq:pi_1_bound} implies that $z\in\lim_{n\to\infty} \Gamma_n(V)$. 
\end{proof}
\section{Numerical Results and the Schiff-Snyder Weinberg Effect}
In this section we present a number of numerical results obtained by implementing our algorithm in MATLAB for potentials from the physical literature.
Therefore we focus on one space dimension. Another motivation is that, for the case $d=1$, there are some simplifications of the algorithm and possibilities to improve convergence which we would like to present.
\subsection{Details of the implementation}
In the case $d=1$ the algorithm \eqref{eq:def_algorithm} can be simplified somewhat while retaining guaranteed convergence. For the remainder of this section, we make an assumption that is satisfied in many physically interesting examples.
\begin{hyp}\label{hyp:additional_assumption_on_V}
	Suppose that $V\in L^1(\R)$ with $\|V\|_{L^1}\leq M$.
\end{hyp}
Under this assumption on $V$, the results in Section \ref{sec:improvements} (see Remark \ref{rem:improved})  guarantee the convergence $\|K(\lambda) - K^{(N_n,R_n,s_n)}(\lambda)\|\to 0$ with rate $n^{- \f 12}$, provided the parameters $r,R,S$ are chosen appropriately. If $d=1$ and e.g.\ $p=4/3$, the choices
\begin{align}\label{eq:improved_params_1d}
	r_n = n^{\f1{2}}, \quad
	R_n = n^{2}, \quad
	s_n 
	= n^2.
\end{align}
achieve a rate of $n^{-\f12}$, as a comparison with Proposition \ref{prop:matrix_error_bound} and Remark \ref{rem:improved} shows.
We have implemented the algorithm \eqref{eq:def_algorithm} in MATLAB in one space dimension $d=1$, following the scheme \eqref{eq:improved_params_1d}. In order to improve numerical performance further, we added several practical augmentations:
\begin{enumi}
	\item 
	Formula \eqref{eq:E_kN1} for $d=N=1$ reads
	\begin{align*}
		E^\lambda_{k,1}(x) &:= (2\pi)^{-\f 12} n^{\f 12} \int_{Q_k} a_\lambda(k) \e^{-\I\xi x}\,d\xi.
	\end{align*}
	As a low-cost way to improve numerical accuracy, we replaced $a_\lambda(k)$ by its first-order Taylor approximation $a(k)+ a_\lambda'(k)(\xi-k)$. This reduces the approximation error without adding any significant computational cost since the integral over $Q_k$ can still be evaluated explicitly.%
	\item Similarly, the accuracy of $\operatorname{Quad}_{R,s}(\overline{E^\lambda_{k,n}} E_{m,n}^0 W_{\!\lambda})$ has been improved by using Simpson's quadrature formula, rather than a simple piecewise constant Riemann sum (this point was touched upon in the discussion after Definition \ref{de:approximate_matrix_elements}).
	\item Instead of computing all points $z\in\cL_n$ in the set \eqref{eq:def_algorithm}, for an $n^{-1}$-fine grid $\cL_n$, we compute the norms $\|(I-K_n(\lambda))^{-1}\|$ on a fixed grid, determine their local minima and use these as starting points for a further gradient descent minimisation. In most applications this procedure produces a more accurate approximation of $\sigma_p(T_V)$ than \eqref{eq:def_algorithm} in reasonable computation time.
\end{enumi}
The next subsection contains some example results computed with  $r_n = \sqrt{n}/2$, $s=n$ and fixed $R$, large enough to ensure $|V|<10^{-6}$ outside $[-R,R]$. 
The operator norm $\|(I-K_n(\lambda))^{-1}\|$ in \eqref{eq:def_algorithm} is computed via singular value decomposition using MATLAB's built-in \texttt{svds()} command.
The MATLAB code that produced the figures below is openly available at \url{https://github.com/frank-roesler/spectral_klein_gordon}.
\subsection{Examples}\label{sec:numerics}
The physics community has seen sustained interest in the spectral behaviour of the one-dimensional Klein-Gordon equation for many years (see e.g.\ \cite{schiff, villalba2006bound, lv2013noncompeting, lv2016bosonic, jiang2021analysis} and the references therein). Specific interest has focused on the so-called \emph{Schiff-Snyder Weinberg effect} -- the fact that under a continuous change of the potential two real eigenvalues of $T_V$ can join up and become a complex conjugate pair. It has been noticed early on that this can happen even for real-valued bounded potentials~\cite{schiff}. 
It is also known that the first pair of real eigenvalues emerging into the spectral gap are simple eigenvalues with strictly positive eigenfunctions (see \cite[Thm.\ 6.1]{MR2819232}. 
In this section we present a range of examples to show that our algorithm yields meaningful results in physically relevant situations. 
\subsubsection{Fixed width Sauter potential}
We begin by studying $\sigma_p(T_V)$ for $m=1$ and
\begin{align}\label{eq:sauter_potential}
	V(x) = -\f{v_0}{2}\bigg(\tanh\bigg(\f{x+D/2}{W}\bigg) - \tanh\bigg(\f{x-D/2}{W}\bigg)\bigg), \quad x\in\R,
\end{align}
where $v_0>0$, $D,W\in\R$ are parameters.
This function can be thought of as a smoothed version of the square well potential and has become known as the \emph{Sauter potential} \cite{sauter1931verhalten, lv2013noncompeting}.

Since $V$ in \eqref{eq:sauter_potential} is smooth and decays exponentially as $x\to\pm\infty$, it satisfies Hypothesis \ref{hyp:assumptions_on_V}.

\begin{figure}[htbp]
	\centering
%
%
\setlength{\figurewidth}{0.7\textwidth}
\setlength{\figureheight}{.3\figurewidth}

\definecolor{mycolor1}{rgb}{0.149039,0.508051,0.55725}%
\begin{tikzpicture}

\begin{axis}[%
width=\figurewidth,
height=\figureheight,
at={(0,0)},
scale only axis,
xmin=-4,
xmax=4,
xlabel style={font=\color{white!15!black}},
xlabel={$x$},
ymin=-4,
ymax=0.5,
ylabel style={font=\color{white!15!black}},
ylabel={$V$},
axis background/.style={fill=white}
]
\addplot [color=mycolor1, line width=1.0pt, forget plot]
  table[row sep=crcr]{%
-4	-4.1638009928846e-07\\
-2.84	-0.000950435027130325\\
-2.67	-0.002950448014861\\
-2.56	-0.00613756399330523\\
-2.48	-0.0104499075994013\\
-2.42	-0.0155678056260342\\
-2.37	-0.021690518114597\\
-2.33	-0.0282686080141064\\
-2.29	-0.0368216668873389\\
-2.26	-0.0448752094171336\\
-2.23	-0.054663917234965\\
-2.2	-0.0665489768229195\\
-2.17	-0.0809607024187118\\
-2.15	-0.0922203785429367\\
-2.13	-0.105000569895788\\
-2.11	-0.119493219317178\\
-2.09	-0.135910585409852\\
-2.07	-0.154486238357461\\
-2.05	-0.175475730656912\\
-2.03	-0.19915680213179\\
-2.01	-0.225828948839399\\
-1.99	-0.255812155121074\\
-1.97	-0.289444559800995\\
-1.95	-0.327078805471177\\
-1.93	-0.369076809520092\\
-1.91	-0.415802704308507\\
-1.89	-0.467613730333833\\
-1.87	-0.524848939799511\\
-1.85	-0.587815687679696\\
-1.83	-0.656774059600829\\
-1.81	-0.731919611837523\\
-1.79	-0.813365071074095\\
-1.77	-0.901121941200284\\
-1.75	-0.995083258329448\\
-1.72	-1.14709441892457\\
-1.69	-1.3110716658613\\
-1.66	-1.48485565623641\\
-1.61	-1.78835616081256\\
-1.54	-2.21514433941711\\
-1.51	-2.38892832937251\\
-1.48	-2.55290557569825\\
-1.45	-2.70491673546656\\
-1.42	-2.84354169227206\\
-1.4	-2.92814844126737\\
-1.38	-3.00643575720175\\
-1.36	-3.07846801503863\\
-1.34	-3.14440167535051\\
-1.32	-3.20446680020115\\
-1.3	-3.25894917366422\\
-1.28	-3.30817382275453\\
-1.26	-3.35249044533531\\
-1.24	-3.39226100081301\\
-1.22	-3.42784952002107\\
-1.2	-3.45961404529241\\
-1.18	-3.48790051611613\\
-1.16	-3.51303836195582\\
-1.14	-3.53533754236877\\
-1.12	-3.55508677633515\\
-1.1	-3.57255271968409\\
-1.08	-3.58797987525998\\
-1.06	-3.60159105019354\\
-1.03	-3.6190392076752\\
-1	-3.63345091338365\\
-0.97	-3.64533594867786\\
-0.94	-3.65512462682056\\
-0.91	-3.66317813310287\\
-0.87	-3.67173113086299\\
-0.83	-3.67830914097035\\
-0.78	-3.68443171859675\\
-0.72	-3.68954938263142\\
-0.65	-3.69343893616103\\
-0.56	-3.69639488843227\\
-0.43	-3.69847969524536\\
-0.22	-3.69960629626187\\
0.21	-3.69962903040408\\
0.46	-3.69814523900984\\
0.59	-3.69559848584381\\
0.68	-3.69198989818802\\
0.75	-3.687243849079\\
0.8	-3.68222199793741\\
0.85	-3.67523615323076\\
0.89	-3.66773434931208\\
0.93	-3.65798588981457\\
0.96	-3.64881260794159\\
0.99	-3.63767054889134\\
1.02	-3.62415356775061\\
1.05	-3.60777954276163\\
1.07	-3.59499936121832\\
1.09	-3.58050672038003\\
1.11	-3.56408936179681\\
1.13	-3.54551371541871\\
1.15	-3.5245242288658\\
1.17	-3.5008431624168\\
1.19	-3.47417102010389\\
1.21	-3.44418781766397\\
1.23	-3.41055541634127\\
1.25	-3.37292117360354\\
1.27	-3.33092317211452\\
1.29	-3.28419727955902\\
1.31	-3.23238625547939\\
1.33	-3.17515104770682\\
1.35	-3.1121843012973\\
1.37	-3.04322593065058\\
1.39	-2.96808037951474\\
1.41	-2.88663492122505\\
1.43	-2.79887805190864\\
1.45	-2.70491673546656\\
1.48	-2.55290557569825\\
1.51	-2.38892832937251\\
1.54	-2.21514433941711\\
1.59	-1.91164383515797\\
1.66	-1.48485565623641\\
1.69	-1.3110716658613\\
1.72	-1.14709441892457\\
1.75	-0.995083258329448\\
1.78	-0.856458300448155\\
1.8	-0.771851550576467\\
1.82	-0.693564233620489\\
1.84	-0.621531974598614\\
1.86	-0.555598312917233\\
1.88	-0.49553318648821\\
1.9	-0.441050811209776\\
1.92	-0.391826160034797\\
1.94	-0.347509535062926\\
1.96	-0.307738976845146\\
1.98	-0.272150454499227\\
2	-0.240385925636377\\
2.02	-0.21209945070357\\
2.04	-0.186961600164045\\
2.06	-0.164662414376831\\
2.08	-0.144913174266085\\
2.1	-0.127447223893268\\
2.12	-0.112020060288955\\
2.14	-0.0984088761794757\\
2.17	-0.0809607024187118\\
2.2	-0.0665489768229195\\
2.23	-0.054663917234965\\
2.26	-0.0448752094171336\\
2.29	-0.0368216668873389\\
2.33	-0.0282686080141064\\
2.37	-0.021690518114597\\
2.42	-0.0155678056260342\\
2.48	-0.0104499075994013\\
2.55	-0.00655993199716054\\
2.64	-0.00360304922255139\\
2.77	-0.00151539850962479\\
2.98	-0.000373808017255506\\
3.52	-1.02148303104954e-05\\
4	-4.1638009928846e-07\\
};
\end{axis}
\end{tikzpicture}%
	\vspace{-4mm}
	\caption{Sauter potential \eqref{eq:sauter_potential} for $v_0=3.7$, $D=3.2$, \vspace{-1mm}$W=0.3$.}
	\label{fig:sauter_potential}
\end{figure}
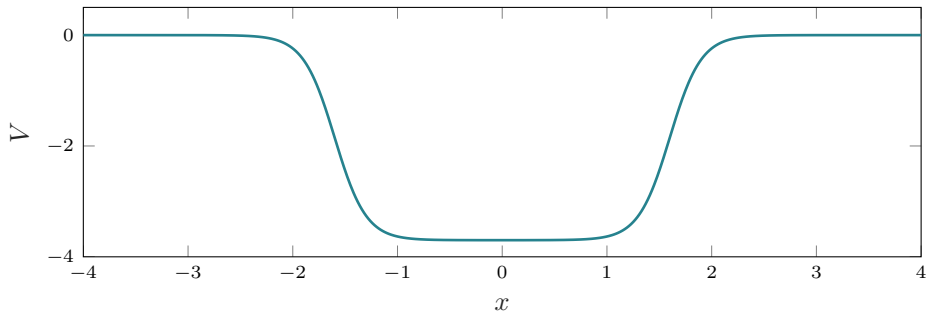

\begin{figure}[htbp]
	\centering
	\includegraphics{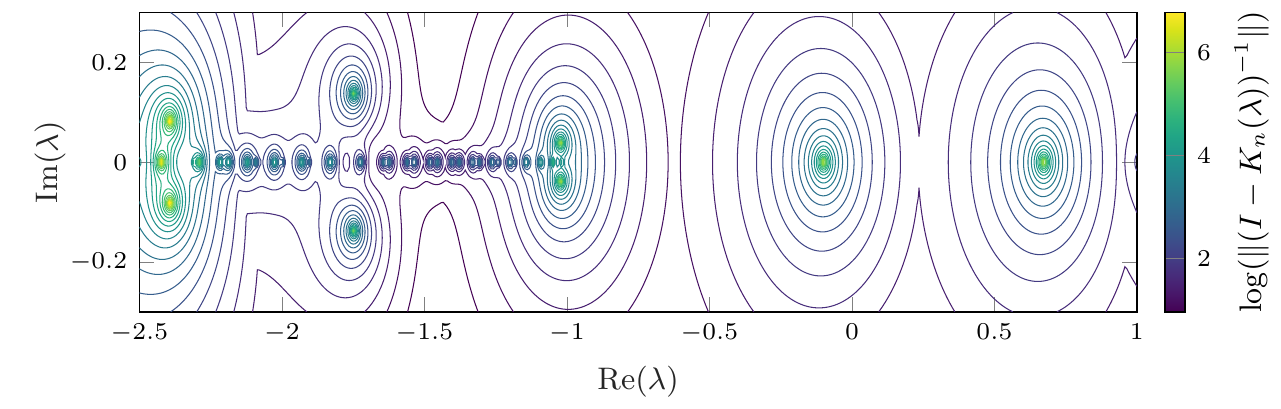}
	\vspace{-4mm}
	\caption{Logarithmic contour plot of $\|(I-K_n(\lambda))^{-1}\|_{L^2\to L^2}$ for Sauter potential with parameters $v_0=3.7$, $D=3.2$, $W=0.3$.}
	\label{fig:sauter_contour}
\end{figure}

As reported in \cite{lv2013noncompeting}, for the potential in Figure \ref{fig:sauter_potential} with parameter values $v_0=3.7$, $D=3.2$, $W=0.3$ the Klein-Gordon operator pencil $T_V$ has three complex conjugate pairs $E_1^\pm$, $E_2^\pm$, $E_3^\pm$ of eigenvalues near the negative real axis, as well as two real eigenvalues $E_4$, $E_5$ in the gap $[-1,1]$ (from left to right). Our method confirms these results: Figure \ref{fig:sauter_contour} shows a contour plot of $\|(I-K_n(\lambda))^{-1}\|$ for $n=200$. Indeed, visual inspection suggests two poles of $\|(I-K_n(\lambda))^{-1}\|$ in the gap $[-m,m]=[-1,1]$ and three pairs of complex conjugate poles below Re$(\lambda)=-1$, which are our candidates for $E_1^\pm$, $E_2^\pm$, $E_3^\pm$. 

\begin{figure}[htbp]
	\centering
%
%
\definecolor{mycolor1}{rgb}{0.00000,0.44700,0.74100}%
\setlength{\figurewidth}{0.3\textwidth}
\setlength{\figureheight}{0.8\figurewidth}

\begin{tikzpicture}

\begin{axis}[%
width=\figurewidth,
height=\figureheight,
at={(0,0)},
scale only axis,
xmin=-2.5,
xmax=-0.9,
xlabel style={font=\color{white!15!black}},
xlabel={\small Re$(\lambda)$},
ymin=-0.19,
ymax=0.19,
ylabel shift = -2mm,
ylabel style={font=\color{white!15!black}},
ylabel={\small Im$(\lambda)$},
axis background/.style={fill=white},
]
\addplot [color=mycolor1, draw=none, mark size=1.5pt, mark=*, mark options={solid, black!20!red}, forget plot]
  table[row sep=crcr]{%
-2.39254314974245	-0.0828812716562908\\
-2.39254314357947	0.0828812720294398\\
-1.74732971351311	-0.137515056766464\\
-1.74732971442845	0.137515056688179\\
-1.02262319780288	-0.0383565857654777\\
-1.02262319757999	0.0383565861048378\\
-0.0992799283624244	1.92850180269488e-10\\
0.673233912866297	-6.76033007351862e-10\\
};
\addplot [forget plot, line width=0.3]
  table[row sep=crcr]{%
-5 0\\
3 0\\
};
\end{axis}

\begin{axis}[%
width=\figurewidth,
height=\figureheight,
at={(\figurewidth+0.1\textwidth,0)},
scale only axis,
xmin=-4,
xmax= 2,
xlabel style={font=\color{white!15!black}},
xlabel={\small Re$(\lambda)$},
ymin=-4,
ymax= 4,
ylabel shift = -2mm,
ylabel style={font=\color{white!15!black}},
ylabel={\small Im$(\lambda)$},
axis background/.style={fill=white}
]
\addplot [color=mycolor1, draw=none, mark size=0.7pt, mark=*, mark options={solid, black!20!red}, forget plot]
  table[row sep=crcr]{%
-2.39254314974245	-0.0828812716562908\\
-2.39254314357947	0.0828812720294398\\
-1.74732971351311	-0.137515056766464\\
-1.74732971442845	0.137515056688179\\
-1.02262319780288	-0.0383565857654777\\
-1.02262319757999	0.0383565861048378\\
};

\addplot[area legend, draw=black, fill=black, draw opacity=0, fill opacity=0.1, forget plot]
table[row sep=crcr] {%
x	y\\
-3.56212347527751	0\\
-3.52614243007268	0.505018035286082\\
-3.49016138486786	0.712388349638287\\
-3.45418033966304	0.870265381483401\\
-3.41819929445821	1.00231593646338\\
-3.38221824925339	1.11773134855404\\
-3.34623720404857	1.22123717080849\\
-3.31025615884374	1.31564729922581\\
-3.27427511363892	1.40279939169075\\
-3.2382940684341	1.48397276911258\\
-3.20231302322927	1.56010094365054\\
-3.16633197802445	1.63189014889571\\
-3.13035093281963	1.69989019954772\\
-3.0943698876148	1.76453921796764\\
-3.05838884240998	1.82619312827128\\
-3.02240779720516	1.88514581943057\\
-2.98642675200033	1.94164335243624\\
-2.95044570679551	1.99589423176055\\
-2.91446466159069	2.04807699793298\\
-2.87848361638586	2.09834594938996\\
-2.84250257118104	2.14683552838877\\
-2.80652152597622	2.19366373388339\\
-2.77054048077139	2.2389348131489\\
-2.73455943556657	2.2827414103391\\
-2.69857839036175	2.32516630033116\\
-2.66259734515692	2.36628380180535\\
-2.6266162999521	2.40616093932826\\
-2.59063525474728	2.44485840693968\\
-2.55465420954245	2.4824313732287\\
-2.51867316433763	2.51893015869209\\
-2.48269211913281	2.55440080933258\\
-2.44671107392798	2.58888558531297\\
-2.41073002872316	2.62242337957385\\
-2.37474898351834	2.65505007832277\\
-2.33876793831351	2.68679887297875\\
-2.30278689310869	2.71770053134077\\
-2.26680584790387	2.74778363431911\\
-2.23082480269905	2.77707478343416\\
-2.19484375749422	2.8055987833815\\
-2.1588627122894	2.83337880323291\\
-2.12288166708457	2.86043651925354\\
-2.08690062187975	2.88679224183539\\
-2.05091957667493	2.91246502865436\\
-2.01493853147011	2.93747278583482\\
-1.97895748626528	2.96183235863843\\
-1.94297644106046	2.9855596129716\\
-1.90699539585564	3.00866950882089\\
-1.87101435065081	3.03117616656997\\
-1.83503330544599	3.05309292702123\\
-1.79905226024117	3.07443240583433\\
-1.76307121503634	3.09520654300054\\
-1.72709016983152	3.11542664789182\\
-1.6911091246267	3.13510344035527\\
-1.65512807942187	3.15424708826566\\
-1.61914703421705	3.17286724189797\\
-1.58316598901223	3.19097306543914\\
-1.5471849438074	3.20857326592035\\
-1.51120389860258	3.22567611981914\\
-1.47522285339776	3.24228949755201\\
-1.43924180819293	3.25842088605396\\
-1.40326076298811	3.27407740961956\\
-1.36727971778329	3.28926584916172\\
-1.33129867257846	3.30399266002725\\
-1.29531762737364	3.31826398849432\\
-1.25933658216882	3.33208568706366\\
-1.22335553696399	3.34546332864413\\
-1.18737449175917	3.35840221972336\\
-1.15139344655435	3.37090741260492\\
-1.11541240134952	3.38298371678594\\
-1.0794313561447	3.39463570954156\\
-1.04345031093988	3.40586774577678\\
-1.00746926573505	3.41668396720012\\
-0.97148822053023	3.42708831086888\\
-0.935507175325406	3.43708451715081\\
-0.899526130120583	3.44667613714335\\
-0.863545084915759	3.45586653958757\\
-0.827564039710936	3.46465891731068\\
-0.791582994506113	3.47305629322818\\
-0.755601949301289	3.4810615259336\\
-0.719620904096466	3.48867731490181\\
-0.683639858891643	3.49590620532899\\
-0.64765881368682	3.50275059263101\\
-0.611677768481996	3.50921272661946\\
-0.575696723277173	3.51529471537324\\
-0.53971567807235	3.5209985288219\\
-0.503734632867526	3.52632600205553\\
-0.467753587662703	3.53127883837453\\
-0.43177254245788	3.53585861209163\\
-0.395791497253056	3.54006677109702\\
-0.359810452048233	3.54390463919671\\
-0.32382940684341	3.54737341823307\\
-0.287848361638587	3.55047418999562\\
-0.251867316433763	3.55320791792931\\
-0.21588627122894	3.55557544864653\\
-0.179905226024117	3.55757751324863\\
-0.143924180819293	3.55921472846169\\
-0.10794313561447	3.56048759759065\\
-0.0719620904096465	3.56139651129539\\
-0.0359810452048235	3.56194174819144\\
0	3.56212347527751\\
0	-3.56212347527751\\
-0.0359810452048235	-3.56194174819144\\
-0.0719620904096465	-3.56139651129539\\
-0.10794313561447	-3.56048759759065\\
-0.143924180819293	-3.55921472846169\\
-0.179905226024117	-3.55757751324863\\
-0.21588627122894	-3.55557544864653\\
-0.251867316433763	-3.55320791792931\\
-0.287848361638587	-3.55047418999562\\
-0.32382940684341	-3.54737341823307\\
-0.359810452048233	-3.54390463919671\\
-0.395791497253056	-3.54006677109702\\
-0.43177254245788	-3.53585861209163\\
-0.467753587662703	-3.53127883837453\\
-0.503734632867526	-3.52632600205553\\
-0.53971567807235	-3.5209985288219\\
-0.575696723277173	-3.51529471537324\\
-0.611677768481996	-3.50921272661946\\
-0.64765881368682	-3.50275059263101\\
-0.683639858891643	-3.49590620532899\\
-0.719620904096466	-3.48867731490181\\
-0.755601949301289	-3.4810615259336\\
-0.791582994506113	-3.47305629322818\\
-0.827564039710936	-3.46465891731068\\
-0.863545084915759	-3.45586653958757\\
-0.899526130120583	-3.44667613714335\\
-0.935507175325406	-3.43708451715081\\
-0.97148822053023	-3.42708831086888\\
-1.00746926573505	-3.41668396720012\\
-1.04345031093988	-3.40586774577678\\
-1.0794313561447	-3.39463570954156\\
-1.11541240134952	-3.38298371678594\\
-1.15139344655435	-3.37090741260492\\
-1.18737449175917	-3.35840221972336\\
-1.22335553696399	-3.34546332864413\\
-1.25933658216882	-3.33208568706366\\
-1.29531762737364	-3.31826398849432\\
-1.33129867257846	-3.30399266002725\\
-1.36727971778329	-3.28926584916172\\
-1.40326076298811	-3.27407740961956\\
-1.43924180819293	-3.25842088605396\\
-1.47522285339776	-3.24228949755201\\
-1.51120389860258	-3.22567611981914\\
-1.5471849438074	-3.20857326592035\\
-1.58316598901223	-3.19097306543914\\
-1.61914703421705	-3.17286724189797\\
-1.65512807942187	-3.15424708826566\\
-1.6911091246267	-3.13510344035527\\
-1.72709016983152	-3.11542664789182\\
-1.76307121503634	-3.09520654300054\\
-1.79905226024117	-3.07443240583433\\
-1.83503330544599	-3.05309292702123\\
-1.87101435065081	-3.03117616656997\\
-1.90699539585564	-3.00866950882089\\
-1.94297644106046	-2.9855596129716\\
-1.97895748626528	-2.96183235863843\\
-2.01493853147011	-2.93747278583482\\
-2.05091957667493	-2.91246502865436\\
-2.08690062187975	-2.88679224183539\\
-2.12288166708457	-2.86043651925354\\
-2.1588627122894	-2.83337880323291\\
-2.19484375749422	-2.8055987833815\\
-2.23082480269905	-2.77707478343416\\
-2.26680584790387	-2.74778363431911\\
-2.30278689310869	-2.71770053134077\\
-2.33876793831351	-2.68679887297875\\
-2.37474898351834	-2.65505007832277\\
-2.41073002872316	-2.62242337957385\\
-2.44671107392798	-2.58888558531297\\
-2.48269211913281	-2.55440080933258\\
-2.51867316433763	-2.51893015869209\\
-2.55465420954245	-2.4824313732287\\
-2.59063525474728	-2.44485840693968\\
-2.6266162999521	-2.40616093932826\\
-2.66259734515692	-2.36628380180535\\
-2.69857839036175	-2.32516630033116\\
-2.73455943556657	-2.2827414103391\\
-2.77054048077139	-2.2389348131489\\
-2.80652152597622	-2.19366373388339\\
-2.84250257118104	-2.14683552838877\\
-2.87848361638586	-2.09834594938996\\
-2.91446466159069	-2.04807699793298\\
-2.95044570679551	-1.99589423176055\\
-2.98642675200033	-1.94164335243624\\
-3.02240779720516	-1.88514581943057\\
-3.05838884240998	-1.82619312827128\\
-3.0943698876148	-1.76453921796764\\
-3.13035093281963	-1.69989019954772\\
-3.16633197802445	-1.63189014889571\\
-3.20231302322927	-1.56010094365054\\
-3.2382940684341	-1.48397276911258\\
-3.27427511363892	-1.40279939169075\\
-3.31025615884374	-1.31564729922581\\
-3.34623720404857	-1.22123717080849\\
-3.38221824925339	-1.11773134855404\\
-3.41819929445821	-1.00231593646338\\
-3.45418033966304	-0.870265381483401\\
-3.49016138486786	-0.712388349638287\\
-3.52614243007268	-0.505018035286082\\
-3.56212347527751	-0\\
}--cycle;
\addplot [forget plot, line width=0.3]
  table[row sep=crcr]{%
-5 0\\
3 0\\
};
\end{axis}

\end{tikzpicture}%
	\vspace{-4mm}
	\caption{Left: approximations of the non-real eigenvalues $E_1^\pm$, $E_2^\pm$, $E_3^\pm$ (red dots, from left to right) for Sauter potential with parameters $v_0=3.7$, $D=3.2$, $W=0.3$. Right: same data, zoomed-out, together with the region of inclusion \vspace{-2mm} \eqref{eq:a-priori_inclusion2}.}
	\label{fig:sauter_spectrum}
\end{figure}
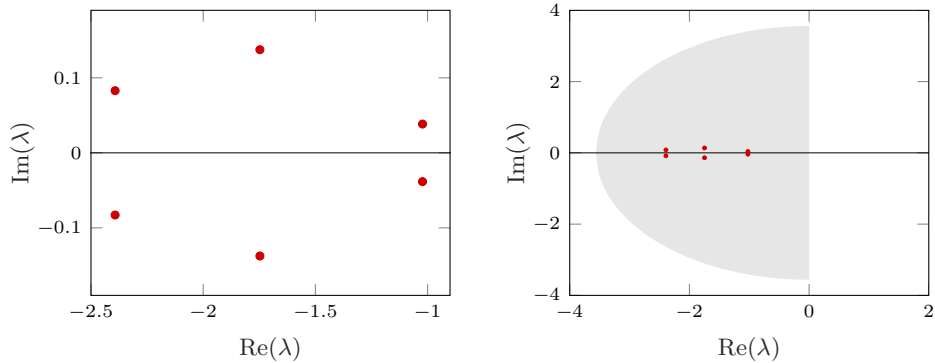

The approximations for the non-real eigenvalues determined from $\|(I-K_n(\lambda))^{-1}\|$ and the region of inclusion \eqref{eq:a-priori_inclusion2} are displayed in Figure~\ref{fig:sauter_spectrum}; the first digits of their numerical values are shown in Table \ref{table:sauter}. The numerically computed values reported in \cite{lv2013noncompeting} agree reasonably well with our eigenvalue approximations which are accompanied by rigorous and explicit error bounds (see Section~\ref{sect:op_norm_bounds}).
\begin{remark}
	Recall that formula \eqref{eq:L_factorization} on which our algorithm is based is only valid for points $\lambda^2\!\notin\!\sigma(H_0)$. This fact is represented by our choice of the lattice $\cL_n$ (cf.\ Section \ref{sec:definition_of_algo}) which explicitly excludes the set $-\sqrt{\sigma(H_0)}\cup \sqrt{\sigma(H_0)}\!=\!(-\infty,-m]\cup[m,\infty)\!=\!(-\infty,-1]\cup[1,\infty)$. Figure~\ref{fig:sauter_contour} shows the numerical~side~of~this: there are many poles on the interval $[-2.5,-1]$ which have to be excluded because they lie \vspace{-2mm} in~$\pm\sqrt{\sigma(H_0)}$.
\end{remark}

\vspace{-2mm} 

\begin{table}[H]
	\begin{tabular}{c r}\toprule
		\textbf{Eigenvalue}\; & \textbf{Approximate value} \\
		\hline\\[-3mm]
		$E_1^\pm$ & $-2.3925 \pm 0.0829\I$\\
		$E_2^\pm$ & $-1.7473 \pm 0.1375\I$\\
		$E_3^\pm$ & $-1.0226 \pm 0.0384\I$\\
		\bottomrule\\
	\end{tabular}
	\caption{Approximate numerical values of the non-real eigenvalues $E_1^\pm$, $E_2^\pm$, $E_3^\pm$ of $T_V$ for the Sauter potential \eqref{eq:sauter_potential} with parameters \vspace{-2mm} $v_0\!=\!3.7$, $D\!=\!3.2$,~$W\!=\!0.3$.}
	\label{table:sauter}
\end{table}
\subsubsection{Varying width Sauter potential}

The Schiff-Snyder-Weinberg effect, and especially its onset, can also be observed if the width $D$ of the Sauter potential \eqref{eq:sauter_potential} is varied, see \cite{jiang2021analysis}. In this subsection we focus on selected results from there, cf.\ \cite[Fig.\ 3(a)]{jiang2021analysis} and chose the range of parameter values accordingly to be $v_0=2.5$, $W=0.1$, $D\in[0,6]$. 

\begin{figure}[htbp]
	\centering
	\includegraphics{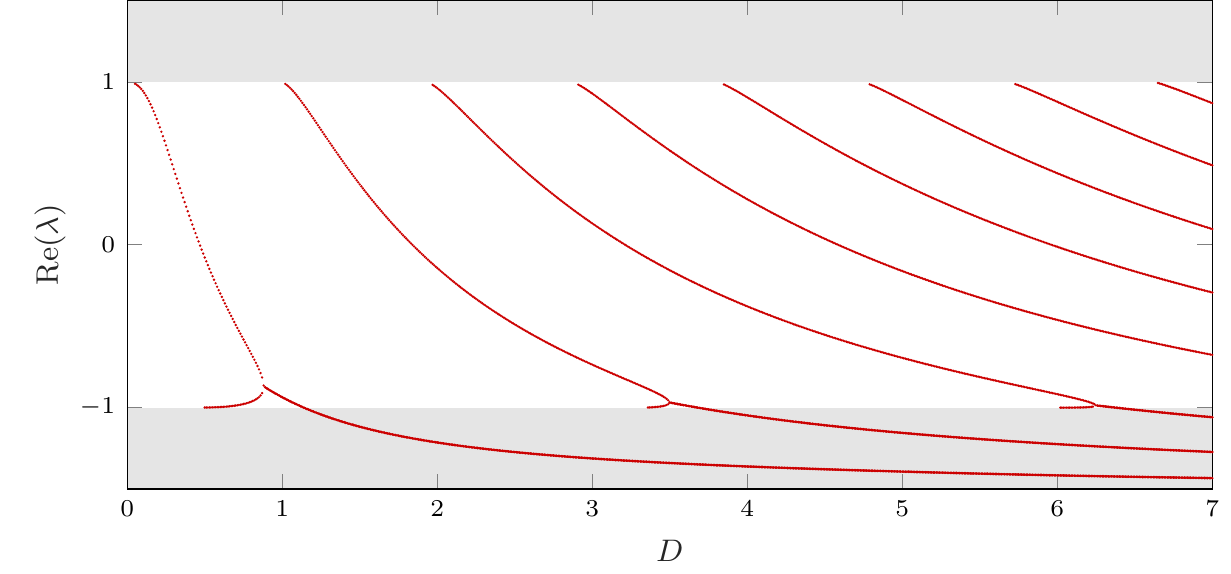}
	\vspace{-4mm}
	\caption{Real part of $\sigma_p(T_V)$ plotted against $D$ for $v_0=2.5$, $W=0.1$, \vspace{-1mm}$D\in[0,6]$.}
	\label{fig:jiang}
\end{figure}

Figure \ref{fig:jiang} shows the real part of the spectrum plotted against the parameter $D$ which controls the width of the potential well. The region shaded in gray indicates the essential spectrum. The figure shows the lowest eigenvalue moving from $+1$ towards $-1$ as $D$ increases from 0. Around $D=0.5$ another eigenvalue appears at $-1$ and moves towards $+1$ until the two join up and form a complex conjugate pair around $D=0.9$. Note that from then on, Figure \ref{fig:jiang} only shows their matching real part, which moves below $-1$. The same fate is met by the next pairs of eigenvalues around $D=3.5$ and $D=6.25$. Our results and our corresponding error bounds thus confirm those obtained recently in \cite{jiang2021analysis} using the so-called \emph{split-operator technique} \cite{braun1999numerical}. 

\begin{figure}[htbp]
	\centering
	\includegraphics[width=\textwidth]{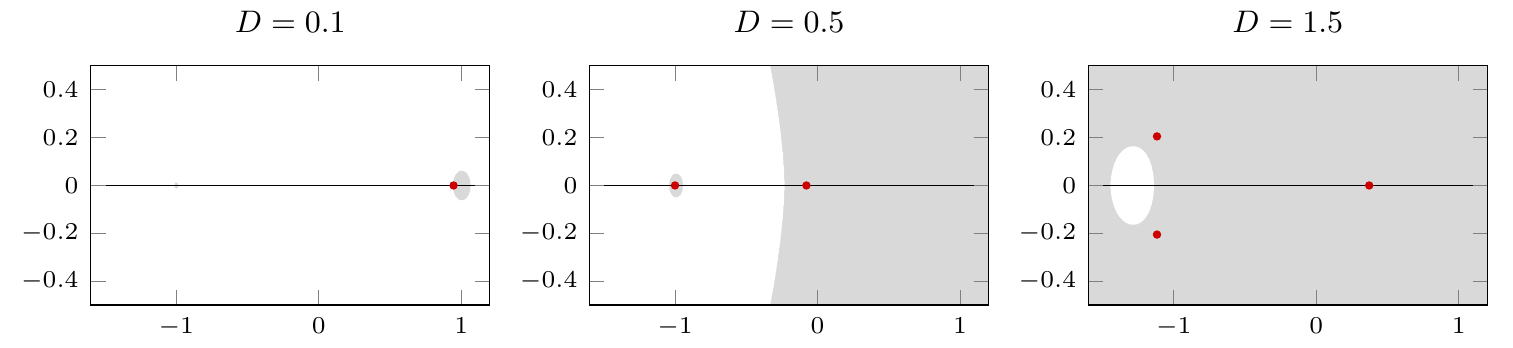}
	\vspace{-5mm}
	\caption{Approximations of the eigenvalues of $T_V$ for Sauter potential with $v_0\!=\!2.5$, $W\!=\!0.1$, $D\!\in\!\{0.1,0.5,1.5\}$ (red), together with region \eqref{eq:AD_bound1} \vspace{-3mm} (gray).}
	\label{fig:AD}
\end{figure}

Figure \ref{fig:AD} shows the eigenvalues in the complex plane for the three values $D\in\{0.1,0.5,1.5\}$, together with the inclusion region $\{\lambda\in\C \,|\,4|m^2-\lambda^2| \leq \|V^2-2\lambda V\|_{L^1(\R)}^2\text{ holds}\}$, see \eqref{eq:AD_bound1}; 
note that, especially for small values of $D$, this enclosure gives rather sharp bounds on the \vspace{-2mm} eigenvalues.

\subsubsection{Varying depth `cusp' potential well}
\begin{figure}[htbp]
	\centering
%
%
\setlength{\figurewidth}{0.7\textwidth}
\setlength{\figureheight}{.3\figurewidth}

\definecolor{mycolor1}{rgb}{0.00000,0.44700,0.74100}%
\definecolor{mycolor1}{rgb}{0.149039,0.508051,0.55725}%

\begin{tikzpicture}

\begin{axis}[%
width=\figurewidth,
height=\figureheight,
at={(0,0)},
scale only axis,
xmin=-4,
xmax=4,
xlabel={\small $x$},
ymin=-3.5,
ymax=0.5,
ylabel={\small $V$},
axis background/.style={fill=white}
]
\addplot [color=mycolor1, line width=1.0pt, forget plot]
  table[row sep=crcr]{%
-4	-0.00100638788370766\\
-3.44	-0.00308443213552412\\
-3.04	-0.00686452995876685\\
-2.8	-0.0110935911494492\\
-2.56	-0.0179280686850181\\
-2.4	-0.0246892411470601\\
-2.24	-0.0340002394640022\\
-2.08	-0.0468226737599489\\
-2	-0.0549469166662027\\
-1.92	-0.06448080403527\\
-1.84	-0.0756689245056821\\
-1.76	-0.0887983055036763\\
-1.68	-0.104205776834216\\
-1.6	-0.122286611935099\\
-1.52	-0.143504668482595\\
-1.44	-0.168404288502401\\
-1.36	-0.197624263279208\\
-1.28	-0.231914221329899\\
-1.2	-0.272153859868237\\
-1.12	-0.319375513137758\\
-1.04	-0.374790636595748\\
-0.96	-0.43982088639105\\
-0.88	-0.516134591469152\\
-0.8	-0.605689553983966\\
-0.72	-0.710783276046365\\
-0.64	-0.834111901359583\\
-0.56	-0.978839383869119\\
-0.48	-1.14867865792534\\
-0.4	-1.34798689235167\\
-0.32	-1.58187727212915\\
-0.24	-1.85635017541842\\
-0.16	-2.17844711122107\\
-0.0800000000000001	-2.55643136689863\\
0	-3\\
0.0800000000000001	-2.55643136689863\\
0.16	-2.17844711122107\\
0.24	-1.85635017541842\\
0.32	-1.58187727212915\\
0.4	-1.34798689235167\\
0.48	-1.14867865792534\\
0.56	-0.978839383869119\\
0.64	-0.834111901359583\\
0.72	-0.710783276046365\\
0.8	-0.605689553983966\\
0.88	-0.516134591469152\\
0.96	-0.43982088639105\\
1.04	-0.374790636595748\\
1.12	-0.319375513137758\\
1.2	-0.272153859868237\\
1.28	-0.231914221329899\\
1.36	-0.197624263279208\\
1.44	-0.168404288502401\\
1.52	-0.143504668482595\\
1.6	-0.122286611935099\\
1.68	-0.104205776834216\\
1.76	-0.0887983055036763\\
1.84	-0.0756689245056821\\
1.92	-0.06448080403527\\
2	-0.0549469166662027\\
2.08	-0.0468226737599489\\
2.24	-0.0340002394640022\\
2.4	-0.0246892411470601\\
2.56	-0.0179280686850181\\
2.8	-0.0110935911494492\\
3.04	-0.00686452995876685\\
3.36	-0.00361961464187388\\
3.84	-0.00138592469634524\\
4	-0.00100638788370766\\
};
\end{axis}

\end{tikzpicture}%
	\vspace{-4mm}
	\caption{`Cusp' potential \eqref{eq:cusp_potential} for $v_0=3$, \vspace{-1mm}$a=1$.}
	\label{fig:cusp_potential}
\end{figure}
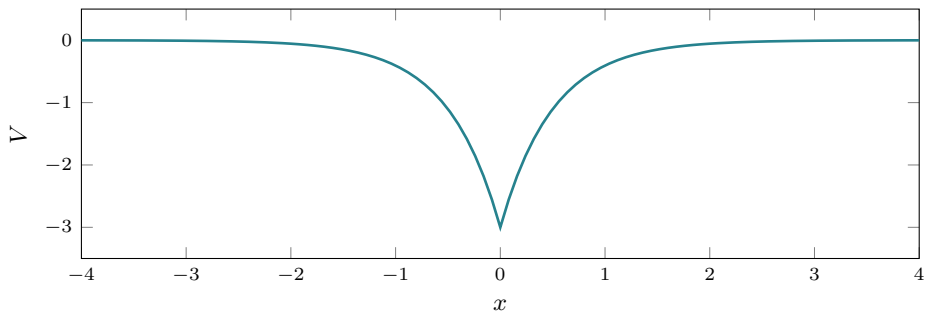

The so-called \emph{cusp potential} is defined as 
\begin{align}\label{eq:cusp_potential}
	V(x) = -v_0\cdot\begin{cases}
		\e^{x/a} & \text{for }x\leq 0,
		\\
		\e^{-x/a} & \text{for }x> 0,
	\end{cases}
\end{align}
where $v_0,a>0$ are parameters (see Figure \ref{fig:cusp_potential}).
Note that, strictly speaking, the point 0 is not a cusp mathematically, but a corner, i.e.\ the directional derivatives from the left and the right are finite at 0. We use the name \emph{cusp potential}, nonetheless, because it has become common in the literature.
It is easily seen that $V\in W^{1,\infty}(\R)$. Since $V$ decays exponentially, it satisfies Hypothesis \ref{hyp:assumptions_on_V} for any $v_0,a>0$.

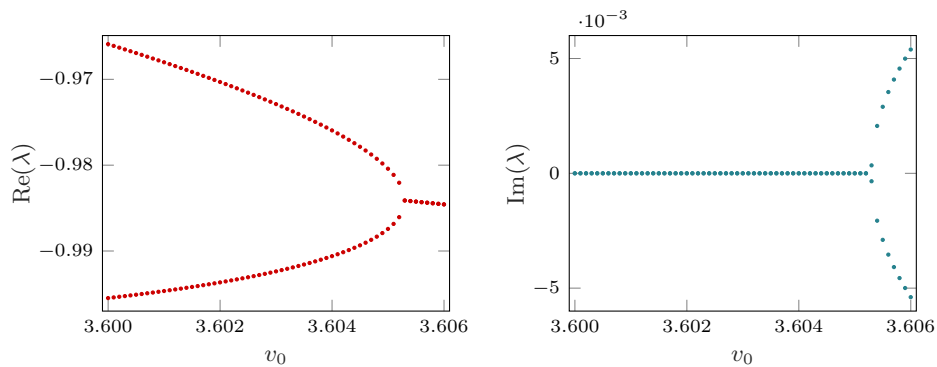
\begin{figure}[htbp]
	\centering
%
%
%
\definecolor{mycolor1}{rgb}{0.00000,0.44700,0.74100}%
\definecolor{mycolor2}{rgb}{0.85000,0.32500,0.09800}%
\definecolor{mycolor3}{rgb}{0.92900,0.69400,0.12500}%
\definecolor{mycolor4}{rgb}{0.49400,0.18400,0.55600}%
\definecolor{mycolor5}{rgb}{0.46600,0.67400,0.18800}%
\definecolor{mycolor6}{rgb}{0.30100,0.74500,0.93300}%
\definecolor{mycolor7}{rgb}{0.63500,0.07800,0.18400}%
\definecolor{mycolor8}{rgb}{0.149039,0.508051,0.55725}%

\setlength{\figurewidth}{0.29\textwidth}
\setlength{\figureheight}{0.8\figurewidth}
\begin{tikzpicture}

\begin{axis}[%
width=\figurewidth,
height=\figureheight,
at={(0,0)},
scale only axis,
xmin=3.5999,
xmax=3.6061,
xlabel style={font=\color{white!15!black}},
xlabel={\small $v_0$},
ymin=-0.997,
ymax=-0.9649,
ylabel style={font=\color{white!15!black}},
x tick label style={/pgf/number format/zerofill,
					/pgf/number format/precision=3},
ylabel={\small Re$(\lambda)$},
ylabel shift = -2mm,
axis background/.style={fill=white}
]
\addplot [color=mycolor1, draw=none, mark size=0.6pt, mark=*, mark options={solid, black!20!red}, forget plot]
  table[row sep=crcr]{%
  3.6	-0.9658985578\\
3.6001	-0.9661032477\\
3.6002	-0.9663092962\\
3.6003	-0.9665167464\\
3.6004	-0.9667256397\\
3.6005	-0.9669360213\\
3.6006	-0.9671479354\\
3.6007	-0.9673614335\\
3.6008	-0.967576567\\
3.6009	-0.9677933899\\
3.601	-0.9680119613\\
3.6011	-0.9682323416\\
3.6012	-0.9684545958\\
3.6013	-0.9686787934\\
3.6014	-0.9689050073\\
3.6015	-0.9691333154\\
3.6016	-0.9693638014\\
3.6017	-0.9695965535\\
3.6018	-0.9698316675\\
3.6019	-0.9700692441\\
3.602	-0.9703093948\\
3.6021	-0.9705522364\\
3.6022	-0.9707978951\\
3.6023	-0.9710465098\\
3.6024	-0.9712982293\\
3.6025	-0.9715532148\\
3.6026	-0.9718116435\\
3.6027	-0.9720737095\\
3.6028	-0.9723396232\\
3.6029	-0.9726096198\\
3.603	-0.9728839554\\
3.6031	-0.9731629183\\
3.6032	-0.9734468263\\
3.6033	-0.9737360382\\
3.6034	-0.9740309582\\
3.6035	-0.974332043\\
3.6036	-0.9746398141\\
3.6037	-0.974954873\\
3.6038	-0.975277914\\
3.6039	-0.9756097525\\
3.604	-0.9759513501\\
3.6041	-0.9763038602\\
3.6042	-0.9766686807\\
3.6043	-0.9770475388\\
3.6044	-0.9774426058\\
3.6045	-0.9778566806\\
3.6046	-0.97829347\\
3.6047	-0.9787580771\\
3.6048	-0.9792578589\\
3.6049	-0.9798041167\\
3.605	-0.9804159683\\
3.6051	-0.9811308692\\
3.6052	-0.9820467229\\
3.6053	-0.9841157823\\
3.6054	-0.9841804013\\
3.6055	-0.9842450231\\
3.6056	-0.9843096445\\
3.6057	-0.9843742644\\
3.6058	-0.984438887\\
3.6059	-0.9845035082\\
3.606	-0.9845681277\\
3.6	-0.9954828428\\
3.6001	-0.9954074175\\
3.6002	-0.99533063\\
3.6003	-0.9952524411\\
3.6004	-0.995172807\\
3.6005	-0.9950916856\\
3.6006	-0.9950090286\\
3.6007	-0.9949247874\\
3.6008	-0.9948389101\\
3.6009	-0.9947513425\\
3.601	-0.9946620251\\
3.6011	-0.9945708984\\
3.6012	-0.9944778984\\
3.6013	-0.9943829519\\
3.6014	-0.9942859896\\
3.6015	-0.9941869323\\
3.6016	-0.994085697\\
3.6017	-0.9939821947\\
3.6018	-0.9938763301\\
3.6019	-0.9937680005\\
3.602	-0.9936570998\\
3.6021	-0.9935435067\\
3.6022	-0.9934270951\\
3.6023	-0.9933077278\\
3.6024	-0.9931852555\\
3.6025	-0.9930595162\\
3.6026	-0.992930334\\
3.6027	-0.9927975136\\
3.6028	-0.9926608455\\
3.6029	-0.9925200948\\
3.603	-0.9923750043\\
3.6031	-0.9922252869\\
3.6032	-0.992070623\\
3.6033	-0.9919106555\\
3.6034	-0.9917449803\\
3.6035	-0.9915731399\\
3.6036	-0.9913946123\\
3.6037	-0.9912088001\\
3.6038	-0.9910150006\\
3.6039	-0.9908124061\\
3.604	-0.990600052\\
3.6041	-0.9903767855\\
3.6042	-0.9901412086\\
3.6043	-0.989891594\\
3.6044	-0.9896257692\\
3.6045	-0.9893409386\\
3.6046	-0.9890333917\\
3.6047	-0.9886980276\\
3.6048	-0.988327492\\
3.6049	-0.9879104739\\
3.605	-0.9874278678\\
3.6051	-0.9868422059\\
3.6052	-0.9860555947\\
3.6053	-0.9841157823\\
3.6054	-0.984180402\\
3.6055	-0.984245027\\
3.6056	-0.9843096446\\
3.6057	-0.9843742664\\
3.6058	-0.9844388872\\
3.6059	-0.9845035082\\
3.606	-0.9845681288\\
};
\end{axis}

\begin{axis}[%
width=\figurewidth,
height=\figureheight,
at={(\figurewidth+0.1\textwidth,0)},
scale only axis,
xmin=3.5999,
xmax=3.6061,
xlabel style={font=\color{white!15!black}},
x tick label style={/pgf/number format/zerofill,
/pgf/number format/precision=3},
xlabel={\small $v_0$},
ymin=-0.006,
ymax=0.006,
ylabel style={font=\color{white!15!black}},
ylabel={\small Im$(\lambda)$},
ylabel shift = -2mm,
]
\addplot [color=mycolor1, draw=none, mark size=0.6pt, mark=*, mark options={solid, mycolor8}, forget plot]
  table[row sep=crcr]{%
3.6	1.00000008274037e-10\\
3.6001	-1.10000009101441e-09\\
3.6002	5.00000041370185e-10\\
3.6003	-2.00000016548074e-10\\
3.6004	-1.00000008274037e-10\\
3.6005	-1.00000008274037e-09\\
3.6006	-2.00000016548074e-10\\
3.6007	-2.00000016548074e-10\\
3.6008	-1.00000008274037e-10\\
3.6009	3.00000024822111e-10\\
3.601	-2.00000016548074e-10\\
3.6011	2.00000016548074e-10\\
3.6012	-2.00000016548074e-10\\
3.6013	-5.00000041370185e-10\\
3.6014	2.00000016548074e-10\\
3.6015	-2.00000016548074e-10\\
3.6016	-0\\
3.6017	-4.00000033096148e-10\\
3.6018	0\\
3.6019	6.00000049644223e-10\\
3.602	-1.00000008274037e-10\\
3.6021	-3.00000024822111e-10\\
3.6022	-0\\
3.6023	-2.00000016548074e-10\\
3.6024	-2.00000016548074e-10\\
3.6025	-1.00000008274037e-10\\
3.6026	4.00000033096148e-10\\
3.6027	-8.00000066192297e-10\\
3.6028	5.00000041370185e-10\\
3.6029	2.00000016548074e-10\\
3.603	-7.0000005791826e-10\\
3.6031	-1.60000013238459e-09\\
3.6032	-2.00000016548074e-10\\
3.6033	-0\\
3.6034	2.00000016548074e-10\\
3.6035	-3.00000024822111e-10\\
3.6036	1.00000008274037e-10\\
3.6037	2.00000016548074e-10\\
3.6038	2.00000016548074e-10\\
3.6039	1.30000010756248e-09\\
3.604	-1.30000010756248e-09\\
3.6041	3.00000024822111e-10\\
3.6042	5.00000041370185e-10\\
3.6043	-3.00000024822111e-10\\
3.6044	1.00000008274037e-09\\
3.6045	-4.09999989514631e-09\\
3.6046	1.80000014893267e-09\\
3.6047	-4.00000033096148e-10\\
3.6048	-1.10000009101441e-09\\
3.6049	7.0000005791826e-10\\
3.605	1.00000008274037e-10\\
3.6051	7.0000005791826e-10\\
3.6052	-9.00000074466334e-10\\
3.6053	0.000345794899999863\\
3.6054	-0.00206328910000009\\
3.6055	0.0028974129999999\\
3.6056	0.0035402105000002\\
3.6057	0.00408306729999985\\
3.6058	0.00456180019999985\\
3.6059	0.00499488519999991\\
3.606	-0.00539332640000012\\
3.6	-1.00000008274037e-10\\
3.6001	-0\\
3.6002	1.00000008274037e-10\\
3.6003	-1.00000008274037e-10\\
3.6004	8.00000066192297e-10\\
3.6005	1.00000008274037e-10\\
3.6006	1.00000008274037e-10\\
3.6007	1.00000008274037e-10\\
3.6008	-1.00000008274037e-10\\
3.6009	-1.10000009101441e-09\\
3.601	-1.00000008274037e-10\\
3.6011	-7.0000005791826e-10\\
3.6012	-2.20000018202882e-09\\
3.6013	1.00000008274037e-10\\
3.6014	1.00000008274037e-10\\
3.6015	-1.00000008274037e-09\\
3.6016	-1.00000008274037e-10\\
3.6017	2.00000016548074e-10\\
3.6018	-1.00000008274037e-10\\
3.6019	-1.40000011583652e-09\\
3.602	-2.00000016548074e-10\\
3.6021	4.00000033096148e-10\\
3.6022	1.00000008274037e-10\\
3.6023	-2.00000016548074e-10\\
3.6024	0\\
3.6025	-1.00000008274037e-10\\
3.6026	-7.0000005791826e-10\\
3.6027	1.00000008274037e-10\\
3.6028	-7.0000005791826e-10\\
3.6029	1.00000008274037e-10\\
3.603	-3.00000024822111e-10\\
3.6031	-2.00000016548074e-10\\
3.6032	-1.00000008274037e-10\\
3.6033	3.00000024822111e-10\\
3.6034	3.00000024822111e-10\\
3.6035	2.00000016548074e-10\\
3.6036	-2.10000017375478e-09\\
3.6037	-2.00000016548074e-10\\
3.6038	3.00000024822111e-10\\
3.6039	-2.00000016548074e-10\\
3.604	-2.00000016548074e-10\\
3.6041	2.00000016548074e-10\\
3.6042	1.80000014893267e-09\\
3.6043	0\\
3.6044	-2.00000016548074e-10\\
3.6045	3.00000024822111e-10\\
3.6046	-0\\
3.6047	5.00000041370185e-10\\
3.6048	3.00000024822111e-10\\
3.6049	1.80000014893267e-09\\
3.605	-9.00000074466334e-10\\
3.6051	-1.00000008274037e-10\\
3.6052	-1.50000012411056e-09\\
3.6053	-0.000345799200000219\\
3.6054	0.00206328940000011\\
3.6055	-0.0028974170999998\\
3.6056	-0.0035402105000002\\
3.6057	-0.0040830677999999\\
3.6058	-0.00456180049999988\\
3.6059	-0.00499488529999992\\
3.606	0.00539332600000009\\
};
\end{axis}

\end{tikzpicture}%
	\vspace{-2mm}
	\caption{Real part (left) and imaginary part (right) of $\sigma(T_V)$ plotted against $v_0$ for a cusp potential well with $a=0.5$ and \vspace{-1mm} $v_0\in[3.600,3.606]$.}
	\label{fig:cusp_spectrum}
\end{figure}

Figure \ref{fig:cusp_spectrum} shows the collision of two initially real eigenvalues as the depth $v_0$ increases from $3.6$ to $3.606$. Similarly to the Sauter potential where the \emph{width} was varied, the plot shows one eigenvalue moving left towards $-1$ and another one moving right towards $+1$ until they collide and form a complex conjugate pair around $v_0=3.6054$. Our computations and accompanying error bounds confirm a result of \cite{villalba2006bound}, which has been obtained previously using a different numerical method (cf.\ \cite[Fig.\ 7]{villalba2006bound}).

\subsection{Convergence and error analysis}
\begin{figure}[htbp]
	\centering
%
%
\begin{tikzpicture}

\setlength{\figurewidth}{0.33\textwidth}
\setlength{\figureheight}{.9\figurewidth}
\definecolor{mycolor1}{rgb}{0.19121,0.40577,0.55600}%
\definecolor{mycolor2}{rgb}{0.21002,0.71984,0.47178}%
\definecolor{mycolor3}{rgb}{0.26819,0.22276,0.51160}%
\definecolor{mycolor4}{rgb}{0.12815,0.56511,0.55089}%
\definecolor{mycolor5}{rgb}{0.56718,0.84279,0.26180}%

\begin{axis}[%
width=\figurewidth,
height=\figureheight,
at={(1.37\figurewidth,0)},
scale only axis,
xmode=log,
xmin=95,
xmax=1700,
xminorticks=true,
xlabel={$n$},
ymode=log,
ymin=4e-06,
ymax=2e-3,
yminorticks=true,
ylabel={\footnotesize $|E^{\pm}_i(2n)-E^{\pm}_i(n)|$},
axis on top,
legend style={legend cell align=left, align=left, draw=white!15!black, only marks, legend columns=2}
]
\addplot [color=mycolor1, line width=0.5pt, mark size=2.0pt, mark=star, mark options={solid, mycolor1}]
  table[row sep=crcr]{%
100	9.87523421433383e-05\\
200	4.55674276997255e-05\\
400	2.11064934214028e-05\\
800	1.03070749978237e-05\\
1600	4.98345515847857e-06\\
};
\addlegendentry{\footnotesize$E^{\pm}_1$\;}

\addplot [color=mycolor2, line width=0.5pt, mark=o, mark size=1.5pt, mark options={solid, mycolor2}]
  table[row sep=crcr]{%
100	0.00161078430723746\\
200	0.000657560725807388\\
400	0.000289895393289275\\
800	0.000111715179210853\\
1600	4.93573019413889e-05\\
};
\addlegendentry{\footnotesize$E_2$}

\end{axis}

\begin{axis}[%
width=\figurewidth,
height=\figureheight,
at={(0,0)},
scale only axis,
xmode=log,
xmin=95,
xmax=1700,
xminorticks=true,
xlabel={$n$},
ymode=log,
ymin=2e-05,
ymax=5e-3,
yminorticks=true,
ylabel={\footnotesize $|E^{\pm}_i(2n)-E^{\pm}_i(n)|$},
axis on top,
legend style={legend cell align=left, align=left, draw=white!15!black,legend columns=3, only marks}
]
\addplot [color=mycolor3, line width=0.5pt, mark size=2.0pt, mark=star, mark options={solid, mycolor3}]
  table[row sep=crcr]{%
100	0.00441937714146339\\
200	0.000357447261461369\\
400	0.000120164262122407\\
800	9.15986985699998e-05\\
1600	4.52602937474415e-05\\
};
\addlegendentry{\footnotesize$E^{\pm}_1$\;}

\addplot [color=mycolor1, line width=0.5pt, mark=o, mark size=1.5pt, mark options={solid, mycolor1}]
  table[row sep=crcr]{%
100	0.000203184023508899\\
200	0.000256536856848163\\
400	0.000130514209321665\\
800	6.53144745578389e-05\\
1600	3.26526026450162e-05\\
};
\addlegendentry{\footnotesize$E^{\pm}_2$\;}

\addplot [color=mycolor4, line width=0.5pt, mark=asterisk, mark size=2pt, mark options={solid, mycolor4}]
  table[row sep=crcr]{%
100	0.00114328755776145\\
200	0.000255986169633974\\
400	0.000115520522930607\\
800	5.74812486907006e-05\\
1600	2.87303479197721e-05\\
};
\addlegendentry{\footnotesize$E^{\pm}_3$}

\addplot [color=mycolor2, line width=0.5pt, mark=square, mark size=1.2pt, mark options={solid, mycolor2}]
  table[row sep=crcr]{%
100	0.00178890236584606\\
200	0.000640156785603719\\
400	0.000296247366336671\\
800	0.000147698263795765\\
1600	7.38531090152817e-05\\
};
\addlegendentry{\footnotesize$E_4$\;}

\addplot [color=mycolor5, line width=0.5pt, mark=diamond, mark size=1.6pt, mark options={solid, mycolor5}]
  table[row sep=crcr]{%
100	0.00412832385124\\
200	0.000705213355087663\\
400	0.0003265571672062\\
800	0.000162910673662031\\
1600	8.14462479239442e-05\\
};
\addlegendentry{\footnotesize$E_5$}

\end{axis}

\end{tikzpicture}%
	\vspace{-2mm}
	\caption{Left: Cauchy errors $e_{2n,n}$ for the Sauter potential \eqref{eq:sauter_potential} with $v_0\!=\!3.7$, $D=3.2$, $W\!=\!0.3$ for $n\in\{100,200,400,800,1600\}$. Right: Cauchy errors $e_{2n,n}$ for a square well potential $V\!=\!-2.11\chi_{[-1,1]}$. Note that each pair $E^{\pm}_k$ of complex conjugate eigenvalues has identical values for $e_{m,n}$, hence only five distinct lines (resp.\ two distinct lines) remain in the figure on the left 
\vspace{-1mm}	(resp.\ right).}
	\label{fig:conv_rates}
\end{figure}
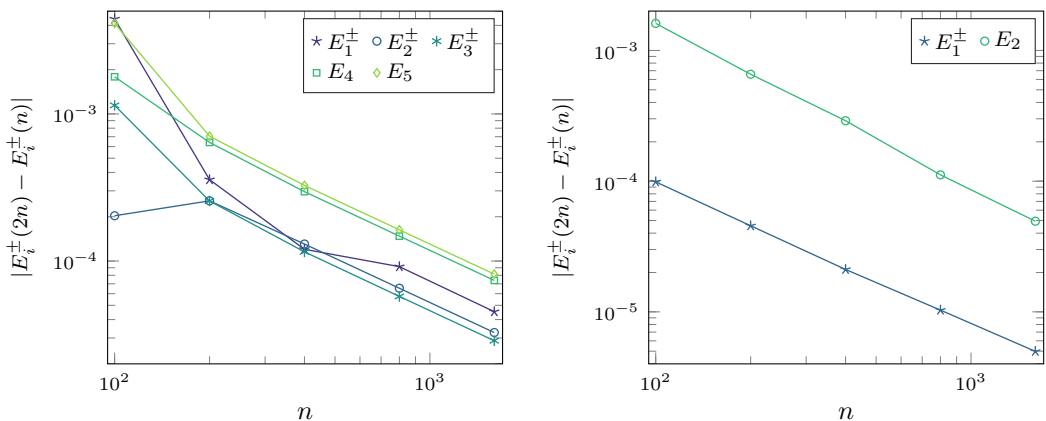

In order to estimate the convergence rate and the approximation error of our algorithm, we study the behaviour of its output as $n\to\infty$. To this end, we compute the approximate eigenvalues $E_k^{(\pm)}(n)$ for different values of $n$ and consider the relative Cauchy errors
\begin{align*}
	e_{mn} := |E_k^{(\pm)}(m) - E_k^{(\pm)}(n)|, \quad m,n\in\N.
\end{align*}
We use the rate at which these differences converge to 0 as a measure for the performance of $\Gamma_n$. The left-hand plot in Figure \ref{fig:conv_rates} shows the evolution of $e_{mn}$ for $m=2n$ and $n\in\{100,200,400,800,1600\}$ for the 8 eigenvalues $E_k^{(\pm)}$ of the Sauter potential (cf.\ Figure \ref{fig:sauter_contour}).

Figure \ref{fig:conv_rates} suggests that all 8 eigenvalues eventually converge with the same rate. The slope of the lines in the left hand plot of Figure \ref{fig:conv_rates} implies a rate of $e_{2n,n}\lesssim n^{-1}$, which suggests that the abstract bound on the convergence rate obtained from \eqref{eq:improved_params_1d} is rather conservative and the algorithm performs better in practice.

\begin{figure}[htbp]
	\centering
	\begin{tikzpicture}
\definecolor{mycolor1}{rgb}{0.00000,0.44700,0.74100}%
\definecolor{mycolor2}{rgb}{0.85000,0.32500,0.09800}%
\definecolor{mycolor3}{rgb}{0.92900,0.69400,0.12500}%
\definecolor{vir}{rgb}{0.149039,0.508051,0.55725}%

\pgfplotsset{scaled y ticks=false}
\setlength{\figurewidth}{0.35\textwidth}
\setlength{\figureheight}{0.6\figurewidth}

\begin{axis}[%
width=\figurewidth,
height=\figureheight,
at={(0,0)},
scale only axis,
xmin=-2,
xmax=2,
xlabel={$x$},
ymin=-2.5,
ymax=0.4,
ylabel={$V$},
]
\addplot [color=vir, line width=1.0pt, forget plot]
  table[row sep=crcr]{%
-2	-0\\
-1.000300030003	-0\\
-0.999499949995	-2.11\\
0.999899989999	-2.11\\
1.000700070007	-0\\
2	-0\\
};
\end{axis}

\begin{axis}[%
width=\figurewidth,
height=\figureheight,
at={(1.4\figurewidth,0)},
scale only axis,
xmin=-1.5,
xmax=1.5,
xlabel={Re$(z)$},
yticklabels = {-0.15,-0.1,-0.05,0,0.05,0.1},
ymin=-0.1,
ymax=0.1,
ylabel={Im$(z)$},
ylabel shift = -0.1cm,
]
\addplot [color=mycolor1, draw=none, mark size=1.6pt, mark=*, mark options={solid, black!20!red}, forget plot]
  table[row sep=crcr]{%
-0.963084546339517  -0.066464998217296\\
-0.963084545650382  +0.066464998009696\\
0.148557009776029  -0.000000000001477\\
};
\end{axis}

\end{tikzpicture}%
	\vspace{-2mm}
	\caption{Left:\! square well potential of width 2 and depth $-2.11$. Right:\! approximate eigen\-values 
	 $E_1^\pm \in\C\setminus \R$ and $E_2\in\R$ returned by $\Gamma_n$ for \vspace{-3mm}$n=3200$.}
	\label{fig:squarewell}
\end{figure}
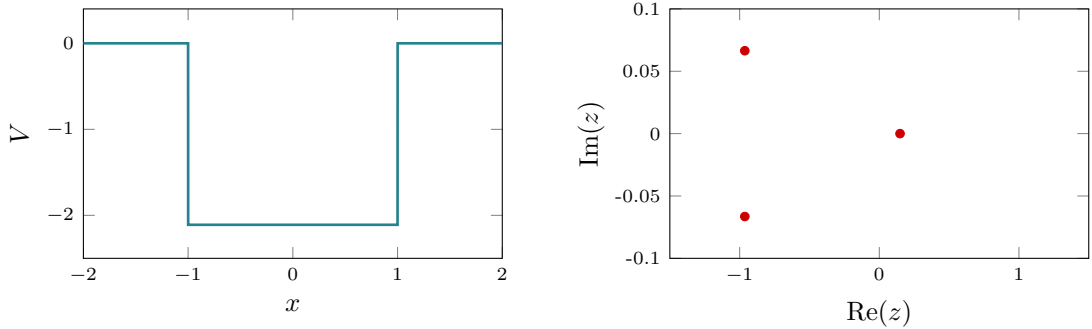

To estimate the approximation error of $\Gamma_{\!n}$, we apply our algorithm to a square well potential~for which the eigenvalues can be computed explicitly using an exponential ansatz. A lengthy but straight\-forward calculation shows that $\lambda\!\in\!\C$ is an eigenvalue of $T_V$ with potential $V\!=\!v_0\chi_{[-a,a]}$ if and only~if
\begin{align}\label{eq:squarewell_exact_eq}
	2\cot\left(2a\sqrt{(v_0-\lambda)^2-m^2}\right) = \frac{\sqrt{(v_0-\lambda)^2-m^2}}{\sqrt{m^2-\lambda^2}} - \frac{\sqrt{m^2-\lambda^2}}{\sqrt{(v_0-\lambda)^2-m^2}}.
\end{align}
Equation \eqref{eq:squarewell_exact_eq} can be solved with very high accuracy (we used 64 digits after the decimal point) using Newton's method. The numerical values of its solutions for $a=1$, $v_0=-2.11$ (these parameter values were previously considered in \cite{MR2268872}) are given in Table \ref{table:squarewell}.
Figure \ref{fig:squarewell} shows the potential and the output of $\Gamma_n$. As the figure suggests, our algorithm returns three eigenvalues $E_1^\pm \in\C\setminus \R$ and $E_2\in\R$ whose numerical values are given in Table \ref{table:squarewell}.
\begin{table}[H]
	\begin{tabular}{c r r}\toprule
		\textbf{Eigenvalue}\; & \textbf{Approximate value}\; & \textbf{Exact value (truncated)} \\
		\hline\\[-3mm]
		$E_1^\pm$ & $-0.9630845 \pm 0.0664649\I$ & $-0.9630888 \pm 0.0664672\I$\\
		$E_2$ & $0.1485570 + 0.0000000\I$ & $0.1485171 + 0.0000000\I$\\
		\bottomrule\\
	\end{tabular}
	\caption{Approximate and exact numerical values of the eigenvalues $E_1^\pm$, $E_2$ of $T_V$ for a square well potential of width 2 and depth $-2.11$. Approximate values are computed using $\Gamma_n$ with $n=3200$; exact values are computed from eq. \eqref{eq:squarewell_exact_eq} via Newton iterations
to 64 digits after the decimal \vspace{-9mm} point.}
	\label{table:squarewell}
\end{table}

A comparison between the approximate and exact values in Table \ref{table:squarewell} shows that our values for $E_1^{\pm}$ agree with the exact ones to within 5 digits after the decimal point, while those for $E_2$ agree to within 4 digits. This suggests approximation errors of $O(10^{-6})$ and $O(10^{-5})$, respectively. The right hand side of Figure \ref{fig:conv_rates} shows the corresponding Cauchy errors $e_{2n,n}$ for $n$ up to 1600. The slopes of the lines for $E_1^\pm$ and $E_2$ again imply a convergence rate of $e_{2n,n}\lesssim n^{-1}$ (as for the Sauter potential). Moreover, the final Cauchy errors for $E_1^{\pm}$ resp.\ $E_2$ are given by 
\begin{align*}
	|E_1^{\pm}(3200) - E_1^{\pm}(1600)| &= e_{3200,1600} \approx 5\cdot10^{-6},
	\\
	|E_2(3200) - E_2(1600)| &= e_{3200,1600} \approx 5\cdot10^{-5},
\end{align*}
which are in agreement with the order of magnitude expected from Table \ref{table:squarewell} and reflect the fact that the approximation error of $E_1^{\pm}$ is one order of magnitude less than that of $E_2$.

\bibliography{mybib}
\bibliographystyle{abbrv}

\end{document}